\documentclass{etds}  
\usepackage{amsfonts} 
\usepackage{amssymb} 
\usepackage{latexsym} 
\usepackage{graphicx}   
\usepackage{amsbsy} 
\usepackage{psfrag}

\newcommand{\D}{\displaystyle}
\newcommand{\restr}[1]{\raisebox{-0.3em}{$\lb|_{#1}\rb.$}} 
\newcommand{\breath}{\medskip} 
\newtheorem{thm}{Theorem}[section]

\newcounter{claimcount}[thm]  
\newcounter{subclaimcount}[claimcount]
\newtheorem{prop}[thm]{Proposition} 

\newtheorem{lemma}[thm]{Lemma} 

\newtheorem{cor}[thm]{Corollary}

\newcommand{\Theorem}[2]{\begin{thm}{\sf #1}  #2 \end{thm}}
\newcommand{\Proposition}[2]{\begin{prop}{\sf #1}  #2 \end{prop}}
\newcommand{\Lemma}[2]{\begin{lemma}{\sf #1}  #2 \end{lemma}}
\newcommand{\Corollary}[2]{\begin{cor}{\sf #1}  #2 \end{cor}} 
\newcommand{\thmfont}[1]{{\sl #1}}    
\newcommand{\example}[1]{        \refstepcounter{thm}                     \begin{list}{} 			{\setlength{\leftmargin}{0em} 			\setlength{\rightmargin}{0em}}        \item {\em Example \thethm.} #1                   \hfill$\diamondsuit$  \end{list}   			}   
\newcommand{\bthmlist}{ \begin{list}{{\bf (\alph{enumii})}}{\usecounter{enumii} 			\setlength{\leftmargin}{0.5em} 			\setlength{\itemsep}{0em} 			\setlength{\parsep}{0em} 			\setlength{\rightmargin}{0em}}}
\newcommand{\ethmlist}{\end{list}}    
\newcommand{\Claim}[1]{\refstepcounter{claimcount}                {\sc Claim \theclaimcount. \ }\thmfont{ #1}} 
\newcommand{\Subclaim}[1]{\refstepcounter{subclaimcount}                 {\sc Claim \theclaimcount.\thesubclaimcount. \ }\thmfont{ #1}} 
\newcommand{\claim}{\Claim}
\newcommand{\subclaim}{\Subclaim}
\newcommand{\bprf}[1][Proof.]{\begin{list}{} 			{\setlength{\leftmargin}{0.7em} 			\setlength{\rightmargin}{0em}}                         \item {\em \hspace{-1em}  #1 \ \ }} 
\newcommand{\eprf}{\end{list}} 
\newcommand{\bthmprf}{\bprf}
\newcommand{\bclaimprf}{\bprf}
\newcommand{\bsubclaimprf}{\bprf}
\newcommand{\ethmprf}{ {\tt \hfill $\Box$ } \eprf  \breath  } 
\newcommand{\eclaimprf}{ \hfill $\Diamond$~{\scriptsize {\tt Claim~\theclaimcount}}\eprf}  
\newcommand{\esubclaimprf}{ \hfill $\triangledown$~{\scriptsize  {\tt Claim~\theclaimcount.\thesubclaimcount}}\eprf}  
\newcommand{\eclaimthmprf}{ \hfill $\Diamond$~{\scriptsize {\tt Claim~\theclaimcount}}~$\Box$\eprf\eprf}  
\newcommand{\QED}{\hfill\ensuremath{\Box}}
\newcommand{\qed}{\QED}     
\newcommand{\beq}{\begin{eqnarray*}}
\newcommand{\eeq}{\end{eqnarray*}} 
\newcommand{\beqn}{ \begin{equation} }
\newcommand{\eeqn}{ \end{equation} }
\newcommand{\bitem}{\begin{itemize}}
\newcommand{\eitem}{\end{itemize}} 
\newcommand{\bdesc}{\begin{description}}
\newcommand{\edesc}{\end{description}}   
\newcommand{\If}{\mbox{\ if \ }} 
\newcommand{\AND}{\mbox{\ and \ }} 
\newcommand{\Cesaro}{Ces\`aro }
\newcommand{\done}{{\mathsf{ 1\!\!1}}} 
\newcommand{\dB}{{\mathbb{B}}}
\newcommand{\dC}{{\mathbb{C}}}
\newcommand{\dE}{{\mathbb{E}}}
\newcommand{\dI}{{\mathbb{I}}}
\newcommand{\dJ}{{\mathbb{J}}}
\newcommand{\dK}{{\mathbb{K}}}
\newcommand{\dL}{{\mathbb{L}}}
\newcommand{\dM}{{\mathbb{M}}}
\newcommand{\dN}{{\mathbb{N}}}
\newcommand{\dP}{{\mathbb{P}}}
\newcommand{\dQ}{{\mathbb{Q}}}
\newcommand{\dR}{{\mathbb{R}}}
\newcommand{\dT}{{\mathbb{T}}}
\newcommand{\dZ}{{\mathbb{Z}}}        
\newcommand{\bard}{{\overline{d}}}
\newcommand{\barlam }{{\overline{\lambda}}}
\newcommand{\barmu}{{\overline{\mu }}}
\newcommand{\barbD}{{\overline{\mathbf{ D}}}}
\newcommand{\barbO}{{\overline{\mathbf{ O}}}}
\newcommand{\barbP}{{\overline{\mathbf{ P}}}}
\newcommand{\barbX}{{\overline{\mathbf{ X}}}}
\newcommand{\barsP}{{\overline{\mathcal{ P}}}}
\newcommand{\barsQ}{{\overline{\mathcal{ Q}}}}
\newcommand{\bB}{{\mathbf{ B}}}
\newcommand{\bC}{{\mathbf{ C}}}
\newcommand{\bD}{{\mathbf{ D}}}
\newcommand{\bE}{{\mathbf{ E}}}
\newcommand{\bF}{{\mathbf{ F}}}
\newcommand{\bG}{{\mathbf{ G}}}
\newcommand{\bI}{{\mathbf{ I}}}
\newcommand{\bJ}{{\mathbf{ J}}}
\newcommand{\bK}{{\mathbf{ K}}}
\newcommand{\bL}{{\mathbf{ L}}}
\newcommand{\bM}{{\mathbf{ M}}}
\newcommand{\bO}{{\mathbf{ O}}}
\newcommand{\bP}{{\mathbf{ P}}}
\newcommand{\bQ}{{\mathbf{ Q}}}
\newcommand{\bR}{{\mathbf{ R}}}
\newcommand{\bS}{{\mathbf{ S}}}
\newcommand{\bT}{{\mathbf{ T}}}
\newcommand{\bU}{{\mathbf{ U}}}
\newcommand{\bV}{{\mathbf{ V}}}
\newcommand{\bW}{{\mathbf{ W}}}
\newcommand{\bX}{{\mathbf{ X}}}
\newcommand{\bY}{{\mathbf{ Y}}}
\newcommand{\bZ}{{\mathbf{ Z}}} 
\newcommand{\ba}{{\mathbf{ a}}}
\newcommand{\bb}{{\mathbf{ b}}}
\newcommand{\bc}{{\mathbf{ c}}}
\newcommand{\bi}{{\mathbf{ i}}}
\newcommand{\bo}{{\mathbf{ o}}}
\newcommand{\bp}{{\mathbf{ p}}}
\newcommand{\bq}{{\mathbf{ q}}}
\newcommand{\br}{{\mathbf{ r}}}
\newcommand{\bs}{{\mathbf{ s}}}
\newcommand{\bt}{{\mathbf{ t}}}
\newcommand{\bv}{{\mathbf{ v}}}
\newcommand{\bw}{{\mathbf{ w}}}
\newcommand{\bx}{{\mathbf{ x}}}
\newcommand{\sA}{{\mathcal{ A}}}
\newcommand{\sB}{{\mathcal{ B}}}
\newcommand{\sD}{{\mathcal{ D}}}
\newcommand{\sK}{{\mathcal{ K}}}
\newcommand{\sL}{{\mathcal{ L}}}
\newcommand{\sM}{{\mathcal{ M}}}
\newcommand{\sO}{{\mathcal{ O}}}
\newcommand{\sP}{{\mathcal{ P}}}
\newcommand{\sQ}{{\mathcal{ Q}}}
\newcommand{\sU}{{\mathcal{ U}}}
\newcommand{\gC}{{\mathfrak{ C}}}
\newcommand{\gF}{{\mathfrak{ F}}}
\newcommand{\gK}{{\mathfrak{ K}}}
\newcommand{\gM}{{\mathfrak{ M}}}
\newcommand{\gO}{{\mathfrak{ O}}}
\newcommand{\gP}{{\mathfrak{ P}}}
\newcommand{\gQ}{{\mathfrak{ Q}}}
\newcommand{\gR}{{\mathfrak{ R}}}
\newcommand{\gS}{{\mathfrak{ S}}}
\newcommand{\gX}{{\mathfrak{ X}}}
\newcommand{\alp }{\alpha}
\newcommand{\bet }{\beta}
\newcommand{\del }{\delta}
\newcommand{\eps }{\epsilon}
\newcommand{\kap }{\kappa}
\newcommand{\lam }{\lambda}
\newcommand{\sig }{\sigma} 
\newcommand{\Del }{\Delta}
\newcommand{\hbM}{{\widehat{\mathbf{ M}}}}
\newcommand{\fb}{{\mathsf{ b}}}
\newcommand{\fe}{{\mathsf{ e}}}
\newcommand{\fj}{{\mathsf{ j}}}
\newcommand{\fm}{{\mathsf{ m}}}
\newcommand{\fp}{{\mathsf{ p}}}
\newcommand{\fu}{{\mathsf{ u}}}
\newcommand{\fv}{{\mathsf{ v}}}
\newcommand{\tl}{\widetilde} 
\newcommand{\tlA}{{\widetilde{A}}}
\newcommand{\tlp}{{\widetilde{p}}}
\newcommand{\tlbK}{{\widetilde{\mathbf{ K}}}}
\newcommand{\tlbM}{{\widetilde{\mathbf{ M}}}}
\newcommand{\tlbP}{{\widetilde{\mathbf{ P}}}}
\newcommand{\tlbU}{{\widetilde{\mathbf{ U}}}}
\newcommand{\tlbW}{{\widetilde{\mathbf{ W}}}}
\newcommand{\tlbp}{{\widetilde{\mathbf{ p}}}}
\newcommand{\tlbq}{{\widetilde{\mathbf{ q}}}}
\newcommand{\tlsP}{{\widetilde{\mathcal{ P}}}}
\newcommand{\undbO}{{\underline{\mathbf{ O}}}}
\newcommand{\undbP}{{\underline{\mathbf{ P}}}}
\newcommand{\lb}{\left}
\newcommand{\rb}{\right} 
\newcommand{\maketall}{\rule[-0.5em]{0em}{1em}}       
\newcommand{\IMPLIES}{\ensuremath{\Longrightarrow}}
\newcommand{\map}{{\longrightarrow}}
\newcommand{\goto}{{\rightarrow}}
\newcommand{\into}{{\map}}
\newcommand{\seilpmi}{{\Longleftarrow}}
\newcommand{\statement}[1]{\lb(  \maketall       \begin{minipage}{40em}       \begin{tabbing}         #1        \end{tabbing}      \end{minipage}  \rb)}     
\newcommand{\oo}{{\infty}}        
\newcommand{\X}{\times}
\newcommand{\x}{\X}
\newcommand{\tensor}{\otimes}
\newcommand{\compl}[1]{#1^\complement}  
\newcommand{\symdif}{{\bigtriangleup}} 
\newcommand{\union}{\cup}
\newcommand{\Union}{\bigcup}
\newcommand{\intsct}{\cap}
\newcommand{\Intsct}{\bigcap}
\newcommand{\disj}{\sqcup}
\newcommand{\Disj}{\bigsqcup}   
\newcommand{\set}[2]{{\left\{ #1 \; ; \; #2 \right\} }} 
\newcommand{\supp}[1]{{\sf supp}\lb(#1\rb)}    
\newcommand{\BINOM}[2]{{  \left( \begin{array}{c} #1 \\ #2 \end{array} \right)}} 
\newcommand{\norm}[2]{{\left\| #1 \right\|_{{#2}} }   }
\newcommand{\inn}[1]{{\left\langle #1 \right\rangle }}       
\newcommand{\pr}[1]{{\mathbf{ pr}_{{#1}}}}
\newcommand{\chr}[1]{{{\done}_{{#1}}}} 
\newcommand{\choice}[1]{{\lb\{ \begin{array}{rcl}                                 #1                                \end{array}  \rb.  }}                              
\newcommand{\eeequals}[1]{\raisebox{-0.9ex}{$\overline{\overline{{\scriptscriptstyle{\mathrm{#1}}}}}$}} 
\newcommand{\leeeq}[1]{\raisebox{-1ex}{${{\D\leq} \atop {\scriptscriptstyle{\mathrm{#1}}}}$}} 
\newcommand{\lt}[1]{\raisebox{-1ex}{${{\D<} \atop {\scriptscriptstyle{\mathrm{#1}}}}$}} 
\newcommand{\grt}[1]{\raisebox{-1ex}{${{\D>} \atop {\scriptscriptstyle{\mathrm{#1}}}}$}}  
\newcommand{\geeeq}[1]{\raisebox{-1ex}{${{\geq} \atop {\scriptscriptstyle{\mathrm{#1}}}}$}} 
\newcommand{\iiimplies}[1]{\eeequals{\ #1}\!\!\!\Rightarrow}
\newcommand{\closeto}[1]{{{\raisebox{-1ex}    {$\widetilde{\ {\scriptstyle #1}\ }$}}}}      
\newcommand{\Fix}[1]{{\sf Fix}\lb[#1\rb]}
\newcommand{\shift}[1]{\sig^{#1}}    
\newcommand{\goesto}[2]{{ -\!\!\!-\!\!\!-\!\!\!-\!\!\!\!\!\!\!\!\!\!\!  ^{{\scriptscriptstyle #2}}_{{\scriptscriptstyle #1}}   \!\!\!\!\!\!\!\!\!\longrightarrow }}                         
\newcommand{\length}[1]{{\sf length}\lb[#1\rb]}
\newcommand{\mtrx}[3]{{\lb[#1  |_{#2}^{#3} \rb]}}         
\newcommand{\trace}[1]{{{\sf trace}\lb[#1\rb]}} 
\newcommand{\Real}{\dR}
\newcommand{\Natur}{\dN}
\newcommand{\Zahl}{\dZ}
\newcommand{\Zahlmod}[1]{{\Zahl_{/#1}}}
\newcommand{\Rat}{\dQ}
\newcommand{\Cplx}{\dC}
\newcommand{\Torus}[1]{{{\dT}^{#1}}} 
\newcommand{\CC}[1]{{\lb[ #1 \rb]}}
\newcommand{\CO}[1]{{\lb[ #1 \rb)}}
\newcommand{\OC}[1]{{\lb( #1 \rb]}}
\newcommand{\OO}[1]{{\lb( #1 \rb)}}   
\newcommand{\Selfmap}[1]{#1\into#1}      
\newcommand{\Lat}{\dL}    
\newcommand{\rot}[1]{\rho^{#1}}  
\newcommand{\tile}{\bw}
\newcommand{\til}[1][]{w^{#1}} 
\newcommand{\stile}{\bv}   
\newcommand{\muae}{{\rm ($\mu$-\ae)}}         
\newcommand{\Nh}{\dB}
\newcommand{\nh}{\fb}
\newcommand{\OP}{\,^o\!\sA^{\Tor}}
\newcommand{\SOTP}{\,^{\ddagger}\!\sA^{\Tor}}
\newcommand{\OTP}{\,^{\perp}\!\!\sA^{\Tor}}
\newcommand{\ZP}[1][0]{\,^o\!\sA^{\Tor}_{#1}}
\newcommand{\DP}{\,^{\frac{1}{2}}\!\sA^{\Tor}}
\newcommand{\MP}{\sA^{\Tor}} 
\newcommand{\ISP}{\,^*\!\sA^{\Tor}} 
\newcommand{\SP}{\,^s\!\sA^{\Tor}}
\newcommand{\TP}{\,^t\!\sA^{\Tor}}
\newcommand{\STP}{\,^{st}\!\sA^{\Tor}}
\newcommand{\PP}{\,^p\!\sA^{\Tor}}
\newcommand{\qP}{\,^{\torq}\!\sA^{\Tor}} 
\newcommand{\QM}[2][\trans{}]{\Upsilon_{#1}\lb(#2\rb)}
\newcommand{\QS}[1][\trans{}]{\gQ\gS_{#1}}
\newcommand{\QSM}{\sM^{\scriptscriptstyle{\mathrm{QS}}}_{\trans{}}\!\lb(\sA^\Lat\rb)}
\newcommand{\ergM}{\sM^{\scriptscriptstyle{\mathrm{erg}}}\!\lb(\sA^\Lat\rb)}  
\newcommand{\Tor}{\mathbf{T}}
\newcommand{\tlTor}{\widetilde{\Tor}}
\newcommand{\tora}{{\tt a}}
\newcommand{\torb}{{\tt b}}
\newcommand{\torc}{{\tt c}}
\newcommand{\torg}{{\tt g}}
\newcommand{\torm}{{\tt m}}
\newcommand{\torp}{{\tt p}}
\newcommand{\tort}{{\tt t}}
\newcommand{\tors}{{\tt s}}
\newcommand{\torr}{{\tt r}}
\newcommand{\toru}{{\tt u}}
\newcommand{\torq}{{\tt q}}
\newcommand{\torz}{{\tt z}}   
\newcommand{\density}[1]{{\sf density}\lb(#1\rb)}
\newcommand{\updensity}[1]{\overline{{\sf density}}\lb(#1\rb)}
\newcommand{\closure}[1]{{\sf closure}\lb(#1\rb)} 
\newcommand{\Boole}[1]{\varsigma^\Lat\lb(#1\rb)} 
\newcommand{\trans}[1]{\varsigma^{#1}} 
\newcommand{\QTor}{\Tor_\Rat}
\newcommand{\LQTor}{\Lat\!\QTor}  
\newcommand{\Qshift}[2][\trans{}]{\Xi_{#1}\lb(#2\rb)} 
\newcommand{\Qseq}[2][\trans{}]{\xi_{#1}\lb(#2\rb)} 
\newcommand{\Splat}{\sP_{\trans{}}}
\newcommand{\Sqlat}{\sQ_{\trans{}}}  
\newcommand{\sss}[1]{\paragraph*{\sc #1:}} 
\newcommand{\as}{\mathrm{a\!s}} 
\newcommand{\wkstlim}{\mathrm{wk}^*\!\!\!-\!\!\!\lim}
\newcommand{\symdiflim}{d_\symdif\!\!\!-\!\!\!\lim}  
\newcommand{\Phitor}{\Phi_{\trans{}}} 
\newcommand{\lamae}{{\rm ($\lam$-\ae)}} 
\newcommand{\topae}{{\rm (top.\ae)}}         
\newcommand{\setsize}[2][]{\lb\lceil #2 \rb\rceil_{#1}} 
\newcommand{\Symdif}{\lefteqn{\mbox{\Large$\bigtriangleup$}}\maketall}   
\newcommand{\Expct}[1]{\dE\lb(#1\rb)}
\newcommand{\given}{\ \lb.\maketall \rb| \ } 
\newcommand{\Sep}[1]{\lb\langle#1\rb\rangle}   
\newcommand{\quoTor}{\overline{\bT}}
\newcommand{\tlquoTor}{\widetilde{\quoTor}}
\newcommand{\quoOP}{\,^o\!\sA^{\quoTor}}
\newcommand{\quoMP}{\sA^{\quoTor}}
\newcommand{\quosP}{\barsP}
\newcommand{\quosQ}{\barsQ}
\newcommand{\quotrans}[1]{\bar{\varsigma}^{#1}}
\newcommand{\quotort}{\bar{\tort}}
\newcommand{\quoSplat}{\quosP_{\quotrans{}}}
\newcommand{\quoQM}[1]{\QM[\quotrans{}]{#1}}
\newcommand{\quoPhitor}{\Phi_{\quotrans{}}}
\newcommand{\quolam}{\barlam}
\newcommand{\quotau}{\overline{\tau}}   
\newcommand{\Ball}{\boldsymbol{\mathsf{B}}} 
\newcommand{\dBlim}{d_B\mbox{-}\!\lim}
\newcommand{\dClim}{d_C\mbox{-}\!\lim}
\newcommand{\symlim}{d_\symdif\!\!-\!\!\lim} 
\newcommand{\presub}[1]{\,_{_{#1}}\!} 
\newcommand{\Lucas}[1]{\sL\lb(#1\rb)}  
\newcommand{\tlgP}{\widetilde{\gP}}
\newcounter{XMPL}[thm]  
\newcommand{\Examples}[3][]{         \refstepcounter{thm} \hspace{-1.5em} {\em Example \thethm.} {\sf #2} #1 \begin{list}{$\langle$\alph{XMPL}$\rangle$}{\usecounter{XMPL}} 			{\setlength{\leftmargin}{0em} 			\setlength{\rightmargin}{0em} 			\setlength{\itemsep}{0em} 			\setlength{\parsep}{0em}}   #3 	\hfill$\diamondsuit$\end{list}   			}   

\renewcommand{\closure}[1]{{\sf closure}\lb(#1\rb)}

\newcommand{\dfn}{\sf\em}

\begin{document}

\ETDS{1}{53}{25}{2005}
\runningheads{M. Pivato}{Cellular Automata vs. Quasisturmian Shifts}

\title{Cellular Automata vs. Quasisturmian Shifts}

\author{Marcus Pivato}
\address{
  Department of Mathematics,  Trent University,
\\ Peterborough, Ontario, Canada, K9L 1Z8.
\\
\email{{\tt pivato@xaravve.trentu.ca}}}

\recd{?}
%\maketitle

\begin{abstract} If $\dL=\Zahl^D$ and $\sA$ is a finite set, then $\sA^\dL$ is a compact space; a {\dfn cellular automaton} (CA) is a continuous transformation $\Phi:\sA^\dL
\into \sA^\dL$ that commutes with all shift maps.  A {\dfn
quasisturmian} (QS) {\dfn subshift} is a shift-invariant subset
obtained by mapping the trajectories of an irrational torus rotation
through a partition of the torus. The image of a QS
shift under a CA is again QS.  We study the topological
dynamical properties of CA restricted to QS shifts, and compare them
to the properties of CA on the full shift $\sA^\Lat$.  We investigate
injectivity, surjectivity, transitivity, expansiveness, rigidity,
fixed/periodic points, and invariant measures.  We also study 
`chopping': how iterating the CA fragments the partition generating the 
QS shift.

{\footnotesize
\breath

\begin{tabular}{rl}
{\bf MSC:}& 37B15 (primary), 68Q80 (secondary)\\
{\bf Keywords:}& Cellular automata, quasiperiodic, Sturmian, 
discrete spectrum. 
\end{tabular}}
\end{abstract}

  Let $D\geq 1$, and let $\Lat=\Zahl^D$ be the $D$-dimensional lattice.
If $\sA$ is a (discretely topologized) finite set, then $\sA^\Lat$ is
compact in the Tychonoff topology. 
  For any $\fv\in\Lat$, let $\shift{\fv}:\Selfmap{\sA^\Lat}$ be the shift
map: $\shift{\fv}(\ba) \ = \ \mtrx{b_\ell}{\ell\in\Lat}{}$, where
$b_\ell = a_{\ell-\fv}$, \ for all $\ell\in\Lat$.  In particular, if
$D=1$, let $\shift{} = \shift{1}$ be the left-shift on $\sA^\Zahl$.
 A {\dfn cellular automaton} (CA) is a continuous map
$\Phi:\Selfmap{\sA^\Lat}$ which commutes with all shifts: \ for any
$\ell\in\Lat$,\ \ \ $\shift{\ell}\circ
\Phi \ = \ \Phi \circ \shift{\ell}$.  A result of Curtis, Hedlund, and
Lyndon \cite{HedlundCA} says any CA is determined by a {\dfn local map}
$\phi:\sA^\Nh \into \sA$ (where $\Nh\subset \Lat$ is some finite
subset), such that, for all $\ell\in\Lat$ and all $\ba\in\sA^\Lat$, if we
define $\ba\restr{\ell+\Nh} \ = \ [a_{\ell+\nh}]_{\nh\in\Nh}
\ \in \ \sA^{\ell+\Nh}$, then $\Phi(\ba)_\ell \ = \
\phi\lb(\ba\restr{\ell+\Nh}\rb)$.

  In \cite{HofKnill}, Hof and Knill studied the action of CA on
`circle shifts', a class of quasiperiodic subshifts similar to
Sturmian shifts \cite{HedlundSturmian,MorseHedlund}.  They showed that
the image of circle subshift under a CA is again a circle shift, and
raised two questions:
\bitem
  \item Empirically, iterating the CA fragments the partition generating
the circle shift. Why?

  \item Is the action of a CA injective when restricted to a circle subshift?
\eitem
  In this paper, we generalize \cite{HofKnill} by studying the action of
CA upon {\dfn quasisturmian} systems.  In \S\ref{S:prelim} we introduce
notation and terminology concerning torus rotations, measurable
partitions of tori, and the Besicovitch metric $d_B$
\cite{BlanchardFormentiKurka}.  We introduce 
$\MP$, the space of measurable partitions of a torus $\Tor$,
endowed with the {\dfn symmetric difference} metric $d_\symdif$. In
\S\ref{S:quasisturm}, we introduce {\dfn quasisturmian} (QS) {\dfn shifts},
 which are natural generalizations of the classical Sturmian shift 
\cite{HedlundSturmian,MorseHedlund},
 obtained by tracking the trajectories of a
torus rotation system $(\Tor,\trans{})$ through
 a fixed open partition of $\Tor$.
Likewise, a {\dfn QS measure} is obtained by projecting
$\trans{}$-trajectories through a {\em measurable} partition of
$\Tor$.  We refer to both QS shifts and QS measures as {\dfn quasisturmian
systems}.  A {\dfn QS sequences} is (roughly speaking) a `typical'
element of a QS shift. We define $\QS$ to be the set of all 
QS sequences in $\sA^\Lat$ generated by the rotation system $(\Tor,\trans{})$.

 In \S\ref{S:CA.on.QS}, we examine the action of a CA on a QS system,
and generalize the results of \cite{HofKnill}, to show that any CA
$\Phi$ induces a natural transformation $\Phitor$ on $\MP$
(Theorem \ref{hof.knill}).  The `induced' dynamical system
$(\MP,d_\symdif,\Phitor)$ is a topological dynamical system, because
$\Phitor$ is Lipschitz relative to $d_\symdif$
(Proposition \ref{Phi.star.continuous}).  There is a natural conjugacy
between a subsystem of $(\MP,d_\symdif,\Phitor)$ and the system
$(\QS,d_B,\Phi)$ (Proposition \ref{Xi.Phi.iso}).

\S\ref{app:MP} contains auxiliary technical results for
\S\ref{S:boundary} and \S\ref{S:injective}. 
 In \S\ref{S:boundary}, we address Hof and Knill's first question, and
show how a partition is `chopped' into many small pieces under
iteration of $\Phitor$.  In \S\ref{S:injective} we address Hof and
Knill's second question, and show that, `generically', a CA restricted
to a QS system is injective (Theorem \ref{comeager.quasisturm.inject}).
\S\ref{app:custom.quasi} is auxiliary to
\S\ref{S:surject} and \S\ref{S:fixedpoints}.
 In \S\ref{S:surject}, we compare the
surjectivity  of $(\MP,d_\symdif,\Phi)$ and $(\sA^\Lat,\Phi)$;
we show that the map $\Phitor\colon\MP\into\MP$ is not generally
surjective, but {\em does} have a $d_\symdif$-dense image in $\MP$ if $\Phi$ is
surjective (Proposition \ref{CA.surjective.on.MP}).  
In \S\ref{S:fixedpoints}
we study QS fixed points, periodic points, and travelling
wave solutions for $\Phi$.

\S \ref{app:LCA} and \S \ref{app:torus} 
contain auxiliary machinery for \S\ref{S:expansive},
\S\ref{S:nilpotence.rigid}, and \S\ref{S:recur}. In \S\ref{S:expansive} we give an example of a CA which acts {\em
expansively} on $(\QS,d_B)$
(Proposition \ref{ledrappier.expansive.on.quasisturm}), thereby refuting a
plausible conjecture arising from
\cite{BlanchardFormentiKurka}.
In \S\ref{S:nilpotence.rigid} we show that {\dfn linear} cellular automata
are either {\dfn niltropic} or {\dfn rigid} when restricted to
a quasisturmian shift (Proposition \ref{LCA.nilpotence.rigid}); 
in \S\ref{S:inv.quasi.measure} we use this to
show that most linear CA have {\em no} quasisturmian invariant measures
(Proposition \ref{no.LCA.inv.quasi.measures}.)
\S\ref{S:recur} constructs a class of quasisturmian measures which
are not asymptotically randomized by a simple linear CA.

   Sections \ref{S:boundary}, \ref{S:injective}, \ref{S:surject},
\ref{S:fixedpoints}, \ref{S:expansive}, 
 \ref{S:nilpotence.rigid}, and \ref{S:recur} are logically independent.
\S\ref{S:inv.quasi.measure} depends on \S\ref{S:nilpotence.rigid}.

\section{Preliminaries and Notation \label{S:prelim}}

  If $\ell,n\in\Zahl$, then $\CO{\ell...n}=\set{m\in\Zahl}{\ell\leq m
< n}$.  If $(\bX,d)$ is a metric space and $x,y\in\bX$, then
``$x \closeto{\eps} y$'' means $d(x,y)\leq \eps$.  Let $\sM(\bX)$ be the
space of (Borel) probability measures on $\bX$.
If $\lam\in\sM(\bX)$, then ``$\forall_\lam \, x$'' means ``$\lam$-almost
every $x$.''  Likewise, ``$\lam$-\ae'' means ``$\lam$-almost everywhere''.
A {\dfn meager} subset of $\bX$ is a nowhere dense set, or any set
obtained through the countable union or intersection of other meager
sets; it is the topological analog of a `set of measure zero'.  
A {\dfn comeager} set is the complement of a meager set.
 A statement holds `topologically almost everywhere'
\topae \  on $\bX$ if it holds for all points in a comeager
subset.

  A {\dfn topological dynamical system} ({\bf TDS}) is a triple
$(\bX,d,\varphi)$ where $(\bX,d)$ is a metric space and $\varphi:\bX\into\bX$
is a continuous transformation.  If $(\barbX,\bard,\bar\varphi)$ is 
another TDS, then a TDS {\dfn epimorphism} is a continuous surjection
$f:\bX\into\barbX$ such that $f\circ\varphi = \bar\varphi\circ f$;
if $f$ is a homeomorphism, then we say $f$ is a TDS {\dfn isomorphism}.

  A {\dfn measure-preserving dynamical system} ({\bf MPDS}) is a triple
$(\bX,\mu,\varphi)$ where $(\bX,\mu)$ is a probability space and
 $\varphi:\bX\into\bX$ is a measurable transformation such that
$\varphi(\mu)=\mu$.   If $(\barbX,\barmu,\bar\varphi)$ is 
another MPDS, then an MPDS {\dfn epimorphism} is a measure-preserving map
$f:\bX\into\barbX$ such that $f\circ\varphi = \bar\varphi\circ f$;  if
$f$ is bijective \muae, then $f$ is an MPDS {\dfn isomorphism}.

  Let $\Lat=\Zahl^D$ be a lattice.  A {\dfn topological $\Lat$-system}
({\bf T$\boldsymbol{\Lat}$S}) is a triple $(\bX,d,\trans{})$ where $(\bX,d)$ is a metric
space and $\trans{}$ is a continuous $\Lat$-action on $\bX$; if
$\ell\in\Lat$, then we write the action of $\ell$ as $\trans{\ell}$.
For example, if $\sA$ is a finite alphabet and $\bX=\sA^\Lat$ is the
space of all $\Lat$-indexed {\dfn configurations} of elements in $\sA$,
then $\Lat$ acts by on $\sA^\Lat$ by shifts; we indicate the shift
action by $\shift{}$, and $(\sA^\Lat,d_C,\shift{})$ is a T$\Lat$S
(where $d_C$ is the Cantor metric --see below).
A {\dfn subshift} is a (Cantor)-closed, $\shift{}$-invariant subset 
$\gX\subset\sA^\Lat$; then $(\gX,d_C,\shift{})$ is also a  T$\Lat$S.

 A {\dfn measure-preserving $\Lat$-system} ({\bf MP$\boldsymbol{\Lat}$S}) is a triple
$(\bX,\mu,\trans{})$ where $(\bX,\mu)$ is a probability space and
$\trans{}$ is a $\mu$-preserving $\Lat$-action on $\bX$.
We define (epi/iso)morphisms of T$\Lat$S and MP$\Lat$S in the obvious way.
A measurable subset $\bU\subset\bX$ is {\dfn
$\trans{}$-invariant} if $\mu\lb[\bU\;\symdif\;\trans{\ell}(\bU)\rb] \
= \ 0$ for all $\ell\in\dL$, and $(\bX,\mu,\trans{})$ is {\dfn ergodic} if,
for any $\trans{}$-invariant set $\bU$, either $\mu[\bU]=0$ or
$\mu[\bU]=1$.
 For any $N>0$, let $\dB(N) = \CO{-N..N}^D \subset\Lat$.  Thus,
$|\dB(N)| = (2N)^D$.

\paragraph*{Generalized Ergodic Theorem:} 
{\sl
 If  $(\bX,\mu,\trans{})$ is ergodic, and
$\bU\subset\bX$ is measurable, then for all $_\mu \, x\in\bX$ are
{\dfn $(\mu,\trans{})$-generic} for $\bU$, meaning that

\hfill 
$\D \mu[\bU]  = \ 
\lim_{N\goto\oo} \ \frac{1}{(2N)^D} \sum_{\fb\in\dB(N)} \ 
\chr{\bU}\lb(\maketall\trans{\fb}(x)\rb)$.\hfill\qed
}\cite{OrnsteinWeissErgodic,Tempelman}
\breath

For example, if $\mu\in\sM(\sA^\dL)$ is {\dfn
$\shift{}$-invariant} (ie. $\shift{\ell}(\mu) \ = \ \mu$, \ 
for all $ \, \ell\in\dL$) then $(\sA^\dL,\mu,\shift{})$ is  an
MP$\Lat$S.  An element $\ba\in\sA^\Lat$ is {\dfn $\mu$-generic}
if $\ba$ is $(\mu,\shift{})$-generic for all cylinder sets of $\sA^\Lat$.

\sss{The Cantor and Besicovitch metrics}

  The standard (Tychonoff) topology on $\sA^\Lat$ is induced by the 
{\dfn Cantor metric:}
\beq
\mbox{For any $\bp,\bq\in\sA^\Lat$,}\quad
  d_C(\bp,\bq) &:=& 2^{-D(\bp,\bq)},
\\
\mbox{where} \quad
D(\bp,\bq) & := &  \min\set{|\ell|}{\ell\in\Lat, \ p_\ell\neq q_\ell}.
\eeq
The topological dynamics of CA in this metric
were characterized in \cite{Kurka,BlanchKurkaMaass}.   For our
purposes, however, it is more appropriate to use the {\dfn Besicovitch
metric} \cite{BlanchardFormentiKurka,CatFormMargMaz} (see also
\cite{FormentiBlanchardCervelle,FormentiBesi}), defined as follows.
 If $\dJ\subset\Lat$, then the {\dfn \Cesaro density} of $\dJ$ is
defined:
\[
  \density{\dJ} \quad:=\quad \lim_{N\goto\oo} 
\frac{\#(\dJ\intsct\dB(N))}{(2N)^D},
\]
 if this limit exists.   If not, then the 
{\dfn upper \Cesaro density} of $\dJ$ is the limsup:
\[
  \updensity{\dJ} \quad:=\quad \limsup_{N\goto\oo} 
\frac{\#(\dJ\intsct\dB(N))}{(2N)^D}.
\]
  We define the {\dfn Besicovitch (pseudo)metric}
on $\sA^\Lat$ as follows: 
\[
  d_B(\bp,\bq) \quad=\quad \updensity{\ell\in\Lat\ ; \ p_\ell \neq q_\ell},
\qquad \mbox{for any $\bp,\bq\in\sA^\Lat$}.
\]
This is only a {\em pseudo}metric because it
is possible for $d_B(\bp,\bq)=0$ while $\bp\neq \bq$, as long as $\bp$
and $\bq$ disagree on a set of  upper \Cesaro density zero.
Thus, we identify any element $\bp\in\sA^\dL$ with its 
Besicovitch equivalence class $\tlbp=\set{\bq\in\sA^\dL}{d_B(\bp,\bq)=0}$.
Let $\widetilde{\sA^\Lat}$ be the set of equivalence classes.
 If $\Phi:\sA^\Lat\into\sA^\Lat$ is a CA, then
 $\Phi$ factors to a $d_B$-continuous map $\tl\Phi:\widetilde{\sA^\dL}
\into \widetilde{\sA^\dL}$ \cite{BlanchardFormentiKurka,CatFormMargMaz}.
We'll usually abuse notation by writing $\tl\Phi$ as ``$\Phi$'' and
$\widetilde{\sA^\Lat}$ as ``$\sA^\Lat$''.  The topological dynamics of
CA in the Besicovitch metric were investigated in \cite{BlanchardFormentiKurka,CatFormMargMaz}.

\sss{Torus Rotation Systems}

  Let $\Torus{1} := \Real/\Zahl$, which we normally identify with $\CO{0,1}$.
 Fix $K\geq1$ and let $\Tor := \Torus{K} = \Torus{1}\x\cdots\x\Torus{1}$ be
the $K$-torus.  For any $\tors\in\Tor$, let $\rot{\tors}:\Tor\ni \tort
\ \mapsto \ (\tort+\tors)\in\Tor$ be the corresponding rotation map.
Suppose $\tau:\Lat\into\Tor$ is a group monomorphism; for any
$\ell\in\Lat$, let $\trans{\ell} = \rot{\tau(\ell)}$ denote the
corresponding rotation of $\Tor$.  This defines a measure-preserving,
topological $\Lat$-system $(\Tor,d,\lam,\varsigma)$ (where $d$ is the
usual metric and $\lam$ is the Lebesgue measure on $\Tor$). We call
this a {\dfn torus rotation system}.  

\Proposition{\label{tau.monothetic.uniqely.ergodic}}
{
Suppose $\tau\colon\Lat\rightarrow\Tor$ is a monomorphism
with dense image.  Then:
\\
{\bf(a)} $(\Tor,d,\trans{})$ is minimal 
{\em (ie. every $\tort\in\Tor$ has dense $\trans{}$-orbit)} and
 uniquely ergodic {\em (ie. $\lam$ is the only $\trans{}$-invariant
probability measure)}.
\\
{\bf(b)}  Let $\bU\subset\Tor$ be open, with $\lam[\partial\bU]=0$.
Then every $\tort\in\Tor$ is $(\lam,\trans{})$-generic for $\bU$.
}
\bthmprf {\bf(a)} {\em Minimal:} 
 Suppose $\bZ:=\tau(\Lat)$ is dense in
$\Tor$.  Then for any $\tort\in\Tor$, \ the $\trans{}$-orbit
 $\{\trans{\ell}(\tort)\}_{\ell\in\Lat} = \tort+\bZ$ is also dense.
 {\em Uniquely Ergodic:}   If $\mu$ is a
 $\trans{}$-invariant probability measure on $\Tor$, then
$\rot{\torz}(\mu)=\mu$ for every $\torz\in\bZ$.  But
$\bZ$ is dense in $\Tor$, so (by a weak*-convergence argument), we get
$\rot{\tort}(\mu)=\mu$ for   all $\tort\in\Tor$.  But the
Haar measure $\lam$  is the {\em only} probability measure on $\Tor$ that
is $\rot{\tort}$-invariant for all $\tort\in\Tor$ \cite[Thm 10.14, p.317]{Folland}; hence $\mu=\lam$.

{\bf(b)}  For any $\tors\in\Tor$ and $\bW\subset\Tor$, let
$\dI(\tors,\bW):= \set{\ell\in\Lat}{\trans{\ell}(\tors)\in\bW}$,
and let $D(\tors,\bW):=\density{\dI(\tors,\bW)}$.
We want to show that  $D(\tort,\bU) = \lam[\bU]$.

\Claim{$\lam[\bU] \leq D(\tort,\bU)$.}
\bclaimprf
 For any $\del>0$, let $\bU_\del :=\set{\toru\in\bU}{d(\toru,\compl{\bU})>\del}$.  Then $\Union_{\del>0} \bU_\del = \bU$ (because if $\toru \in \bU$, but
$\toru\not\in\bU_\del$ for any $\del>0$, then there is no $\del$-ball around
$\toru$ contained in $\bU$, contradicting that $\bU$ is open).  Thus,
$\D\lim_{\del\searrow 0} \  \lam[\bU_\del]=\lam[\bU]$.

  Given $\eps>0$, find $\del>0$ such that $\lam[\bU_\del]> \lam[\bU]-\eps$.
Let $\bG:=\{\tort\in\Tor$; $\tort$ is generic for both
$\bU$ and $\bU_\del$ $\}$.  The Generalized Ergodic Theorem says
that $\lam[\bG]=1$; hence $\bG$ is dense in $\Tor$.
  Find  $\torg\in\bG$ with $\torg\closeto{\del}\tort$.  Thus,
for any $\ell\in\Lat$,
$\statement{$\trans{\ell}(\torg)\in\bU_\del$}
\IMPLIES
\statement{$\trans{\ell}(\tort)\in\bU$}$
(because $\trans{\ell}$ is an isometry).
Thus, $\dI(\torg,\bU_\del)\subseteq\dI(\tort,\bU)$.
Hence $D(\tort,\bU) \ \geq \ D(\torg,\bU_\del) = \lam[\bU_\del] >
\lam[\bU]-\eps$.
  Since $\eps$ is arbitrary, we conclude that $D(\tort,\bU)\geq \lam[\bU]$.
\eclaimprf
  Now, let $\bV=\mathrm{int}(\compl{\bU})$.  Then $\bV$ is also open, and
\[ 
 \lam[\bU] \quad \leeeq{(\dagger)} \quad D(\tort,\bU) \quad \leeeq{(*)} \quad
 1-D(\tort,\bV)  \quad \leeeq{(\ddagger)} \quad
1 - \lam[\bV] \quad\eeequals{(\diamond)}\quad \lam[\bU].
\]
Thus $D(\tort,\bU)=\lam[\bU]$, as desired.  
$(\dagger)$ is by Claim 1.
$(*)$ is because  $\dI(\tort,\bV)$ is disjoint from $\dI(\tort,\bU)$.
$(\ddagger)$ is because
Claim 1 (applied to $\bV$) yields $D(\tort,\bV)\geq \lam[\bV]$.
$(\diamond)$ is because  $\compl{(\bV\disj\bU)}
= \partial\bU$, so $\lam[\bV]+\lam[\bU] = 
\lam[\bV\disj\bU] = 1-\lam[\partial\bU]=1-0=1$.\nopagebreak
\ethmprf

{\sc Note:}  {\em We assume throughout this paper
that the hypothesis of Proposition {\rm\ref{tau.monothetic.uniqely.ergodic}} 
 is satisfied.}

\example{\label{X:circle.rotation}
Let $D=1$, so $\Lat=\Zahl$.  Let $K=1$, so 
$\Tor=\Torus{1}$.  Identify $\Tor$ with $\CO{0,1}$ in the obvious way.
Let $\tora\in\CO{0,1}$ be irrational; define
$\tau:\Zahl\into\Torus{1}$ by $\tau(z)= z\cdot \tora$ (mod $1$).  
Thus, $\trans{z}(\tort) = \tort+z\tora, \ \forall z\in\Zahl, \, \forall \tort\in\Tor$.
Thus, $(\Torus{1},\lam,\varsigma)$ is an irrational rotation of a
circle.}

\sss{Measurable partitions}

An {\dfn ($\sA$-labelled) measurable partition} of $\Tor$ is a
finite collection of disjoint measurable sets $\sP = \{\bP_a\}_{a\in\sA}$ so
that, if $\bP_* := \D\Disj_{a\in\sA} \bP_a$, then $\lam(\bP_*)=1$. 
 
 We'll often treat the partition $\sP$ as a measurable function
$\sP:\bP_*\into\sA$, where $\sP^{-1}\{a\} := \bP_a$.  Let $\MP$ be the
set of all measurable $\sA$-labelled partitions of $\Tor$, which
we topologize with the {\dfn symmetric difference} metric, defined:
\beqn
\label{sym.dif.metric}
 d_\symdif (\sP,\sQ)
\quad = \quad \sum_{a\in\sA} \lam\lb(\bP_a\symdif\bQ_a\rb),
\qquad \mbox{for any $\sP,\sQ\in\MP$}.
\eeqn
   $(\MP,d_\symdif)$ is a complete and bounded metric space, but not
compact (Proposition  \ref{MP.complete}).

An {\dfn $\sA$-labelled open partition} of $\Tor$ is a finite family of
disjoint open sets $\sP = \{\bP_a\}_{a\in\sA}$, such that $\bP_*$ 
is a dense open subset of
$\Tor$, and $\lam(\bP_*)=1$. If $\OP$ is the set of 
$\sA$-labelled open partitions of $\Tor$, then 
$\OP\subset\MP$, and $\OP$ is $d_\symdif$-dense in
$\MP$ (Corollary \ref{OP.dense.in.MP}).

\section{Quasisturmian systems \label{S:quasisturm}}

 Quasisturmian systems are $\shift{}$-invariant subsets or measures in
$\sA^\Lat$ which  generalize the classical Sturmian shift  of
\cite{HedlundSturmian,MorseHedlund}
(see also \cite{Arnoux,Berstel,Berthe,Brown,BlanchardKurka,BertheVuillon}).
  \nocite{MaassBlanchardNogueira,PytheasFogg}

\sss{Quasisturmian shifts}

Let $\sP= \{\bP_a\}_{a\in\sA}$ be an open partition of $\Tor$.
Let $\bP_* = \D\Disj_{a\in\sA} \bP_a$, and
let $\tlTor=\tlTor(\sP)  :=  
 \D\Intsct_{\ell\in\Lat}\trans{\ell}(\bP_*)$; then $\tlTor$ is a $\trans{}$-invariant, dense $G_\delta$ subset of $\Tor$, and $\lam(\tlTor)=1$.   (We will
write $\tlTor(\sP)$ as $\tlTor$ when the partition $\sP$ is 
clear from context.)

\example{\label{X:sturm1}
Let $K=1$ so $\Tor \ \cong \ \CO{0,1}$ as in Example
\ref{X:circle.rotation}.  Let $\sA=\{0,1\}$,
and let $\bP_0 = \OO{0,\tora}$ and $\bP_1 = \OO{\tora,1}$.
  Then $\bP_* \ \cong \ \OO{0,1}\setminus\{\tora\}$, and
$\tlTor \ \cong \ \CO{0,1} \setminus \set{z\cdot\tora \bmod{1}}{z\in\Zahl}$.}
 
 Recall that $\sP:\tlTor\into\sA$ is defined by $\sP^{-1}\{a\} =
\bP_a$. For any $\tort\in\tlTor$, let $\Splat(\tort)\in\sA^\Lat$ be
the $\sP$-trajectory of $\tort$.  That is, for all $\ell\in\Lat$, \quad
$\Splat(\tort)_\ell = \sP\lb(\trans{\ell}(\tort)\rb)$.  This defines a
function $\Splat:\tlTor\into\sA^\Lat$.

\Proposition{\label{splat.homeo}}
{
  Let $\sP\in\OP$ be a nontrivial open partition of $\Tor$.  Then:
\\
{\bf(a)} \  $\Splat\circ\trans{\ell} = \shift{\ell}\circ \Splat$ for all
$\ell\in\Lat$. 
\\
{\bf(b)} \ $\Splat:\tlTor \into \sA^\Lat$ is continuous with respect
to both the $d_C$ and $d_B$ metrics on $\sA^\Lat$.
}
\bthmprf {\bf(a)} is by definition of $\Splat$.
 
{\bf(b)}\quad
{\em $d_C$-continuous:}  \
Let $\dM\subset\Lat$ be any finite subset; we want a neighbourhood $\bU$
around $\tort$ such that, if $\tort'\in\bU$, then 
$p_\fm = p'_\fm$ for all $\fm\in\dM$.
Let $\sQ=\bigvee_{\fm\in\dM} \trans{\fm}(\sP)$.
The atoms of $\sQ$ are  finite intersections
of open sets, hence open. Let $\bU$ be the $\sQ$-atom containing $\tort$.
If $\tort'\in\bU$ then for all $\fm\in\dM$, \ 
$\sP(\trans{\fm}(\tort)) =\sP(\trans{\fm}(\tort'))$,
i.e. $p_\fm=p'_\fm$, as desired.

{\em $d_B$-continuous:} [See also Proposition \ref{splat.is.lipschitz}(a)]
 \ Fix $\eps>0$.  
We want $\del>0$ such that, if $\tort\closeto{\del}\tort'$, then
$d_B(\bp,\bp')<\eps$. For any $\del>0$ let $\bP_{\del} :=
\set{\torp\in\bP_*}{d\lb(\torp,\compl{(\bP_*)}\rb)>\del}$.
Find $\del$ small enough that $\lam[\bP_{\del}]>1-\eps$
(see second paragraph in Claim 1 of Proposition \ref{tau.monothetic.uniqely.ergodic}).
If $\tort\closeto{\del}\tort'$, then
\beq
d_B(\bp,\bp') &=&
 \density{\ell\in\Lat \ ; \ p_\ell\neq p'_\ell}
 \quad \leeeq{(\dagger)} \quad 
 \density{\ell\in\Lat \ ; \ \trans{\ell}(\tort)\not\in\bP_\del}
\\
& \eeequals{(*)} &  \lam\lb[\compl{(\bP_\del)}\rb]
 \quad < \quad \eps.
\eeq
 $(*)$ is by
Proposition \ref{tau.monothetic.uniqely.ergodic}(b).
 $(\dagger)$ is because, for any $\ell\in\Lat$,
 $ \statement{$p_\ell \neq p'_\ell$}
\IMPLIES
\statement{$\trans{\ell}(\tort)\not\in\bP_\del$}$.
To see this, suppose $p_\ell=a$
ie. $\trans{\ell}(\tort)\in\sP_a$.  If $\trans{\ell}(\tort)\in\bP_\del$,
then $
d\lb(\trans{\ell}(\tort),\trans{\ell}(\tort')\rb)   \leq   \del 
 < 
d\lb(\trans{\ell}(\tort), \compl{(\bP_*)}\rb)
 \leq  
d\lb(\trans{\ell}(\tort), \compl{(\bP_a)}\rb)$.
Hence $\trans{\ell}(\tort')\in\sP_a$ also;  hence
$p'_\ell = a = p_\ell$. By contradiction, if $p_\ell \neq p'_\ell$,
then $\trans{\ell}(\tort)\not\in\bP_\del$.
\ethmprf

 The {\dfn $\trans{}$-quasisturmian} (or $\mathbf{QS_\trans{}}$) {\dfn shift} 
induced by $\sP$ is the $d_C$-closed, $\shift{}$-invariant subset:
\[
  \Qshift{\sP} \quad:=\quad d_C\!\!-\!\!\closure{\Splat(\tlTor)} \quad\subset\quad \sA^\dL.
\]
\example{\label{X:sturm2}
If $K=1$ and $\sP = \lb\{\OO{0,\tora},\OO{\tora,1}\rb\}$
as in Example \ref{X:sturm1}, then $\Qshift{\sP}\subset\sA^\Zahl$
is the classical Sturmian shift of \cite{HedlundSturmian,MorseHedlund}.}

\sss{Quasisturmian Sequences}

We define the set of {\dfn $\trans{}$-quasisturmian sequences}:

\breath

\centerline{$\QS
 \quad:=\quad \set{\Splat(\tort)}{\sP\in\OP, \ \tort\in\tlTor(\sP)}.$}

Then $\QS$ is a $d_C$-dense subset of $\sA^\Lat$ (Corollary \ref{QP.is.dense}).

\Proposition{\label{Xi.iso}}
{
Let $\ZP  =  \set{\sP\in\OP}{0 \in \tlTor(\sP)}$.
For  all $\sP\in\ZP$, let 
$\Qseq{\sP} = \Splat(0)$.
\bthmlist
  \item  $\ZP$ is a comeager, $\trans{}$-invariant subset of $\OP$,
and  $\Qseq{\ZP}=\QS$.

  \item  $\xi_\trans{}:\ZP\into\QS$ is a distance-halving isometry.  
That is,  if
 $\sP,\sQ\in\ZP$, then $d_\symdif(\sP,\sQ) \ = \ 
2\cdot d_B\lb(\maketall \Qseq{\sP}, \Qseq{\sQ}\rb)$.

  More generally, if  $\sP,\sQ\in\OP$ and $\tort\in\tlTor(\sP)\intsct\tlTor(\sQ)$, then $d_\symdif(\sP,\sQ)  =  
2\cdot d_B\lb(\maketall \Splat(\tort),\Sqlat(\tort)\rb)$.

  \item   If $\bp,\bq\in\QS$, then $\statement{$d_B(\bp,\bq)=0$}\iff\statement{$\bp=\bq$}$.  
 Hence, $d_B$ is a true metric when restricted to $\QS$.

  \item   For any $\ell\in\dL$, \ 
$\xi_\trans{}\circ \trans{\ell}  =  \shift{\ell}\circ \xi_\trans{}$.
 Thus, $\QS$ is a $\shift{}$-invariant subset of $\sA^\dL$.

  \item  $\xi_\trans{}$ is a TDS isomorphism from 
$(\ZP,d_\symdif,\trans{})$ to $(\QS,d_B,\shift{})$.
\ethmlist
}
\bthmprf {\bf(a)} It is clear that $\ZP$ is a comeager and
$\trans{}$-invariant set.   Let $\bq\in\QS$; hence there is some
 $\sQ\in\OP$ and $\tort\in\tlTor$ such that $\bq=\Sqlat(\tort)$.  Define
partition $\sP\in\OP$ such that $\bP_a \ := \ \set{\torq-\tort}{\torq\in\bQ_a}$
for all $a\in\sA$.  Then $\sP\in\ZP$, and $\Splat(0)=\Sqlat(\tort)=\bq$.

{\bf(b)}
Let $\bp=\Qseq{\sP}$ and $\bq=\Qseq{\sQ}$.
If $a\in\sA$, the Generalized  Ergodic Theorem says:
\beq
\lam\lb(\bP_a\setminus\bQ_a\rb)  
&=& \lim_{N\goto\oo} 
\frac{\#\set{\fb\in\dB(N)}{p_\fb \ = \ a \ \neq \  q_\fb}}{(2N)^D}.\\
\mbox{Thus,}\quad
\sum_{a\in\sA} \ \lam\lb(\bP_a\setminus\bQ_a\rb)  
&=&
\sum_{a\in\sA} \ \lim_{N\goto\oo} \ 
\frac{\#\set{\fb\in\dB(N)}{p_\fb \ = \ a \ \neq \  q_\fb}}{(2N)^D}
\\
\lefteqn{\hspace{-5em}=\quad
\lim_{N\goto\oo}  \
\frac{\#\set{\fb\in\dB(N)}{p_\fb \neq \  q_\fb}}{(2N)^D}
\quad=\quad d_B\lb(\maketall \Qseq{\sP}, \Qseq{\sQ}\rb).}
\eeq
Likewise, \ $\D \sum_{a\in\sA} \lam\lb(\bQ_a\setminus\bP_a\rb)  
\ = \ d_B\lb(\maketall \Qseq{\sP}, \Qseq{\sQ}\rb)$.

  Hence, \ 
$\D d_\symdif(\sP,\sQ)  \ = \  
\sum_{a\in\sA} \lam\lb(\bP_a\setminus\bQ_a\rb)  
+ \sum_{a\in\sA} \lam\lb(\bQ_a\setminus\bP_a\rb) 
\ = \  2\cdot  d_B\lb(\maketall \Qseq{\sP}, \Qseq{\sQ}\rb)$.

{\bf(c)} follows from {\bf(b)}.  {\bf(d)} follows from the
definitions, and {\bf(e)} follows from {\bf(d)}.\nolinebreak
\ethmprf

\sss{Quasisturmian Measures}

  Let $\lam$ be the Lebesgue measure on $\Tor$.  Let $\sP\in\MP$ be a
  {\em measurable} partition of $\Tor$.  Define $\tlTor$ as
  before, then $\tlTor$ is a $\trans{}$-invariant, measurable subset
  of $\Tor$, and $\lam(\tlTor)=1$.  Define
  $\Splat:\tlTor\into\sA^\dL$ as before; then $\Splat$ is a measurable
  function defined $\lam$-\ae\ on $\Tor$.  Define
  $\Upsilon_\trans{}:\MP\into\sM(\sA^\Lat)$ by $\QM{\sP} :=
  \Splat(\lam)$, \ for all $ \ \sP\in\MP$.  Then
  $\Upsilon_\trans{}(\sP)$ is a $\shift{}$-invariant measure on
  $\sA^\dL$, called the {\dfn $\trans{}$-quasisturmian} 
(or $\mathbf{QS_\trans{}}$) {\dfn measure} induced by $\sP$.
If $\QSM$ is the
set of $\trans{}$-quasisturmian measures, 
then $\QSM$ is weak*-dense in the space of $\shift{}$-ergodic
probability measures on $\sA^\Lat$ (Corollary \ref{MQP.is.wkst.dense}).

  Quasisturmian measures, sequences, and shifts are related as
follows: If $\sP$ is an open partition, then
$\supp{\Upsilon_\trans{}(\sP)}=\Qshift{\sP}$.  Also, if
$\tort\in\tlTor$, then $\Qshift{\sP}$ is the $\shift{}$-orbit closure
of $\Splat(\tort)$, and $\Splat(\tort)$ is $\shift{}$-generic for
$\Upsilon_\trans{}(\sP)$.

\Proposition{\label{QM.continuous}}
{
 $\Upsilon_\trans{}:\MP\into\QSM$ is continuous relative to
$d_\symdif$ and the weak* topology.
}
\bthmprf
  Suppose $\{\sP^{(n)}\}_{n=1}^\oo\subset\MP$ is a sequence of partitions, and
$\D \symdiflim_{n\goto\oo}\sP^{(n)} \ = \ \sP$.
  Let $\mu^{(n)}=\QM{\sP^{(n)}}$
for all $n$, and let $\mu=\QM{\sP}$;
we claim that $\D\wkstlim_{n\goto\oo} \mu^{(n)} \ = \ \mu$.

  Suppose $\dM\subset\dL$ is finite, and let $\bw\in\sA^\dM$;
we must show that $\D\lim_{n\goto\oo}\mu^{(n)}[\bw] \ = \ \mu[\bw]$.  
Let $\bw := [w_\fm]_{\fm\in\dM}$.  If $\sP=\{\bP_a\}_{a\in\sA}$ and
$\sP^{(n)} \ = \ \{\bP^{(n)}_a\}_{a\in\sA}$ then
\[
 \mu[\bw] \quad = \quad \lam\lb[\Splat^{-1}\{\bw\}\rb]
\quad = \quad \lam\lb[\Intsct_{\fm\in\dM} \trans{\fm} (\bP_{w_\fm})\rb].
\]
Likewise,
 $\D \mu^{(n)}[\bw] 
 =   \lam\lb[\Intsct_{\fm\in\dM} \trans{\fm} (\bP^{(n)}_{w_\fm})\rb]$.
Let  $M:=\#(\dM)$.  If $d_\symdif(\sP,\sP^{(n)})  <  \eps$,  
then 
\[ 
\lam\lb[\Intsct_{\fm\in\dM} \trans{\fm} (\bP_{w_\fm})\rb]
\quad \closeto{2 M\cdot\eps} \quad 
\lam\lb[\Intsct_{\fm\in\dM} \trans{\fm} (\bP^{(n)}_{w_\fm})\rb];
\] this can be seen by
setting $J:=1$ and $K:=M$
in Lemma \ref{Phi.star.continuous.C3}(c) below.
\ethmprf

  The proof of Proposition \ref{QM.continuous} 
(and later,  Theorem \ref{Phi.star.continuous}) uses  the following 
lemma:

\Lemma{\label{Phi.star.continuous.C3}}
{
Let $\{\bP_i\}_{i=1}^I$ and
$\{\bO_i\}\,_{i=1}^I$ be measurable sets, with $\bP_i\subset\bO_i$ and
$\lam(\bO_i\setminus \bP_i)<\eps$, for all $i\in\CC{1..I}$.
\bthmlist
  \item If $ \barbP \ = \  \Union_{i=1}^I \bP_i$ and 
$ \barbO \ = \  \Union_{i=1}^I \bO_i$, then  $\barbP\subset\barbO$,
and $\lam(\barbO\setminus\barbP) \ \leq I\cdot\eps$.

  \item If $ \undbP  \ = \  \Intsct_{i=1}^I \bP_i$ and 
$ \undbO \ = \  \Intsct_{i=1}^I \bO_i$, then  $\undbP\subset\undbO$,
and $\lam(\undbO\setminus\undbP) \ \leq I\cdot\eps$.

  \item
Let $\{\bP^j_k\}_{j=1}^J\,_{k=1}^K$ and
$\{\bQ^j_k\}_{j=1}^J\,_{k=1}^K$ be measurable sets, with
$\lam(\bQ^j_k\symdif \bP^j_k)<\eps$, for all $j\in\CC{1..J}$ 
and $k\in\CC{1..K}$.

If \
$\D\bP  \ = \ \Union_{j=1}^J  \Intsct_{k=1}^K \bP^j_k$
\ and \
$\D\bQ \ = \ \Union_{j=1}^J \Intsct_{k=1}^K \bQ^j_k$, \ then
$\lam(\bQ \symdif \bP) \ < \ 2 JK\cdot\eps$.
\ethmlist
}
\bthmprf
{\bf(a):}  $\D \barbO\setminus\barbP \ = \
\Union_{i=1}^I (\bO_i\setminus\barbP)
 \ \subseteq \ \Union_{i=1}^I (\bO_i\setminus\bP_i)$,
so $\lam[\barbO\setminus\barbP]
\ \leq \  \D \sum_{i=1}^I \lam[\bO_i\setminus\bP_i]
 \leq  I \eps$.

{\bf(b):}\quad Let $\bO'_i := \bP_i^c$ and
$\bP'_i := \bO_i^c$ for $i\in\CC{1..I}$;  then $\bP'_i\subset\bO'_i$
and $\lam(\bO'_i\setminus\bP'_i) < \eps$.  
Now let  $\barbP' := \Union_{i=1}^I \bP_i'$ and
 $\barbO' := \Union_{i=1}^I \bO_i'$.  Then
$\barbO'\setminus\barbP'
\ = \ \undbP\setminus\undbO$, and
 {\bf(a)} implies $\lam(\barbO'\setminus\barbP')  \ < \ I\eps$. 

{\bf(c)}\quad
For all $j$ and $k$, let $\bO^j_k \ := \ \bP^j_k \union \bQ^j_k$.
Thus, $\bP^j_k\subset \bO^j_k$ and  $\lam(\bO^j_k\setminus \bP^j_k)<\eps$.
Likewise $\bQ^j_k\subset \bO^j_k$ and 
$\lam(\bO^j_k\setminus \bQ^j_k)<\eps$.
\quad  Now, for each $j\in\CC{1..J}$,
 let 
\[
\bO^j \ := \ \Intsct_{k=1}^K \bO^j_k;\qquad
\bP^j \ := \ \Intsct_{k=1}^K \bP^j_k;\quad\AND\quad
\bQ^j \ := \ \Intsct_{k=1}^K \bQ^j_k. 
\]
Thus, setting $I=K$ in part {\bf(b)} implies
 $\bP^j\subset\bO^j$ and
$\lam(\bO^j\setminus\bP^j)  \ < \ K\cdot \eps$.
Likewise, $\bQ^j\subset\bO^j$ and
$\lam(\bO^j\setminus\bQ^j)  \ < \ K\cdot \eps$.

Now let $\D\bO  :=  \Union_{j=1}^J \bO^j  =  \Union_{j=1}^J
\Intsct_{k=1}^K \bO^j_k$, and observe that $\D\bP  = 
\Union_{j=1}^J \bP^j$ and $\D\bQ  =  \Union_{j=1}^J \bQ^j$.

Thus,  setting $I=J$ in part {\bf(a)} implies that that $\bP\subset\bO$ and
$\lam(\bO\setminus\bP)  \ < \ JK\cdot \eps$.
Likewise, $\bQ\subset\bO$ and
$\lam(\bO\setminus\bQ)  \ < \ JK\cdot \eps$.  Thus, by the triangle
inequality, $\lam(\bP\symdif\bQ)=2 JK\eps$.
\ethmprf

\section{CA on QS systems:  induced dynamics on $\MP$
\label{S:CA.on.QS}}

  We begin  by generalizing a result of Hof and Knill \cite{HofKnill}.

\Theorem{\label{hof.knill}}
{ 
Let $\Phi:\Selfmap{\sA^\Lat}$ be a cellular automaton. 
\bthmlist
  \item $\Phi(\QS)\subseteq\QS$.  That is:
if $\bp\in\sA^\Lat$ is a ${\rm QS}_\trans{}$ sequence, then $\Phi(\bp)$ is also a ${\rm QS}_\trans{}$  sequence.   

  \item  If $\gP\subset\sA^\Lat$ is a ${\rm QS}_\trans{}$ shift, then $\Phi(\gP)$ is also a ${\rm QS}_\trans{}$ shift.

  \item  $\Phi[\QSM]\subseteq\QSM$. That is:
  if $\mu$ is a ${\rm QS}_\trans{}$ measure, then
$\Phi(\mu)$ is also  a ${\rm QS}_\trans{}$ measure.
\qed
\ethmlist
} 

\breath

To prove Theorem \ref{hof.knill}, suppose $\Phi$ has local map
$\phi:\sA^\Nh\into\sA$ (where $\Nh\subset\Lat$ is finite).
Suppose $\sP\in\MP$ is a measurable partition of $\Tor$.   For
each $a\in\sA$, define
 \beqn
\label{partition.image}
  \bQ_a \quad = \quad 
\Union_{{{\scriptstyle \bc\in\sA^\Nh}\atop{\scriptstyle \phi(\bc)=a}}} \ \ \Intsct_{\nh\in\Nh}
\ \trans{-\nh} (\bP_{c_\nh}) \quad \subset \quad \Tor.
\eeqn

  Now define measurable partition $\sQ:=\{\bQ_a\}_{a\in\sA}$.
We write:  `$\sQ = \Phitor(\sP)$'. 
Thus, $\Phi$ induces a map $\Phitor:\MP\into\MP$.   It is easy to
verify:

\Lemma{\label{induced.map.on.partitions}}
{
Let $\sP\in\MP$ and let $\sQ=\Phitor(\sP)$.  Then
$\tlTor(\sP)  =  \tlTor(\sQ)$.  Also:
\bthmlist
\item If $\sP$ is open then $\sQ$ is also open. \
 If $\tort \in\tlTor$,  then
$\Phi\lb(\Splat(\tort)\rb)=\Sqlat(\tort)$.

\item  $\Phi\lb(\Qshift{\sP}\rb)=\Qshift{\Phitor(\sP)}$.

\item  $\Phi\lb(\QM{\sP}\rb)=\QM{\Phitor(\sP)}$.
\qed
\ethmlist}
Theorem \ref{hof.knill} follows: 
set $\bp=\Splat(\tort)$, \  $\gP=\Qshift{\sP}$, or $\mu=\QM{\sP}$ 
 in Lemma \ref{induced.map.on.partitions}.

\begin{figure}[hptf]
\centerline{
\begin{tabular}{rl}
\multicolumn{2}{c}{
\psfrag{T}[][]{$\bT\cong [0,1)$}
\psfrag{P1d}[][]{$\bP_1 \ := \ \OO{0,\frac{1}{2}}$}
\psfrag{P0d}[][]{$\bP_0 \ := \ \OO{\frac{1}{2},1}$}
\psfrag{P1}[][]{$\bP_1$}
\psfrag{sP1}[][]{$\trans{1}(\bP_1)  \ = \ \OO{\tora,\frac{1}{2}+\tora}$}
\psfrag{Q1}[][]{$\bQ_1 \ = \ \bP_1\symdif  \trans{1}(\bP_1) \ = \  \OO{0,\tora} \union  \OO{\frac{1}{2},\frac{1}{2}+\tora}$} 
\includegraphics[scale=0.7]{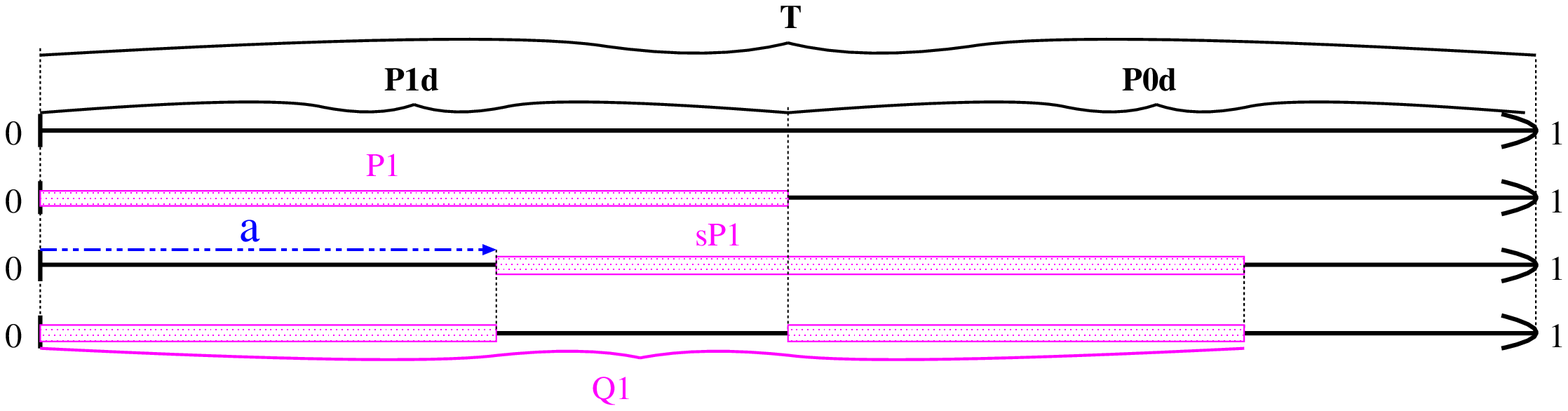}} \\
\includegraphics[height=8.5em,width=18em]{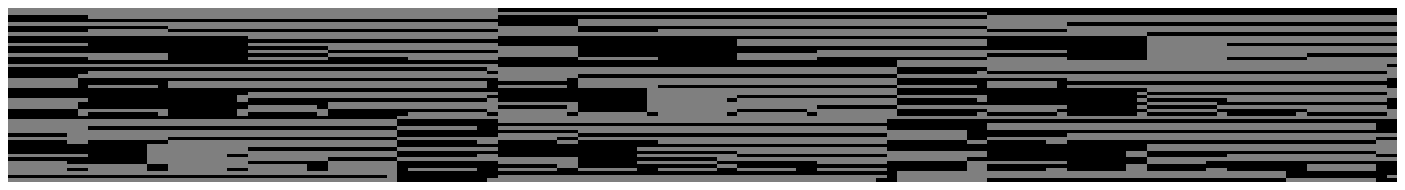}&
\includegraphics[height=8.5em,width=22em]{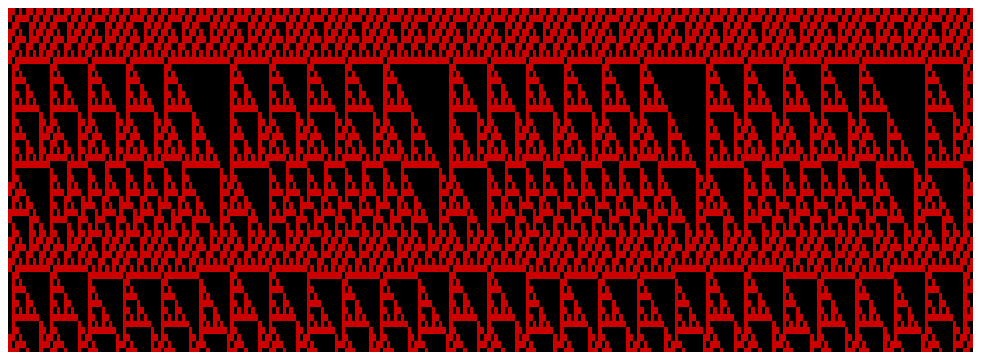} 
\end{tabular}}
\caption{\footnotesize 
Let $\Phi$ be the linear CA from Example \ref{X:ledrappier}, and
let $\Tor=\Torus{1}$.
{\bf Top:} \ The action of $\Phitor$ on $\sP=\{\bP_0,\bP_1\}$, where
 $\bP_1:=\OO{0,\frac{1}{2}}$ and 
$\bP_0:=\OO{\frac{1}{2},1}$.
{\bf Bottom:} \ 
Let $\sP$ be the partition of 
from Example \ref{X:sturm1}, with $\tora\approx 0.3525...$.
Let $\bp\in\Qshift{\sP}$ be an element of the
corresponding Sturmian shift (Example \ref{X:sturm2}).
{\bf Left:} \ Fifty iterations of $\Phitor$ on $\sP$.
{\bf Right:} \  Fifty iterations of $\Phi$ on $\bp$. 
\label{fig:X:ledrappier}}
\end{figure}

\breath

\Examples{Linear Cellular Automata \ {\rm(see \S\ref{app:LCA})}}
{
\item \label{X:ledrappier}
 Let $\Lat=\Zahl$.
Let $\sA=\Zahlmod{2}$ and let $\Phi$ have
local map $\phi(a_0,a_1)=
a_0+a_1 \pmod{2}$.  Then $\bQ_1 \ = \ \bP_1 \,\symdif\,
\varsigma(\bP_1)$ and $\bQ_0 \ = \ \compl{\bQ_1} \ = \ 
\lb[\bP_0 \intsct
\varsigma(\bP_0)\rb]\, \disj \, \lb[\bP_1 \intsct 
\varsigma(\bP_1)\rb]$.  (See Figure \ref{fig:X:ledrappier}).

\item  \label{X:BLCA.phitor}
More generally, let $\Lat=\Zahl^D$ and let $\dB\subset\Lat$ be finite.
Let  $\sA=\Zahlmod{2}$, and suppose $\Phi$ has local map
$\phi:\sA^\Nh\into\sA$ defined:
$\phi(\ba) \ := \ \D \sum_{\nh\in\Nh} a_\nh \pmod{2}$,
for any $\ba\in\sA^\Nh$.
  Then  $\Phitor(\sP) = \{\bQ_0,\bQ_1\}$,
where $\D\bQ_1 \ = \ \D \Symdif_{\nh\in\Nh} \trans{\nh}(\bP_1)$, and $\bQ_0 \ = \ \Tor\setminus \bQ_1$.
  
\item \label{X:LCA.phitor} 
Let $n\in\Natur$ and let $\sA=\Zahlmod{n}$.
Let $\varphi_\nh\in\Zahlmod{n}$ be constants for all $\nh\in\Nh$.
Suppose $\Phi$ has local map $\phi:\sA^\Nh\into\sA$ defined
\ 
$\D  \phi(\ba) \ := \ \sum_{\nh\in\Nh} \varphi_\nh a_\nh \pmod{p}$,
for any $\ba\in\sA^\Nh$.
   Treat any  $\sP\in\MP$ as a function
 $\sP:\Tor\into\sA$.  Then 
$\D  \Phitor(\sP) \ = \  \sum_{\nh\in\Nh} \varphi_\nh\cdot\trans{\nh}(\sP)$.
}

  If $(\bX,d)$ is a metric space, recall that a function $\varphi:\bX\into\bX$ is
{\dfn Lipschitz} with constant $K>0$ if $\varphi$ is continuous, and furthermore,
for any $x,y\in\bX$,
\quad $d\lb(\varphi(x),\varphi(y)\maketall \rb) \ < \ K\cdot d(x,y)$.

\Theorem{\label{Phi.star.continuous}}
{
$(\MP,d_\symdif,\Phitor)$ is a topological dynamical system, and  
$\Phitor\colon\MP\goto\MP$ is $d_\symdif$-Lipschitz.
 
}
\bthmprf   
  Suppose $\Phi$ has local map $\phi:\sA^\Nh\into\sA$. 
Let  $B = |\Nh|$ and $A=|\sA|$; hence  $|\sA^\Nh| = A^B$.
We claim that $\Phitor$ has Lipschitz constant $2B\cdot A^B$. 
Suppose $\sP,\sP'\in\MP$, and $d_\symdif(\sP,\sP') < \eps$.
If $\Phitor(\sP)=\sQ$ and $\Phitor(\sP')=\sQ'$, then for all $a\in\sA$,
\[
  \bQ_a \quad = \quad 
\Union_{{{\bc\in\sA^\Nh}\atop{\phi(\bc)=a}}} \  \Intsct_{\nh\in\Nh}
\ \trans{-\nh} (\bP_{c_\nh}) 
\quad\AND\quad
  \bQ'_a \quad = \quad 
\Union_{{{\bc\in\sA^\Nh}\atop{\phi(\bc)=a}}} \  \Intsct_{\nh\in\Nh}
\ \trans{-\nh} (\bP'_{c_\nh}).
\]
Set $J:=A^B$ and
$K:=B$ in Lemma \ref{Phi.star.continuous.C3}(c)  to conclude 
$d_\symdif(\sQ, \ \sQ')  <  2B A^B\cdot \eps$. \nolinebreak
\ethmprf

\Proposition{\label{Xi.Phi.iso}}
{
Let $\ZP$ and $\xi_\trans{}:\ZP\into\QS$ be as in 
{\rm Proposition \ref{Xi.iso}}.
Then $\xi_\trans{}\circ \Phitor \ = \ \Phi \circ \xi_\trans{}$.
Thus, $\QS$ is a $\Phi$-invariant subset of $\sA^\dL$, and
$\xi_\trans{}$ is an  isomorphism from the topological dynamical system 
$(\ZP,d_\symdif,\Phitor)$ to the system $(\QS,d_B,\Phi)$.
}
\bthmprf
  Combine Lemma \ref{induced.map.on.partitions}(a) with Proposition
\ref{Xi.iso}(a,b).
\ethmprf  

A topological dynamical system  $(\bX,d,\varphi)$ is {\dfn equicontinuous} if,
for every $\eps>0$, there is $\del>0$ such that,
for any $x,y\in\bX$,
$\statement{ $d(x,y)<\del$}
\IMPLIES\statement{$d\lb(\varphi^n(x),\varphi^n(y)\rb)<\eps$ for all $n\in\Natur$}$.

\Proposition{\label{Phi.star.equicont}}
{
  If $(\sA^\Lat,d_C,\Phi)$ is equicontinuous, then 
$(\ZP,d_\symdif,\Phitor)$ is equicontinuous.
}
\bthmprf
 If  $(\sA^\Lat,d_C,\Phi)$ is equicontinuous, then 
Proposition 7 of \cite{BlanchardFormentiKurka} 
says  $(\sA^\Lat,d_B,\Phi)$ is equicontinuous.  Thus, the
subsystem $(\QS,d_B,\Phi)$  is also equicontinuous.  
Now apply  Proposition \ref{Xi.Phi.iso}.
\ethmprf

\section{\label{app:MP} 
The space of measurable partitions\protect\footnotemark}

\footnotetext{This section contains technical results which are used 
in \S\ref{S:boundary},  \S\ref{S:injective} and \S\ref{S:surject}.} 

The next result is used to prove Theorem  \ref{comeager.quasisturm.inject}.

\Proposition{\label{MP.complete}}
{
$\MP$ is complete and bounded in the $d_\symdif$ metric. 
}
\bthmprf  We'll embed $\MP$ as a closed subset of $\bL^2(\Tor,\Cplx)$,
so that the $d_\symdif$ metric is equivalent to the (complete) $\bL^2$ metric.

\Claim{\label{MP.complete.C1}
\hspace{-1.3em} If $\sP=\{\bP_a\}_{a\in\sA}$ and 
$\sQ=\{\bQ_a\}_{a\in\sA}$,  then  
$d_\symdif(\sP, \sQ)  =  2  \D\sum_{{a,b\in\sA}\atop{a\neq b}} 
\lam[\bP_a\intsct\bQ_b]$.}
\bclaimprf
For any $a\in\sA$,
\quad
$\lam[\bP_a\setminus\bQ_a] \ = \ \D \sum_{b\neq a} \ \lam[\bP_a\intsct\bQ_b]$.
Thus,
\[  \sum_{a\in\sA} \lam[\bP_a\setminus\bQ_a] 
\ = \quad \sum_{a\in\sA}  \ \sum_{b\neq a} \ \lam[\bP_a\intsct\bQ_b]
\ = \quad  \sum_{a\neq b\in\sA} 
\lam[\bP_a\intsct\bQ_b].
\]
Thus, $\D d_\symdif(\sP, \sQ) \ = \ 
\sum_{a\in\sA} \lam[\bP_a\setminus\bQ_a] +
\sum_{a\in\sA} \lam[\bQ_a\setminus\bP_a]
\ = \ 2  \sum_{a\neq b\in\sA} 
\lam[\bP_a\intsct\bQ_b]$.
\eclaimprf
  Suppose $|\sA| = A$, and identify $\sA$ with the $A$th roots of
unity in some arbitrary way:
\[
  \sA \quad\cong\quad \set{e^{\frac{2\pi k}{A}\bi}}{0\leq k < A}.
\]
  Any $\sP\in\MP$  defines a function in $\bL^2(\Tor;\Cplx)$ 
(also denoted by $\sP$), such that $\sP^{-1}\{a\} \ = \ \bP_a$
for any $a\in\sA$. 
We can then measure the distance between partitions in the
$\bL^2$ metric:
\ $\D
\norm{\sP-\sQ}{2}  =  
 \lb( \int_{\Tor} \lb|\sP(t)-\sQ(t)\rb|^2 \ d\lam[t] \rb)^{1/2}$. \ 
  Let $m =\D \min_{a\neq b\in\sA} |a-b|^2$ and
$M = \D \max_{a\neq b\in\sA} |a-b|^2$.

\Claim{\label{MP.complete.C2}
For any $\sP,\sQ\in\MP$, \quad
$\frac{m}{2}\cdot   d_\symdif (\sP,\sQ) \ < \
\norm{\sP-\sQ}{2}^2 \  < \  \frac{M}{2}\cdot  d_\symdif (\sP,\sQ)$.
}
\bclaimprf
$\D\norm{\sP-\sQ}{2}^2 \ = \ 
\int_\Tor \lb|\sP(\tort)-\sQ(\tort)\rb|^2 \ d\lam[\tort]
\ = \ 
\sum_{a\neq b\in\sA} \int_{\bP_a\intsct\bQ_b} |a-b|^2 \ d\lam$

$\D = \ \sum_{a\neq b\in\sA}  |a-b|^2 \cdot \lam\lb[\bP_a\intsct\bQ_b\rb]$.
\ Thus,
$\D\sum_{a\neq b\in\sA} m \cdot \lam\lb[\bP_a\intsct\bQ_b\rb]
\quad \leq \quad 
\norm{\sP-\sQ}{2}^2
\ \leq \ 
\sum_{a\neq b\in\sA} M \cdot \lam\lb[\bP_a\intsct\bQ_b\rb]$. \quad
Now apply Claim \ref{MP.complete.C1}.
\eclaimprf
  Thus, the $d_\symdif$ and $\bL^2$ metrics are equivalent, so
$\MP$ is bounded and complete in $d_\symdif$ if and only if $\MP$
is bounded and complete in $\bL^2$.  It remains to show:

\Claim{\label{MP.complete.C3} $\MP$ is a closed, bounded subset of
$\bL^2(\Tor,\Cplx)$.} 

 \bclaimprf $\MP$ is bounded because it is a subset of the unit ball
in $\bL^2$.  To see that $\MP$ is closed, suppose
$\{\sP_n\}_{n=1}^\oo\subset\MP$ was a sequence of $\sA$-labelled
partitions, and that $\D\bL^2\!-\!\lim_{n\goto\oo} \sP_n \ = \ \sP$, where
$\sP\in\bL^2(\Tor,\Cplx)$.  We must show that $\sP$ is also an
$\sA$-labelled partition ---in other words, that the essential image of
$\sP$ is $\sA\subset\Cplx$.

  Suppose not.  Then there is some $\eps>0$ and some subset
$\sU\subset\Cplx$ with $d(\sU,\sA)=\eps$ such that, if $\bU = \sP^{-1}(\sU)\subset\Tor$, then
$\lam[\bU]>0$.  But then for any $n\in\Natur$,
\beq
 \norm{\sP_n-\sP}{2}^2 
& = &\int_\Tor \lb|\sP_n-\sP\rb|^2 \ d\lam
\quad \geq \quad \int_\bU \lb|\sP_n-\sP\rb|^2 \ d\lam
\\& \geq & \int_\bU |\eps|^2 \ d\lam
\ = \  \eps^2 \cdot \lam[\bU].
\eeq
This contradicts the hypothesis that $\D\lim_{n\goto\oo} \norm{\sP_n-\sP}{2} \ = \ 0$.
\eclaimthmprf
{\em Remark.}  Although $(\MP,d_\symdif)$ is complete and bounded, it is
not compact.  For example, let $\sA=\{0,1\}$ and $\Tor=\CO{0,1}$.
Fix $n\in\Natur$, and for each $j\in\CC{1..2^n}$, let 
$\bI_n \ = \ \OO{\frac{j-1}{2^n}, \frac{j}{2^n}}$.  Now
define partition $\sP^{(n)}\ = \lb\{\bP_0^{(n)},\bP_1^{(n)}\rb\}$, where
$\D \bP_0^{(n)}  =   \Disj_{\mathrm {even} \ j=2}^{2^n} \bI_j$
and
 $\D \bP_1^{(n)}  =   \Disj_{\mathrm {odd} \ j=1}^{2^n-1} \bI_j$.
  It is easy to check that $d_\symdif(\sP^{(n)},\sP^{(m)}) \ = \ 1$
for any $n\neq m$.  Hence, $\{\sP^{(n)}\}_{n=1}^\oo$ is an infinite
sequence of partitions with no convergent subsequence, which would be
impossible if $\MP$ were compact.

\subsection{Dyadic partitions and the density of $\OP$ in $\MP$
\label{app:open.dyadic}}

Identify $\Tor = \CO{0,1}^K$ for some $K\in\Natur$.  For any $n>0$, an
{\dfn $n$-dyadic number} is a number $d = \frac{k}{2^n}$ for some
$k\in\CC{0..2^N}$. An {\dfn $n$-dyadic interval} is a closed interval
$\bI=\CC{d_1,d_2}$, where $d_1,d_2$ are $n$-dyadic numbers; an {\dfn
$n$-dyadic cube} in $\Tor$ is a set $\bC = \bI_1 \x \bI_2 \x \ldots \x
\bI_K$, where $\bI_1,\ldots,\bI_K$ are $n$-dyadic intervals.  An {\dfn $n$-dyadic set} is an open set $\bD\subset\Tor$ such that $\bD =
int\lb(\Disj_{j=1}^J \bC_J\rb)$ for some collection
$\{\bC_1,\ldots,\bC_J\}$ of $n$-dyadic cubes.
An ($\sA$-labelled) {\dfn $n$-dyadic partition} is an open partition
$\sD= \{\bD_a\}_{a\in\sA}$ where $\bD_a$ is $n$-dyadic for all $a\in\sA$.
 A {\dfn dyadic} number (cube, set,
partition, etc.) is one which is $n$-dyadic for some $n>0$.
Let $\DP\subset\OP$ be the set of all $\sA$-labelled dyadic partitions.

The next result is used to prove Lemma \ref{SOT.dense},
Lemma \ref{primitive.dense}, and Proposition \ref{CA.surjective.on.MP}.

\Proposition{\label{DP.dense.in.MP} \label{OP.dense.in.MP}}
{
 $\DP$ is $d_\symdif$-dense in $\MP$.
\quad
Thus,  $\OP$ is also $d_\symdif$-dense in $\MP$.

Furthermore $\MP$ is separable.
\qed
}

To prove Proposition \ref{DP.dense.in.MP}, we use:

\Lemma{\label{dyadic.set.props}}
{
 {\bf(a)} \  If $\bD$ is an $n$-dyadic set, then $\bD$ is also $N$-dyadic for any
$N>n$.

  {\bf(b)} \   The union or intersection of any two $N$-dyadic sets is $N$-dyadic.

  {\bf(c)} \   If $\bD\subset\Tor$ is an $N$-dyadic set, then $\Tor\setminus\barbD$
is also $N$-dyadic.

  {\bf(d)} \  Let $\bM\subset\Tor$ be a measurable set.  For any $\del>0$, 
there is some $N>0$ and some $N$-dyadic set $\tlbM\subset\Tor$ with 
$\lam(\tlbM\symdif\bM) < \del$.

  {\bf(e)} \  Let $\bD\subset\Tor$ be an $n$-dyadic set, and let $\bM\subset\bD$ be
a measurable subset. For any $\del>0$, 
there is some $N>n$ and some $N$-dyadic $\tlbM\subset\bD$ with 
$\lam(\tlbM\symdif\bM) < \del$.

  {\bf(f)} \  If $\bM$ is measurable and $\bD$ is an $N$-dyadic set, then
$\lam(\bM \intsct \barbD) = \lam(\bM \intsct \bD)$.
}
\bthmprf {\bf(a)}, {\bf(b)} and {\bf(c)} are immediate from the definition.

{\bf(d):} \quad Recall \cite[Thm. 2.40(c), p.68, for example]{Folland}
 that $\bM$ can be $\del/2$-approximated by a
finite union of open cubes $\bK_1,\ldots,\bK_J$ for some $J>0$. Next, for
each $j\in\CC{1..J}$, there is some $N_j$ such that $\bK_j$ can be
$(\del/2J)$-approximated by an $N_j$-dyadic set $\tlbK_j$.  Let
$N=\max\{N_1,\ldots,N_J\}$.  Then {\bf(a)} says that
$\tlbK_1,\ldots\tlbK_J$ are all $N$-dyadic sets.  If $\bD =
\Union_{j=1}^J \tlbK_j$, then $\bD$ is also an $N$-dyadic set (by
{\bf(b)}), and $\bM \ \closeto{\del/2} \ \Union_{j=1}^J \bK_J \
\closeto{\del/2} \ \bD$.

  {\bf(e):}\quad First use part {\bf(d)} to find some $N$-dyadic set
$\hbM\subset\Tor$ with $\lam(\hbM\symdif\bM)<\eps$.  Now let
$\tlbM = \hbM\intsct\bD$.  Then $\tlbM$ is also an $N$-dyadic set
(by {\bf(a)} and {\bf(b)}), and since $\bM\subset\bD$, we have
$\lam(\tlbM\symdif\bM) \leq \lam(\hbM\symdif\bM)<\del$.  

 {\bf(f):}\quad $\partial\bD$ is a finite union of $(K-1)$-dimensional
hyperfaces, so $\lam(\partial\bD) \ = \ 0$.
\ethmprf

\bthmprf[Proof of Proposition \ref{DP.dense.in.MP}]
Let $\sP= \{\bP_a\}_{a\in\sA}$ be a measurable partition,
and let $\eps>0$.  We want a dyadic partition
$\sD= \{\bD_a\}_{a\in\sA}$ such that $d_\symdif(\sP,\sD)<\eps$.
  For simplicity, let $\sA:=\{1,2,\ldots,A\}$ for some
$A\in\Natur$.

Fix $\del>0$. Lemma \ref{dyadic.set.props}{\bf(d)} yields
$N_1>0$ and an $N_1$-dyadic set $\bD_1 \subset\Tor$ such that
$\lam(\bP_1 \symdif \bD_1) \ < \ \del$.
  Lemma \ref{dyadic.set.props}{\bf(c)} says that $\Tor' =
\Tor\setminus\barbD_1$ is also $N_1$-dyadic.  For all
$a\in\CC{3..A}$, let $\bP'_a \ = \ \bP_a \intsct \Tor'$, while $\bP'_2
\ = \ (\bP_1 \intsct \Tor')\disj (\bP_2 \intsct \Tor')$.

\Claim{\label{OP.dense.in.MP.C2}
$\lam(\bP_a\symdif \bP'_a) \ < \ \del$ for any $a\in\CC{3..A}$,
while $\lam(\bP'_2\symdif\bP_2) \ < \ 2\del$.}
\bclaimprf
  $\bP'_a\subset\bP_a$ and $\bP_a\setminus\bP'_a \ = \ \bP_a \intsct
\barbD_1$, so that  
\ $\lam(\bP_a\symdif \bP'_a) \  = \ 
 \lam(\bP_a\setminus\bP'_a) \  = \  \lam(\bP_a \intsct \barbD_1)
\  \eeequals{(*)} \  \lam(\bP_a \intsct \bD_1)$,
\ where $(*)$ is by Claim 
\ref{dyadic.set.props}{\bf(f)}. But $\bP_a \intsct \bD_1 \ \subset \ \bD_1 \setminus \bP_1$, so  $\lam(\bP_a \intsct \bD_1)  \leq  \lam(\bD_1 \setminus \bP_1)
  \leq  \lam(\bD_1 \symdif \bP_1)  =  \del$.  The proof for $\bP_2$ is similar.
\eclaimprf

  Lemma \ref{dyadic.set.props}{\bf(e)} yields some $N_2\geq N_1$
and an $N_2$-dyadic $\bD_2\subset\Tor'$
such that $\lam(\bP'_2 \symdif \bD_2) \ < \  \del$.
Hence,  \
$\lam(\bP_2 \symdif \bD_2) \ \leq \quad
 \lam(\bP_2\symdif\bP'_2) + \lam(\bP'_2 \symdif \bD_2)
\quad < \quad 2\del + \del \ = \quad 3\del$.

  Now, Lemma \ref{dyadic.set.props}{\bf(c)} says
 $\Tor'' = \Tor'\setminus\barbD_2$ is an $N_2$-dyadic set. 
For all $a\in\CC{4..A}$,
let $\bP''_a \ = \ \bP_a \intsct \Tor''$; hence $\bP''_a\subset\bP'_a$
and $\lam(\bP_a\setminus\bP'_a) \ < \ \del$, as in Claim  
\ref{OP.dense.in.MP.C2}.
 Also, let $\bP''_3 \ = \ (\bP'_2 \intsct \Tor'')\disj (\bP'_3\intsct \Tor'')$.
Then $\lam(\bP'_3\symdif\bP''_3) \ < \ 2\del$ as in Claim
 \ref{OP.dense.in.MP.C2}.

  Proceeding inductively, we obtain the triangle shown below:
\[
\begin{array}{ccccccccccccccccccc}
\bP_1 &\closeto{\del} & \bD_1 \\
\bP_2 &\closeto{2\del} & \bP'_2  &\closeto{\del} & \bD_2 \\
\bP_3 &\closeto{\del} & \bP'_3 &\closeto{2\del} & \bP''_3  &\closeto{\del} & \bD_3 \\
\bP_4 &\closeto{\del} & \bP'_4  &\closeto{\del} & \bP''_4  &\closeto{2\del} & \bP'''_4  &\closeto{\del} & \bD_4 \\

\vdots &\vdots& \vdots&\vdots & \vdots & \vdots&  \vdots & \vdots& \ddots & \ddots & \ddots  \\
\bP_A &\closeto{\del} & \bP'_A  &\closeto{\del} & \bP''_A  &\closeto{\del} & \bP'''_A  &\closeto{\del} &  \ldots &  \closeto{2\del} &\bP''''''_A &\closeto{\del} & \bD_A 
\end{array}
\]
  Hence, $\bP_a \closeto{(a+1)\del} \bD_a$ for any $a\geq 2$.  
  Now, $\bD_1,\ldots,\bD_A$ are disjoint dyadic sets, and
\[
d_\symdif(\sP,\sD) \quad<\quad
 \lam(\bP_1\symdif\bD_1) \ + \ \sum_{a=2}^A \lam(\bP_a\symdif\bD_a) 
\quad<\quad 
\del \ + \  \sum_{a=2}^A (a+1)\del 
\quad=\quad M\cdot \del,
\]
 where $M := \frac{A(A+1)}{2} - 2$.  So, choose $\del < \eps/M$.
\ethmprf

\subsection{Simple Partitions and the Injectivity of $\Splat$
\label{app:simple.inject}}

  If $\sP\in\MP$, and $\tors\in\Tor$, then $\tors$ is a
(translational) {\dfn symmetry} of $\sP$ if $\rot{\tors}(\sP) = \sP$
\lamae.  If $\sP\in\OP$,  then we also require that
 $\rot{\tors}(\sP) = \sP$ \topae.
 The (translational) symmetries of $\sP$ form a
closed subgroup of $\Tor$; if this group is trivial, we say $\sP$ is
{\dfn simple}.

\example{
 Identify $\Torus{1} \cong \CO{0,1}$.
Let $0<\bet_1<\bet_2<\cdots<\bet_N<1$ be irrational numbers.  Then
 $\sP \ \colon= \ \lb\{ \CO{0,\bet_1}, \ \CO{\bet_1,\bet_2}, \ \ldots, \
\CO{\bet_N,1}\rb\}$ is a simple partition of $\Torus{1}$.

To see this, suppose $\tors\in\Torus{1}$ is a symmetry.  Then
$\rot{\tors}\CO{0,\bet_1} \ = \ \CO{\bet_n,\bet_{n+1}}$ for some
$n\in\CO{1...N}$.  Hence, $\tors=\bet_n$.  But we can assume WOLOG that
$\tors$ is a rational number (see Lemma \ref{C:symmetry.rational}) 
Hence, either $\bet_n$ is rational (a contradiction) or $\tors=0$.}

 The next result is used to prove Theorem
\ref{comeager.quasisturm.inject}(a) and  Lemma
\ref{splat.quasiperiodic.factor.map}.

\Lemma{\label{simple.partn.inj}}
{ 
If $\sP\in\MP$ is simple, 
then the map $\Splat:\tlTor\into\sA^\Lat$ is injective \lamae.
}
\bthmprf
  We must show that, for $\forall_\lam \ \tort_1,\tort_2\in\Tor$,
 if $\tort_2\neq\tort_1$ then $\Splat(\tort_2)\neq \Splat(\tort_1)$.
Let $\tors=\tort_2-\tort_1$, and let $\bU=\set{\tort\in\Tor}{\sP(\tort)
\neq \sP(\tort+\tors)}$.   Since $\sP$ has trivial symmetry group, it follows
that $\lam[\bU]>0$.  Let $\tlbU=\D\Union_{\ell\in\Lat}\trans{\ell}(\bU)$;
 then $\lam[\tlbU]>0$ and $\tlbU$ is $\varsigma$-invariant.  Since $\varsigma$ is ergodic,
it follows that $\lam[\tlbU]=1$.  Thus, generically, $\tort_1\in\tlbU$,
which means that $\tort_1\in\trans{-\ell}(\bU)$ for some $\ell\in\Lat$.  This
means that $\trans{\ell}(\tort_1)\in\bU$.  Thus,
\[
\Splat(\tort_1)_\ell 
 \ = \ 
\sP\lb(\trans{\ell}(\tort_1)\rb) \  \neq \ 
 \sP\lb(\trans{\ell}(\tort_1) + \tors\rb) 
\ =\ 
 \sP\lb(\trans{\ell}(\tort_1 + \tors)\rb) 
\ =\  \sP\lb(\trans{\ell}(\tort_2)\rb)
 \ = \ \Splat(\tort_2)_\ell.
\]
Hence $\Splat(\tort_1) \ \neq \ \Splat(\tort_2)$. 
\ethmprf

In Proposition \ref{expansive.implies.chopping}, we'll need to replace a
 {\em nonsimple} partition with a simple `quotient' partition.

\Lemma{\label{quotient.partition}}
{
If $\sP\in\OP$ is not simple, let
$\bS\subset\Tor$ be the translational symmetry group of $\sP$.

\breath

\  {\bf(a)} \   $\quoTor=\Tor/\bS$ is also a torus (possibly of lower dimension),
and the quotient homomorphism $q:\Tor\into\quoTor$ is continuous.

\  {\bf(b)} \    There is a unique $\sA$-valued, open partition $\quosP\in\quoOP$
such that $\quosP\circ q \ = \ \sP$.  

\  {\bf(c)} \  There is a natural
$\Lat$-action on $\quoTor$ ---denoted $\quotrans{}$ ---such that
for all $ \ell\in\Lat$, 
\quad $\quotrans{\ell}\circ q \ = \ q\circ \trans{\ell}$.   

\breath

 Define $\tlquoTor$ and 
$\quoSplat:\tlquoTor\into\sA^\Lat$ in the obvious way.  Then:

\breath

\  {\bf(d)} \   If $\tort\in\tlTor$, and $\quotort = q(\tort)$, then
$\quotort\in\tlquoTor$, and
$\quoSplat(\quotort) \ = \ \Splat(\tort)$.

\  {\bf(e)} \   Hence, $\quoSplat(\tlquoTor) \ = \ \Splat(\tlTor)$,
\ $\Qshift[\quotrans{}]{\quosP} \ = \ \Qshift{\sP}$,
\ and $\QM[\quotrans{}]{\quosP} \ = \ \QM{\sP}$.  

\  {\bf(f)} \  $\quoSplat:\tlquoTor\into\sA^\Lat$ is injective \lamae.
}  
\bthmprf {\bf(a)} is because $\bS$ is a closed subgroup of $\Tor$.\quad
{\bf(b)}  is because $\bS$ is a group of symmetries of $\sP$ (resp. $\sQ$).\quad
To see {\bf(c)}, suppose $\tau:\Lat\into\Tor$ is the homomorphism
such that $\trans{\ell}(\tort) \ = \ \tort+\tau(\ell)$ for all $\tort\in\Tor$
and $\ell\in\Lat$.  Define $\quotau \ = \ q\circ\tau:\Lat\into\quoTor$, and
then define  $\quotrans{\ell}(\quotort) \ = \ \quotort+\quotau(\ell)$
 for all $\quotort\in\quoTor$ and $\ell\in\Lat$.
\quad
 {\bf(d)} follows from the defining properties
of $\quosP$ in {\bf(b)} and $\quotrans{}$ in {\bf(c)},
and {\bf(e)} follows immediately from {\bf(d)}.
\quad To see {\bf(f)}, 
observe that the symmetry group of $\quosP$ is $q(\bS)=\{0\}$.
Now apply Lemma \ref{simple.partn.inj}.
\ethmprf

If $\sP$ is an {\em open} partition, we can
strengthen Lemma \ref{simple.partn.inj} to get a homeomorphism;
this is used to prove Proposition \ref{expansive.implies.chopping} and
 Theorem \ref{comeager.quasisturm.inject}(b).

\newcommand{\Pinv}{\P}

\Proposition{\label{open.simple.partn.inj}}
{
Let  $\sP\in\OP$ be a simple open partition. Let
$\tlgP=\Splat(\tlTor)\subset\sA^\Lat$.  Then:
\bthmlist
\item  $\Splat:\tlTor\into\tlgP$ is a uniform homeomorphism with respect to the $d_C$ metric on $\tlgP$.

\item  $\Splat:\tlTor\into\tlgP$ is a uniform homeomorphism with respect to the $d_B$ metric on $\tlgP$.

\item  Thus, the Cantor topology and the Besicovitch topology agree on $\tlgP$.

\item  Let $\gP=\Qshift{\sP}$; then there is a continuous surjection $\Pinv:\gP\into\Tor$
such that:
\bdesc
  \item[\quad {[i]}] For all $\ell\in\Lat$, \ 
$\Pinv\circ\shift{\ell} \ =  \ \trans{\ell}\circ\Pinv$.

  \item[\quad {[ii]}]
 If $\tort\in\tlTor$, then $\Pinv\circ\Splat(\tort) \ = \ \tort$.
 If $\bp\in\tlgP$, then $\Splat\circ\Pinv(\bp) \ = \ \bp$.
\edesc
\ethmlist
}
\bthmprf 
{\bf(a)} \ $\Splat$ is continuous by Proposition \ref{splat.homeo}(b).
To show $\Splat$ is  injective and uniformly open,
fix $\eps>0$.  We claim there is some $M\in\Natur$
such that
\beqn
\label{open.simple.partn.inj.e1}
\mbox{For any $\tort,\toru\in\tlTor$,} \quad 
\statement{$\Splat(\tort)\restr{\dB(M)} \ = \ 
\Splat(\toru)\restr{\dB(M)}$}
\IMPLIES
\statement{$d(\tort,\toru)<\eps$}.
\eeqn
Suppose not; then  for any $n\in\Natur$, there  are $\tort_n,\toru_n\in\tlTor$
such that $\Splat(\tort_n)\restr{\dB(n)}  =  
\Splat(\toru_n)\restr{\dB(n)}$, but
$d(\tort_n,\toru_n)>\eps$.  Let $\tors_n=\tort_n-\toru_n$.
By dropping to a subsequence if necessary,
we can assume that the sequence $\{\tors_n\}_{n=1}^\oo$ converges
to some $\tors\in\Tor$, with $d(\tors,0)>\eps$ (because $\Tor$ is compact).  
We'll show that $\tors$ is a symmetry of $\sP$ (thereby contradicting
simplicity).

\Claim{
For any small enough $\del>0$, 
there is a $\del$-dense subset $\bU_\del\subset \Tor$, 
such that $\sP(\rot{\tors}(\toru))  =  \sP(\toru)$ for all
$\fu\in\bU_\del$.}
\bclaimprf
 For any $\tort\in\Tor$ and all $\ell\in\Lat$ 
let $\tort^\ell := \trans{\ell}(\tort)$.  The torus rotation system
$(\Tor,d,\trans{})$ is minimal and isometric, so
 if $n$ is large enough, then, for any $\tort\in\Tor$, the
set $\{\tort^\fb\}_{\fb\in\dB(n)}$ is $\del$-dense in
$\Tor$. 

Recall (\S\ref{S:prelim}) that  
$\bP^*:=\Disj_{a\in\sA} \bP_a$ is open and dense in $\Tor$;
thus $\partial\sP:=\Union_{a\in\sA} \partial\bP_a$ 
is nowhere dense.  Thus, $\rot{-\tors}(\partial\sP)$ is 
nowhere dense (since $\rot{\tors}$ is an isometry).
So, if $\del$ is small enough, the set 
$\bU_\del  :=  
\set{\toru_n^\fb}{\fb\in\dB(n)
\AND d(\rot{\tors}(\toru_n^\fb),\partial\sP)\geq\del}$
is $\del$-dense in $\Tor$.

If $n$ is large enough, then $\tors_n \ \closeto{\del} \ \tors$.  Thus,
$\tort_n \ = \  \toru_n+\tors_n \ \closeto{\del}\  \toru_n+\tors 
\ = \ \rot{\tors}(\toru_n)$;
thus, for all $\fb\in\dB(n)$, \ $\tort^\fb_n  \ \closeto{\del} \ 
\rot{\tors}(\toru_n^\fb)$.  Thus, for any $\toru_n^\fb\in\bU_\del$, 
\
$\sP(\rot{\tors}(\toru_n^\fb)) \ \eeequals{(*)} \ \sP(\tort_n^\fb) 
  \ \eeequals{(\dagger)} \ \sP(\toru_n^\fb)$.

 $(*)$ is because $\tort^\fb_n \ \closeto{\del} \ \rot{\tors}(\toru_n^\fb)
\ \not\!\!\!\!\closeto{\ \del} \ \partial\sP$;\quad
$(\dagger)$ is because $\Splat(\tort_n)\restr{\dB(n)} =
\Splat(\toru_n)\restr{\dB(n)}$.
\eclaimprf
  Now, let $\{\del_n\}_{n=1}^\oo$ be a sequence tending to zero,
and for each $n\in\Natur$, let $\bU_{\del_n}$ be as in Claim 1.
Let $\bU := \Union_{n=1}^\oo \bU_{\del_n}$;  then $\bU$ is a dense
subset of $\Tor$ so that 
$\sP(\rot{\tors}(\toru)) \ = \ \sP(\toru)$ for all $\fu\in\bU$.
In other words, for each $a\in\sA$, \ $\rot{\tors}(\bP_a)\intsct\bU
\ = \ \bP_a\intsct\bU$, which means that 
$\rot{\tors}(\bP_a)\intsct\bU \intsct  \compl{\barbP_a}
\ = \ \emptyset$,
where $\compl{\barbP_a}:= \Tor\setminus\barbP_a$.
But $\D \rot{\tors}(\bP_a) \intsct \compl{\barbP_a}$
is an open subset of $\Tor$, so if it is disjoint from the dense set $\bU$,
 it
must be empty.  But if
$\rot{\tors}(\bP_a) \intsct \compl{\barbP_a} \ = \ \emptyset$,
then $\rot{\tors}(\bP_a) \subseteq\barbP_a$.
By symmetric reasoning, $\rot{\tors}(\barbP_a) \supseteq\bP_a$;
hence $\rot{\tors}(\bP_a) =\bP_a$ \topae.  This holds for all $a\in\sA$,
so we conclude that $\rot{\tors}(\sP)=\sP$ \topae.  Thus,
$\tors$ is a symmetry of $\sP$, contradicting the simplicity of $\sP$.

{\bf(b)} \quad 
 We claim that, for any
$\eps>0$, there is some $\del>0$
such that, for any $\tort,\toru\in\tlTor$,
\[
\statement{$\Splat(\tort)\restr{\dM}  =  
\Splat(\toru)\restr{\dM}$ for some $\dM\subset\Lat$ with $\density{\dM}>1-\del$}
\IMPLIES
\statement{$d(\tort,\toru)<\eps$}.
\]
 The proof  the same as {\bf(a):}
replace $\dB(n)$ with $\dM$, prove
the appropriate version of Claim 1, deduce a symmetry, 
and derive a contradiction.

{\bf(c)} \quad Follows from {\bf(a)} and {\bf(b)}.

{\bf(d)}\quad 
Recall that $\gP$ is the $d_C$-closure of $\tlgP=\Splat(\tlTor)$.  Thus,
if $\bp\in\gP$, then there exists a sequence $\{\tort_n\}_{n=1}^\oo\subset\tlTor$ 
such that $\bp = \D\dClim_{n\goto\oo} \ \Splat(\tort_n)$.  
We can drop to a subsequence such that $\{\tort_n\}_{n=1}^\oo$
converges to some $\tort\in\Tor$ (because $\Tor$ is compact).  
We define $\Pinv(\bp):=\tort$. 

{\em $\Pinv(\bp)$ is well-defined:} 
 Suppose $\{\tort'_n\}_{n=1}^\oo\subset\tlTor$
was another sequence with $\bp = \D\dClim_{n\goto\oo} \ \Splat(\tort'_n)$;
we claim $\D \lim_{n\goto\oo} \tort'_n =  \lim_{n\goto\oo} \tort_n$.
Let $\eps>0$, and let $M>0$ be as in 
assertion (\ref{open.simple.partn.inj.e1}).  
Find $N$ large enough that, if $n>N$, then  
$\Splat(\tort'_n)\restr{\dB(M)} \ = \ 
\bp\restr{\dB(M)} \ = \ \Splat(\tort_n)\restr{\dB(M)}$.
 Then assertion (\ref{open.simple.partn.inj.e1})
says that $\tort'_n \ \closeto{\eps} \ \tort_n$.
Since $\eps$ is arbitrary, we conclude that
 $\D \lim_{n\goto\oo} \tort'_n =  \lim_{n\goto\oo} \tort_n$.

{\em Continuous:}\quad Fix $\eps>0$. Let $M$ be as in 
 assertion (\ref{open.simple.partn.inj.e1});
and suppose $\bp\restr{\dB(M)} = \bp'\restr{\dB(M)}$. If
$\{\tort_n\}_{n=1}^\oo\subset\tlTor$ 
is such that $\D \bp = \D\dClim_{n\goto\oo} \ \Splat(\tort_n)$ and
$\{\tort'_n\}_{n=1}^\oo\subset\tlTor$ 
is such that $\D \bp' = \D\dClim_{n\goto\oo} \ \Splat(\tort'_n)$,
then  assertion (\ref{open.simple.partn.inj.e1}) says that, 
for large  $n$, we must have $\tort_n \closeto{\eps} \tort'_n$.
Hence $\Pinv(\bp) \ = \ \D\dClim_{n\goto\oo} \tort_n \ \closeto{\eps}  \
\D\dClim_{n\goto\oo} \tort'_n  \ = \ \Pinv(\bp')$.

{\em Surjection:}\quad Let $\tort\in\Tor$; find a sequence
$\{\tort_n\}_{n=1}^\oo\subset\tlTor$ such that $\D \tort=\lim_{n\goto\oo}
\tort_n$.  Let $\bp = \D\dClim_{n\goto\oo} \Splat(\tort_n)$;  then
by definition, $\tort=\Pinv(\bp)$.

{\bf[i]:}\quad If 
$\bp = \D\dClim_{n\goto\oo} \ \Splat(\tort_n)$, then
\[
\shift{\ell}(\bp) \quad = \quad 
 \D\dClim_{n\goto\oo} \ \shift{\ell}(\Splat(\tort_n))  
\quad \eeequals{(*)} \quad
  \D\dClim_{n\goto\oo} \ \Splat(\trans{\ell}(\tort_n)),
\]
where $(*)$ is Proposition \ref{splat.homeo}(a).
Hence $\D\Pinv( \shift{\ell}(\bp))
 =   \lim_{n\goto\oo} \ \trans{\ell}(\tort_n)
 =    \trans{\ell}( \lim_{n\goto\oo} \tort_n)
 =   \trans{\ell}(\Pinv(\bp))$.

{\bf[ii]:} \    Let $\tort\in\tlTor$ and
$\bp=\Splat(\tort)$. Let $\tort_n=\tort$, \ $\forall \, n\in\Natur$;
then $\D\Pinv(\bp) = \lim_{n\goto\oo}\tort_n  =  \tort$.\nolinebreak
\ethmprf
  Thus, in a sense, the topological $\Lat$-system $(\gP,d_C,\shift{})$ is
`almost' isomorphic to the system $(\Tor,d,\trans{})$.  The only caveat
is that the function $\Pinv$ is many-to-one on the elements of
$\Tor\setminus\tlTor$.

\section{Boundary growth \& Chopping  \label{S:boundary}}

Suppose $\Tor=\Torus{1}$, and $\sP$ is an open partition of $\Tor$
such that each element of $\sP$ is a finite collection of open
intervals; \ then $\partial\sP \ := \ \Union_{a\in\sA}\ \partial \bP_a$
is a finite set of points in $\Tor$.  Let $\Phi$ be a cellular
automaton.  Hof and Knill \cite{HofKnill} observed empirically that
$\D\#\lb(\maketall \partial \lb( \Phitor^n[\sP]\rb) \rb)$ grows
polynomially like $n^\alp$ as $n\goto\oo$, for some exponent $\alp\leq D$
(where $\Lat=\Zahl^D)$.  They asked: is $\#\lb(\maketall \partial \lb(
\Phitor^n[\sP]\rb)\rb)$ really growing?  Is the growth polynomial?  What's the exact value of $\alp$?

  If $\#\lb(\maketall \partial \lb( \Phitor^n[\sP]\rb)\rb)$ gets large
as $n\goto\oo$, then each cell of $\Phitor^n[\sP]$ is `chopped' into
many tiny separate intervals; \  we say $\Phitor$ is {\dfn chopping}
(we'll make this precise later).  In this section, we investigate
chopping, answer Hof and Knill's questions, and generalize these ideas to
$\Torus{K}$.

\paragraph*{\sc Partition boundary size in $\Torus{K}$:}
 Suppose $\Tor=\Torus{K}$ for $K\geq 1$, and let $\sP\in\OP$.
 To characterize the growth
of $\partial\lb( \Phitor^n(\sP)\rb)$, we first need a way to measure its size.
Let $\gC = \set{\bC\subset\Tor}{\mbox{$\bC$ closed}}$, and
let $\setsize{\bullet}:\gC\into\CC{0,\oo}$ be some `pseudomeasure',
satisfying:
\bdesc
  \item[(M1)] {\em Monotonicity:} \ If $\bC_1\subset\bC_2$, then $\setsize{\bC_1}
\leq \setsize{\bC_2}$.

  \item[(M2)] {\em Additivity:}\  $\setsize{\bC_1 \disj \bC_2} = \setsize{\bC_1} + \setsize{\bC_2}$.

  \item[(M3)] {\em Translation Invariance:} \ For any $\tort\in\Tor$,\quad
  $\setsize{\rot{\tort}(\bC)} = \setsize{\bC}$.

  \item[(M4)] {\em Nontriviality:} \ $0 < \setsize{\partial\sP} < \oo$,
and  $0<\setsize{\partial\lb(\Phitor^n(\sP)\rb)} <\oo$ for all $n\in\Natur$.
\edesc
 If $\Tor=\Torus{1}$, then $\partial\sP$ is usually a discrete
subset of $\Tor$, and the obvious function satisfying {\bf(M1)}-{\bf(M4)}
is $\setsize{\bC} = \#(\bC)$ (modulo multiplication by some constant).
However, if $\Tor=\Torus{K}$ for $K\geq 2$, then condition {\bf(M4)}
makes the choice of $\setsize{\bullet}$ dependent on the geometry of
$\sP$:
\bitem
  \item If  $\partial \sP$ is a union of
piecewise smooth $(K-1)$-dimensional submanifolds of $\Torus{K}$,
then let $\setsize[*]{\bullet}$ be the $(K\!-\!1)$-dimensional Lebesgue
measure.  For example, if $K=1$, $2$, or $3$, then 
$\setsize[*]{\bullet}$ measures cardinality, 
length, or surface area, respectively.

  \item  If $(K-1)\leq \kap < K$, and $\partial\sP$
has Hausdorff dimension $\kap$, then
let $\setsize[\kap]{\bullet}$ be the $\kap$-dimensional Hausdorff measure.
If $\kap = (K-1)$, then $\setsize[\kap]{\bC} = \setsize[*]{\bC}$.

  \item For any 
$\bC\subset\Tor$ and $\eps>0$, let $\Ball(\bC,\eps) \ = \
\set{\tort\in\Tor}{d(\tort,\torc)<\eps \mbox{ \ for some \ }
\torc\in\bC}$. 
Define 
$\D  \setsize[L]{\bC} \ = \   \lim_{\eps\goto 0} \
\frac{1}{2\eps} \lam\lb(\maketall \Ball(\bC,\eps)\rb)$.
\   If $\bC$ is a $(K\!-\!1)$-dimensional submanifold, then $\setsize[L]{\bC}
=\setsize[*]{\bC}$.  We call $\setsize[L]{\bullet}$ the {\dfn Lipschitz}
 pseudomeasure because of Proposition \ref{splat.is.lipschitz} below.
\eitem

  We say $\Phitor$ {\dfn chops} $\sP$ {\dfn on average} if
$\D  \lim_{N\goto\oo} \frac{1}{N} \sum_{n=1}^N \setsize{\partial\lb(\Phitor^n[\sP]\rb)} 
 =  \oo$.
  Equivalently, there is a subset
$\dJ\subset\Natur$ of \Cesaro density 1 such that
$\D \lim_{\dJ\ni j\goto\oo} \ \setsize{\partial\lb(\Phitor^n[\sP]\rb)} \ = \ \oo$.

  We say  $\Phitor$ {\dfn chops} $\sP$ {\dfn intermittently} if
$\D \limsup_{n\goto\oo} \setsize{\partial\lb(\Phitor^n[\sP]\rb)} \ = \ \oo$.
Equivalently, there is a (possibly zero-density) subset
$\dJ\subset\Natur$ such that
$\D \lim_{\dJ\ni j\goto\oo} \ \setsize{\partial\lb(\Phitor^n[\sP]\rb)} \ = \ \oo$.
Clearly, if $\Phitor$ chops $\sP$ on average, then it does so
intermittently (but not conversely).

  Note that these definitions depend upon the 
pseudomeasure $\setsize{\bullet}$.  For a fixed choice of 
$\setsize{\bullet}$, we say $\Phitor$ is {\dfn
$\setsize{\bullet}$-chopping on average} (resp. {\dfn intermittently})
if, for any $\sP\in\OP$ with $0 < \setsize{\partial\sP} <
\oo$, \quad $\Phitor$ chops $\sP$ on average (resp. intermittently)
with respect to $\setsize{\bullet}$.

Whenever $\Phitor$ chops $\sP$,
the growth rate of  $\setsize{\partial\lb(\Phitor^n[\sP]\rb)}$ must be
(sub)polynomial:

\Proposition{\label{boundary.grows.polynomially}}
{
  Let $\Lat=\Zahl^D$. Let 
$\Phi:\sA^\Lat\into\sA^\Lat$ be a CA.  There is a
constant $C>0$ such that, if $\sP\in\OP$ and $n\in\Natur$, then
$\setsize{\maketall \partial \lb( \Phitor^n[\sP]\rb)}
\ \leq \  C\cdot n^D\cdot  \setsize{\partial \sP}$.
}
\bthmprf 
  Suppose $\Phi$ has local map $\phi:\sA^\Nh\into\sA$ for some finite
$\Nh\subset\Lat$.  It follows from eqn.(\ref{partition.image}) that
$\partial \lb(\Phitor(\sP)\rb) \subset \D \Union_{\nh\in\Nh}
\trans{\nh}(\partial \sP)$.  Hence, if $B = \#(\Nh)$, then
\beqn
\label{boundary.grows.polynomially.e1}
\setsize{\maketall \partial \lb(\Phitor[\sP] \rb)} 
\ \ \leeeq{(M1)}  \ \
  \setsize{\Union_{\nh\in\Nh} \trans{\nh}(\partial \sP)} 
\ \ \leeeq{(M2)} \ \
 \sum_{\nh\in\Nh} \setsize{ \trans{\nh}(\partial \sP) } 
\ \ \eeequals{(M3)} \ \ 
\sum_{\nh\in\Nh} \setsize{ \partial \sP }
 \ \ = \ \ B\cdot \setsize{\partial \sP}.
\eeqn
  Let $\Nh_n  =  \set{\nh_1 + \cdots + \nh_n}{\nh_1,...,\nh_n\in\Nh}$.
Then $\Phi^n$ has a local map $\phi^{(n)}\colon\sA^{\Nh_n}\rightarrow\sA$; \ hence,
by reasoning similar to (\ref{boundary.grows.polynomially.e1}), we have:
\beqn
\label{boundary.grows.polynomially.e2}
\setsize{\maketall \partial \lb(
\Phitor^n[\sP] \rb)} \quad \leq \quad 
 B_n\cdot \setsize{ \partial \sP },
\eeqn
where $B_n \ = \  \#(\Nh_n)$.  Find $R>0$ such that
 $\Nh\subseteq\CC{-R...R}^D$.  Then $\Nh_n\subseteq\CC{-nR\ldots nR}^D$, so
\beqn
\label{boundary.grows.polynomially.e3}
B_n \quad \leq \quad \#\lb(\CC{-nR\ldots nR}^D\rb) \quad = \quad C\cdot n^D,
\eeqn
 where $C=(2R+1)^D$. Combine (\ref{boundary.grows.polynomially.e2}) and
(\ref{boundary.grows.polynomially.e3}) to get: 
$\setsize{\maketall \partial \lb( \Phitor^n[\sP] \rb)}
\ \leq \ C n^D \cdot \setsize{\partial \sP}$.\nolinebreak
\ethmprf

\paragraph*{\sc Chopping in Boolean Linear Cellular Automata:}
   Let $\sA=\Zahlmod{2}=\{0,1\}$.   A CA with
local map $\phi:\sA^\Nh\into\sA$ is a 
{\dfn boolean linear cellular automaton} (BLCA) if 
\beqn
\label{BLCA2}
  \phi(\ba) \ = \ \sum_{\nh\in\Nh} a_\nh \pmod{2},
\quad \mbox{for any $\ba\in\sA^\Nh$.}
\eeqn
 We assume that $\Phi$ is `nontrivial'
in the sense that $\#(\dB)>1$.
Thus, Examples \ref{X:ledrappier} and \ref{X:BLCA.phitor} were BLCA.
We'll show that BLCA are chopping on
average.  Then we'll characterize the asymptotic growth rate of
$\setsize{\maketall \partial \lb(\Phitor^j[\sP]\rb)}$ in a special case.

\Proposition{\label{boolean.LCA.on.partn}}
{
  Let  $\sP = \{\bP_0,\bP_1\}$ be an $\sA$-indexed open 
partition of $\Tor$. Then
\[
  {\bf(a)} \ \ \partial(\sP) \ = \ \partial(\bP_1)
 \quad\AND\quad{\bf(b)} \ \
\Symdif_{\nh\in\Nh} \trans{\nh}(\partial \sP)
\quad\subseteq\quad
\partial \lb(\maketall \Phitor(\sP) \rb)
\quad\subseteq\quad
\Union_{\nh\in\Nh} \trans{\nh}(\partial \sP).
\]
}
\bthmprf  {\bf(a)} is immediate.  For {\bf(b)}, recall from  
Example \ref{X:BLCA.phitor} that 
$\Phitor(\sP) = \{\bQ_0,\bQ_1\}$, where $\bQ_1  =  \D \Symdif_{\nh\in\Nh} \trans{\nh}(\bP_1)$, and $\bQ_0  =  \Tor\setminus \bQ_1$.  Now apply {\bf(a)}.
\ethmprf

\Examples{{\rm Let $\dL=\Zahl$, $\Tor=\Torus{1}$, and $\tora\in\OO{0,1}$,
as in Example \ref{X:circle.rotation}.  Assume $\tora<\frac{1}{2}$.}}
{\item
\label{X.boolean.LCA.vs.partn}
  If $\bP_1 = \OO{0,\tora}$ and
$\bP_0 = \OO{\tora,1}$, then $\partial\sP = \{0,\tora\}$.  If
$\Phi$ is as in Example \ref{X:ledrappier}, and $\sQ=\Phitor(\sP)$,
then $\bQ_1 = \OO{0,2\tora}$
and $\bQ_0 = \OO{2\tora,1}$.  Thus, \
$\partial \sQ \quad =  \quad \{0,2\tora\}
\quad = \quad \{0,\tora\} \ \symdif \ \{\tora,2\tora\} \quad = \quad  \partial\sP \ \symdif \ \trans{1}(\partial\sP)$.  (See Figure \ref{fig:X:ledrappier},
bottom left).

\item  Let  $\torb<\tora$, and let
 $\bP_1 = \OO{0,\torb}$ and $\bP_0 = \OO{\torb,1}$.
  Then  $\bQ_1 = \OO{0,\torb}\disj\OO{\tora,\torb+\tora}$
and $\bQ_0 = \OO{\torb,\tora}\disj \OO{\torb+\tora,1}$.  Thus, $\partial \sQ 
 =  \{0,\torb,\tora,\torb+\tora\}
 = \{0,\torb\} \ \symdif \ \{\tora,\torb+\tora\}
 =  \partial\sP \ \symdif \
\trans{1}(\partial\sP)$.
}

 If $\tort\in\Tor$, recall that the {\dfn $\varsigma$-orbit} of $\tort$
is the set $\bO_\tort \ = \ 
\lb\{\trans{\ell}(\tort)\rb\}_{\ell\in\Lat}$.  
Let $\gO$ be the set of all $\varsigma$-orbits in $\Tor$;
\ then $\D\Tor \ = \ \Disj_{\bO\in\gO} \bO$, \ so
we can write $\partial \sP$ as a disjoint union:
\beqn
\label{orbit.decomposition}
  \partial\sP \quad = \quad
 \Disj_{\bO\in\gO} \partial_\bo\sP,
\eeqn
where $\partial_\bo\sP \ := \ \bO\,\intsct\,\partial\sP$ for all $\bO\in\gO$.
The decomposition in eqn.(\ref{orbit.decomposition})
commutes with the action of $\Phitor$ on
$\partial \sP$.  That is, $\partial\lb(\maketall \Phitor(\sP) \rb) 
\ = \ \D \Disj_{\bO\in\gO} \partial_\bo\lb( \Phitor(\sP) \rb)$, and
Proposition \ref{boolean.LCA.on.partn}(b) says:
\beqn
\label{phitor.vs.orbit.decomposition}
\mbox{For all $\bO\in\gO$,}\qquad
 \partial_\bo\lb( \Phitor(\sP) \rb) \quad \supseteq \quad 
\Symdif_{\nh\in\Nh} \trans{\nh}(\partial_\bo \sP).
 \eeqn
 For each $\bO\in\gO$, fix a representative $t_\bo\in\bO$, and 
define  $\beta_\bo:\OP\into\sA^\Lat$ as follows:
for any $\sP\in\OP$, let $\beta_\bo(\sP):=\bb$, where
$ \statement{$b_\ell =1$}
\iff \statement{$\trans{\ell}(t_\bo)\in\partial_\bo\sP$}$.
Then eqn.(\ref{phitor.vs.orbit.decomposition}) implies:
\beqn
\label{phi.commutes.with.beta}
\beta_\bo\lb(\Phitor(\sP)\maketall \rb)\quad \geq \quad
\Phi\lb(\beta_\bo(\sP)\maketall \rb)
\qquad \mbox{(componentwise).}
\eeqn

\example{\label{X.boolean.LCA.vs.partn.2}
Let $\Phi$, $\sP$ and $\sQ=\Phitor(\sP)$ be as in Example \ref{X.boolean.LCA.vs.partn}.
 Then all  elements of
$\partial\sP = \{0,\tora\}$ and $\partial\sQ = \{0,2\tora\}$ 
belong to the orbit $\gO_0$ of zero, and
\beq
\beta_0(\sP) & = & [\ldots0\,0\,0\,\underline{1}\,1\,0\,0\,0\ldots]\\
\mbox{while} \ 
\beta_0(\sQ) & = & [\ldots0\,0\,0\,\underline{1}\,0\,1\,0\,0\ldots]
\quad=\quad \Phi\lb(\beta_0(\sP)\rb).
\eeq
(where the zeroth element of each sequence is underlined)
}
  If $\bb\in\sA^\Lat$, let $\supp{\bb} \ := \
 \set{\ell\in\Lat}{b_\ell=1}$.  Thus, 
\beqn
\label{boundary.equiv.support}
\partial_\bo\sP \quad = \quad
 \set{\trans{\ell}(t_\bo)}{\ell\in\supp{\beta_\bo(\sP)}}.
\eeqn 
 Hence, the growth  of partition boundaries under the action of $\Phitor$
is directly related to the growth in the support of boolean configurations
under the action of $\Phi$.

\Proposition{\label{linear.CA.are.diffusive}}
{
 Let $\Phi$ be any BLCA. \nopagebreak
\bthmlist
\item  For all $n\in\Natur$, let $\Nh_n\subset\Lat$ be such that
$\Phi^n$ has local map 
$\phi_n(\ba)  =  \sum_{\nh\in\Nh_n} a_\nh \pmod{2}$,
There is a subset $\dJ\subseteq\Natur$ with $\density{\dJ}=1$ such that
$\D \lim_{\dJ\ni j \goto \oo } \#(\Nh_j)  =  \oo$.

\item  For all $\ba\in\sA^\dL$ with finite
support, there is a subset $\dJ\subseteq\Natur$ with $\density{\dJ}=1$
such that
$\D  \lim_{\dJ\ni j\goto \oo}\  \#\lb(\supp{\Phi^j(\ba)\maketall}\rb)  \ = \  \oo$.\ethmlist
}
\bthmprf  See Theorem 15 of \cite{PivatoYassawi1}. 
\ethmprf

\Corollary{\label{linear.CA.are.chopping}}
{
 If $\Tor=\Torus{1}$ and $\Phi$ is any nontrivial BLCA,
 then $\Phitor$ is \#-chopping on average.
}
\bthmprf  Let $\sP\in\OP$, and suppose  $\partial\sP$ is finite.
  For each $\bO\in\gO$, let $\bb_\bo := \beta_\bo(\sP)$; then 
$\supp{\bb_\bo}$
is finite, and is nontrivial for only
finitely many $\bO\in\gO$.  Thus, for any $j\in\Natur$,
\begin{eqnarray}\nonumber
 \#\lb(\maketall \partial\lb(\Phitor^j[\sP]\rb)\rb)
&\eeequals{(*)}&
\sum_{\bO\in\gO}  \#\lb(\maketall \partial_\bo\lb(\Phitor^j[\sP]\rb)\rb)
\quad\eeequals{(\dagger)}\quad
\sum_{\bO\in\gO} 
\# \lb[\supp{\maketall \beta_{\bo}(\Phitor^j[\sP])}\rb] 
\\& \geeeq{(\ddagger)}&
\sum_{\bO\in\gO} \#\lb(\maketall \supp{\Phi^j\lb[\bb_\bo\rb]}\rb),
\label{linear.CA.are.chopping.e1}
\end{eqnarray}
where $(*)$ is by eqn.(\ref{orbit.decomposition}); \quad
$(\dagger)$ is by eqn.(\ref{boundary.equiv.support});\quad and
$(\ddagger)$ is by eqn.(\ref{phi.commutes.with.beta}).

For each $\bO\in\gO$, 
Proposition \ref{linear.CA.are.diffusive} yields
some $\dJ_\bo\subseteq\Natur$ such that
$\density{\dJ_\bo} \ = \ 1$ and
$\D\lim_{\dJ_\bo\ni j\goto\oo}\  \#\lb(\maketall\supp{ \Phi^j(\bb_\bo)}\rb)
\ = \ \oo$.  Let $\dJ := \D \Union_{\bO\in\gO} \dJ_\bo$;  then
$\density{\dJ} \ = \ 1$, and eqn.(\ref{linear.CA.are.chopping.e1}) implies
that $\D \lim_{\dJ\ni j\goto\oo}\ 
\#\lb(\maketall \partial \lb(\Phitor^j[\sP]\rb)\rb)
\ = \ \oo$; hence $\Phitor$ chops $\sP$ on average.
\ethmprf

  To generalize Corollary \ref{linear.CA.are.chopping} to $\Torus{K}$ ($K>1$),
we need some notation.
If $\bS\subset\Tor$ is some subset, then many
 points in $\bS$ may share the same $\trans{}$-orbit.
Define
\[
  \trans{\perp}(\bS) \quad:=\quad \bS \setminus \Union_{0\neq \ell\in\Lat}
\trans{\ell}(\bS)
\quad=\quad\set{\tors\in\bS}{\bO_\tors \intsct \bS \ = \ \{\tors\} \maketall }.
\]
Let
$\OTP \ := \ \set{\sP\in\OP}{\setsize{\maketall\trans{\perp}(\partial\sP)} \ > \ 0}$.
  For example, if $\Tor=\Torus{1}$ and
$\setsize{\bullet}=\#(\bullet)$, then
\[
\statement{$\sP\in\OTP$} \iff 
\statement{$\trans{\perp}(\partial\sP)\neq \emptyset$}
 \iff   \statement{There is some $\tors\in\partial\sP$ 
which is not \\ in the orbit of any other $\tort\in\partial\sP$}.
\]

\Lemma{\label{orbital.transversality.yields.chop.lemma}}
{
If $\Phi$ is the BLCA {\rm(\ref{BLCA2})}, then 
$\#(\Nh) \cdot \setsize{\maketall\trans{\perp}(\partial\sP)}
\ \leq \  \setsize{\partial\lb(\Phitor(\sP)\rb)}
\ \leq \ 
\#(\Nh) \cdot \setsize{\partial\sP}$.
}
\bthmprf
Let $\bS = \trans{\perp}(\partial\sP)$ and let 
$\bU=\partial\sP \setminus \bS$.  By definition of $\bS$,
the sets $\{\trans{\torb}(\bS)\}_{\nh\in\Nh}$ are disjoint both from
one another and from the set $\D \Union _{\nh\in\Nh} \, \trans{\torb}(\bU)$.
Thus, Proposition \ref{boolean.LCA.on.partn}(b) says:
\beqn
\label{orbital.transversality.yields.chop.lemma.e1}
  \Disj_{\nh\in\Nh} \trans{\torb}(\bS)
\quad \subseteq \quad
 \Symdif_{\nh\in\Nh} \trans{\nh}(\partial \sP)
 \quad \subseteq \quad 
\partial\lb(\maketall \Phitor(\sP)\rb)
 \quad \subseteq \quad 
\Union_{\nh\in\Nh} \trans{\nh}(\partial \sP).
\eeqn
 Thus, $\begin{array}[t]{rcl} \#(\Nh) \cdot \setsize{\bS}
& \eeequals{(M3)} &  \D
\sum_{\nh\in\Nh} \setsize{\trans{\torb}(\bS)}
\quad\eeequals{(M2)} \quad
\setsize{  \Disj_{\nh\in\Nh} \trans{\torb}(\bS)}
\quad\leeeq{(M1)} \quad
\setsize{\maketall \partial\lb(\Phitor(\sP)\rb)}
\\ & \leeeq{(M1)} &\D
\setsize{  \Union_{\nh\in\Nh} \trans{\torb}(\partial\sP)}
\quad\leeeq{(M2)} \quad
\sum_{\nh\in\Nh} \setsize{\trans{\torb}(\partial\sP)}
\quad\eeequals{(M3)} \quad 
  \#(\Nh) \cdot \setsize{\partial\sP}.
\end{array}$

The {\bf(M1)} inequalities follow from 
eqn.(\ref{orbital.transversality.yields.chop.lemma.e1})
and property {\bf(M1)} of $\setsize{\bullet}$.
\ethmprf

\Proposition{\label{orbital.transversality.yields.chop}}
{
 Let $\Tor=\Torus{K}$. If $\sP\in\OTP$,
then any nontrivial BLCA chops $\sP$ on average.
}
\bthmprf   Let $\Phi$ be a BLCA. For all
 $n\in\Natur$, let $\Nh_n\subset\Lat$ be such that
$\Phi^n$ has local map 
$\phi_n(\ba) \ = \ \sum_{\nh\in\Nh_n} a_\nh \pmod{2}$.
Thus, Lemma \ref{orbital.transversality.yields.chop.lemma} says that
$\setsize{\partial\lb(\Phitor^n(\sP)\rb)}
\ \geq \ 
\#(\Nh_n) \cdot \setsize{\trans{\perp}(\partial\sP)}$.
Recall  $\sP\in\OTP$, so $\setsize{\trans{\perp}(\partial\sP)} >0$.
 Proposition \ref{linear.CA.are.diffusive}(a) yields 
a subset $\dJ\subset\Natur$ of density one such that
$\D \lim_{\dJ\ni j \goto \oo } \#(\Nh_j) \ = \ \oo$.
Thus,
$\D\lim_{\dJ\ni j \goto \oo } \setsize{\partial\lb(\Phitor^j(\sP)\rb)}
\ = \ \oo$.
\ethmprf

Let
$\SOTP := \ \set{\sP\in\OTP}{\setsize{\maketall\trans{\perp}(\partial\sP)}
 =  \setsize{\partial\sP}}.$  For example, 
if $\Tor=\Torus{1}$ and
$\setsize{\bullet}=\#(\bullet)$, then
\[
\statement{$\sP\in\SOTP$}\iff
\statement{ $\trans{\perp}(\partial\sP)=\partial\sP$}
\iff
\statement{Every element of 
$\partial\sP$ \\ occupies  a distinct $\trans{}$-orbit}.
\]

\Lemma{\label{SOT.dense}}
{ 
 $\SOTP$ is a $d_\symdif$-dense subset of $\OP$.
}
\bthmprf
Let $\DP$ be the set of dyadic open partitions of $\Tor$
(see \S\ref{app:open.dyadic}).  
 Proposition \ref{DP.dense.in.MP} says that 
$\DP$ is a $d_\symdif$-dense subset of $\OP$.  Thus,
it suffices to show that $\DP\subseteq\SOTP$.

 To see this, 
suppose $\sD\in\DP$ is an $n$-dyadic partition for some $n\in\Natur$.
Identify $\Tor\cong\CC{0,1}^K$ as usual.
Let $\bE=\{\frac{j}{2^n}\}_{j=0}^{j=2^n}$ be the set of $n$-dyadic numbers,
and for each $k\in\CC{1..K}$, let $\bC_k \ = \ \CC{0,1}^{k-1}\x \bE
\x \CC{0,1}^{K-k}$.  If $\bC =  \Union_{k=1}^K \bC_k$,
then $\partial\sD\subset \bC$. 

For any nonzero $\ell\in\Lat$, \quad
$\setsize{\trans{\ell}(\bC)\intsct\bC}=0$, because $\trans{\ell}(\bC)$ and
$\bC$ intersect transversely.  But $\partial\sD\subset\bC$, so
$\setsize{\trans{\ell}(\partial\sD)\intsct(\partial\sD)}=0$ also.
Thus,  $\setsize{\maketall\trans{\perp}(\partial\sD)} \ = \ 
\setsize{\partial\sD}$.  Hence $\sD\in\SOTP$.  

  Since this holds for any $\sD\in\DP$, we conclude that 
$\DP\subseteq\SOTP$, as desired.
\ethmprf
We'll now precisely characterize the growth rate of 
$ \setsize{\maketall\partial \lb(\Phitor^n[\sP]\rb)}$ 
for a particular BLCA.

\Proposition{\label{X.boolean.LCA.vs.partn.growth.thm}}
{

Let $\Phi:\sA^\Zahl\into\sA^\Zahl$ be the BLCA of Example
{\rm\ref{X:ledrappier}}.  Let $\trans{}$ be a $\Zahl$-action on
$\Tor=\Torus{K}$, and let $\sP\in\OTP$.  Then, as $n\goto\oo$...
\bthmlist

  \item ...the maximum of $\setsize{\partial\lb(\Phitor^n[\sP]\rb)}$ grows linearly.  That is:
\[ 
0 \quad < \quad \setsize{\trans{\perp}(\partial\sP)\maketall }
\quad  \leq\quad 
\limsup_{n\goto\oo}  \frac{1}{n}\ \setsize{\partial\lb(\Phitor^n[\sP]\rb)}
\quad \leq\quad  
\setsize{\partial\sP}.
\]
  \item ...the minimum of $\setsize{\partial\lb(\Phitor^n[\sP]\rb)}$
 remains constant:
\[
 \liminf_{n\goto\oo} \ \setsize{\partial\lb(\Phitor^n[\sP]\rb)} 
\quad  \leq \quad  2\setsize{\partial\sP}.
\]
  \item ...the average  of $\setsize{\maketall\partial \lb(\Phitor^n[\sP]\rb)}$
grows like $n^\alp$, where $\alp := \log_2\lb(\frac{3}{2}\rb)$. That is:
\[
\mbox{If} \ 
A(N) \quad := \quad \D\frac{1}{N} \sum_{n=0}^{N-1} \setsize{\maketall\partial \lb(\Phitor^n[\sP]\rb)} ,
\qquad \mbox{then}
\qquad
 \lim_{N\goto\oo} \frac{\log\lb(\maketall A(N)\rb) }{\log(N)} \quad=\quad \alp.
\]

  \item Finally, both {\bf(a)} and {\bf(b)} become equalities if $\sP\in\SOTP$.
\ethmlist
}
\bthmprf
  For all $n\in\Natur$, let $\Nh_n\subset\Lat$ be
such that $\Phi^n$ has local map 
$\phi_n(\ba) \ = \ \sum_{\nh\in\Nh_n} a_\nh \pmod{2}$.
Let $\nu(n)$ be the number of $1$'s in the binary expansion of
$n$.   

\Claim{\label{X.boolean.LCA.vs.partn.growth.thm.C2}
For any $n\in\Natur$, \quad $
2^{\nu(n)}\setsize{\trans{\perp}(\partial\sP)\maketall } \ \leq \ \setsize{\partial\lb(\Phitor^n[\sP]\rb)} 
\ \leq \ 2^{\nu(n)}\setsize{\partial\sP}$, with equality
when $\sP\in\SOTP$.
}
\bclaimprf
For any $n\in\Natur$, Lucas' Theorem (see \S\ref{app:LCA})
 implies that
$ \#\,(\Nh_n) \ = \  2^{\nu(n)}$.
(For example, $n=75$ has binary expansion $1001011$, so 
 $\#\lb(\Nh_{75}\rb) = 2^{\nu(75)} = 2^4= 16$.) 
 Now apply Lemma \ref{orbital.transversality.yields.chop.lemma}.
\eclaimprf

 {\bf(a)}\quad To see that
$\D\limsup_{n\goto\oo} \
\frac{1}{n}\ \setsize{\maketall\partial\lb(\Phitor^n[\sP]\rb)}
\ \leq  \ \setsize{\partial\sP}$,
observe that $\Nh_n\subset \CC{0..n}$.
Thus,   $\#\lb(\Nh_n\rb) \ \leq \ \#\CC{0..n} \ = \ n+1$.  Thus,
\beq
\limsup_{n\goto\oo} \ \frac{1}{n}\setsize{\partial \lb(\Phitor^n[\sP]\rb)}  
& \leeeq{(*)} &
 \limsup_{n\goto\oo} \ \frac{1}{n} \ \#\lb(\Nh_n\rb)\cdot\setsize{\partial\sP} 
\quad \leq \quad
\lim_{n\goto\oo} \ \frac{n+1}{n} \cdot\setsize{\partial\sP} 
\\&=&
\setsize{\partial\sP},
\eeq
where $(*)$ is by
Lemma \ref{orbital.transversality.yields.chop.lemma}.

To see that $\D \limsup_{n\goto\oo} \ \frac{1}{n}\setsize{\partial \lb(\Phitor^n[\sP]\rb)}\ \geq \ \setsize{\trans{\perp}(\partial\sP)\maketall } $,
 let $n = 2^{m}-1$ for some $m\in\Natur$.  Then
$\nu(n) \ = \ m$, so  Claim \ref{X.boolean.LCA.vs.partn.growth.thm.C2}
says \ $\D \setsize{\maketall\partial \lb(\Phitor^n[\sP]\rb)}
  \ \geq \  2^{\nu(n)} \cdot \setsize{\trans{\perp}(\partial\sP)\maketall } \ = \  2^m \cdot \setsize{\trans{\perp}(\partial\sP)\maketall }$.
Thus,
\beq
 \limsup_{n\goto\oo} \ \frac{\setsize{\partial \lb(\Phitor^n[\sP]\rb)}}{n}
& \geq &
 \limsup_{m\goto\oo} \ \frac{\setsize{\partial \Phitor^{(2^m-1)}(\sP)}}{2^m-1}
\quad \geq \quad \lim_{m\goto\oo} \ \frac{2^m \cdot \setsize{\trans{\perp}(\partial\sP)\maketall }}{2^m-1}
\\&=&
\setsize{\trans{\perp}(\partial\sP)\maketall }.
\eeq
{\bf(b)}  \quad
If $n = 2^{m}$ for some $m\in\Natur$, then
$\nu(n) \ = \ 1$, so  Claim \ref{X.boolean.LCA.vs.partn.growth.thm.C2}
 says \ $\D \setsize{\maketall\partial \lb(\Phitor^n[\sP]\rb)}
  \ \leq \  2^{\nu(n)} \cdot \setsize{\partial\sP} 
\ = \  2 \cdot \setsize{\partial\sP}$.
Thus, \ $\D
 \liminf_{n\goto\oo} \setsize{\maketall\partial \lb(\Phitor^n[\sP]\rb)}
\ \leq \ \liminf_{m\goto\oo}  \setsize{\partial \Phitor^{(2^m)}(\sP)}
\ = \ 2 \cdot \setsize{\partial\sP}$.

{\bf(d)} \quad follows from the second part of Claim
\ref{X.boolean.LCA.vs.partn.growth.thm.C2}.

{\bf(c)}  \quad For any $n\in\Natur$, let  $f(n) = 2^{\nu(n)}$, and
for any $N\in\Natur$, let
$\D \tlA(N) \ = \ \frac{1}{N} \sum_{n=0}^{N-1} f(n)$.  
Thus, Claim \ref{X.boolean.LCA.vs.partn.growth.thm.C2} implies that
\
$\tlA(N)\cdot\setsize{\trans{\perp}(\partial\sP)\maketall }
 \leq \  A(N)\  \leq  \  \tlA(N)\cdot\setsize{\partial\sP}$.
Hence,
\[
1 \ = \ 
\lim_{N\goto\oo} 1 + 
 \frac{\log\setsize{\trans{\perp}(\partial\sP)\maketall}}{\log(\tlA(N))}
 \ \ \leq \ \
\lim_{N\goto\oo} \frac{\log(A(N))}{\log(\tlA(N))}
\ \ \leq \ \
 \lim_{N\goto\oo} 1 +  \frac{\log\setsize{\partial\sP}}{\log(\tlA(N))}
\ \ = \ \ 1.
\]
Thus, it suffices to examine the asymptotics of $\tlA(N)$.

Suppose $N= 2^{M_0} +2^{M_1} + \cdots + 2^{M_J}$
for some  $M_0 > M_1  > \cdots > M_J$.

\Claim{\label{X.boolean.LCA.vs.partn.growth.thm.eqn2}  
$\D \tlA(N) \ = \ \frac{1}{N}\sum_{j=0}^J 2^j\cdot 3^{M_j}$}
\bclaimprf
  Let $n$ be a random element of $\CO{0..N}$ (with uniform distribution).
Then $f(n)$ is also a random variable, and  $\D \tlA(N) \ = \ \Expct{f(n)}$
is the expected value of $f(n)$.
Let $\bI_0 := \CO{0..2^{M_0}}$, \ $\bI_1 := 
\CO{2^{M_0}\ldots 2^{M_0}\!+\!2^{M_{1}}}$, and in general,
\[
\bI_j
 \quad := \quad 
\CO{\sum_{i=0}^{j-1} 2^{M_{i}} \ldots \sum_{i=0}^{j} 2^{M_{i}}},
\qquad \mbox{for all $j\in\CC{0..J}$.}
\]
Then $\CO{0..N} \ = \ \bI_0 \disj \bI_1 \disj \cdots \disj \bI_J$, and
for all $j\in\CC{0..J}$,\quad
$\D\dP\lb( n\in\bI_j\rb) \quad=\quad {2^{M_j}}/{N}$.

\subclaim{\label{X.boolean.LCA.vs.partn.growth.thm.C3}
Suppose $n\in\bI_0$.   Then $\Expct{f(n) \given n\in\bI_0} \ = \ 
 \lb(\frac{3}{2}\rb)^{M_0}$.}
\bsubclaimprf
Write $n$ in binary notation.  Then the $M_0$ binary digits of $n$ are
independent, equidistributed boolean random variables, so $\nu(n)$
is a random variable with a binomial distribution:
 \beqn
\label{X.boolean.LCA.vs.partn.growth.thm.eqn1}
  \dP\lb(\nu(n) = m\maketall\rb) \quad=\quad \frac{1}{2^{M_0}}\BINOM{M_0}{m}
\qquad \mbox{for any $m\in\CC{0...M_0}$.}
\eeqn
Thus, 
$\begin{array}[t]{rcl}
\Expct{f(n)} &=&
\D\sum_{m=0}^{M_0} 2^m \cdot \dP\lb(\maketall f(n) = 2^m\rb) 
\quad=\quad
\sum_{m=0}^{M_0} 2^m \cdot \dP\lb(\maketall\nu(n) = m\rb) 
\\ &\eeequals{(*)} &
\D \frac{1}{2^{M_0}}\sum_{m=0}^{M_0}  2^m \cdot  \BINOM{M_0}{m}
\ \ \eeequals{(B)} \ \
\frac{1}{2^{M_0}} (1+2)^{M_0}
\  = \  
\lb(\frac{3}{2}\rb)^{M_0},
\end{array}$

where $(*)$ is by
eqn.(\ref{X.boolean.LCA.vs.partn.growth.thm.eqn1}), and {\bf(B)}
is the Binomial Theorem.
\esubclaimprf

\subclaim{\label{X.boolean.LCA.vs.partn.growth.thm.C4}
For any $j\in\CC{0..J}$, \quad
$\Expct{f(n) \given n\in\bI_j} \ = \ 2^j\cdot \lb(\frac{3}{2}\rb)^{M_j}$.
}
\bsubclaimprf
 If $j=0$, this is Claim \ref{X.boolean.LCA.vs.partn.growth.thm.eqn2}.\ref{X.boolean.LCA.vs.partn.growth.thm.C3}.
Let  $j\geq 1$.  If $n\in\bI_j$, then
$n = 2^{M_0} + \cdots + 2^{M_{j-1}} + n_1$, for some
$n_1\in\CO{0...2^{M_j}}$.
Thus, $\nu(n) = j + \nu(n_1)$, so that $f(n) \ = \  2^j\cdot f(n_1)$.
Thus 
\[
  \Expct{f(n) \given n\in\bI_j} \quad = \quad 
 2^j\cdot \Expct{f(n_1) \given n_1\in\CO{0...2^{M_j}}} 
\quad \eeequals{(*)} \quad 
2^j\cdot (\mbox{$\frac{3}{2}$})^{M_j}. 
\]
where $(*)$
is like Claim \ref{X.boolean.LCA.vs.partn.growth.thm.eqn2}.\ref{X.boolean.LCA.vs.partn.growth.thm.C3}.
\esubclaimprf
\beq
\mbox{It follows that}\qquad
\tlA(N) &=& \Expct{f(n)}
\quad=\quad
\sum_{j=0}^J \Expct{f(n) \given n\in\bI_j} \cdot\dP\lb( n\in\bI_j\rb)
\\ 
&\eeequals{(*)}&
\sum_{j=0}^J 2^j\cdot \lb(\frac{3}{2}\rb)^{M_j} \cdot \frac{2^{M_j}}{N}
\quad=\quad
\frac{1}{N}\sum_{j=0}^J 2^j\cdot 3^{M_j}, 
\eeq
where $(*)$ is by Claim 
\ref{X.boolean.LCA.vs.partn.growth.thm.eqn2}.\ref{X.boolean.LCA.vs.partn.growth.thm.C4}.
\eclaimprf

\Claim{\label{X.boolean.LCA.vs.partn.growth.thm.C5}
$\tlA(N) \ > \ \frac{1}{3} N^\alp$.}
\bclaimprf
$\D \tlA(N)
 \quad \eeequals{(C\ref{X.boolean.LCA.vs.partn.growth.thm.eqn2})} \quad 
  \frac{1}{N}\sum_{j=0}^J 2^j\cdot 3^{M_j}
\ \geq \ 
\frac{3^{M_0}}{N} 
\quad \grt{(\star)} \quad
\frac{1}{3} \frac{3^{M_0+1}}{2^{M_0+1}}
\ = \ \frac{1}{3} 2^{\alp\cdot(M_0+1)} 
\quad \grt{(\star)} \quad \frac{1}{3} N^\alp$.
 
{\bf(C\ref{X.boolean.LCA.vs.partn.growth.thm.eqn2})} is by
Claim \ref{X.boolean.LCA.vs.partn.growth.thm.eqn2},
 and the $(\star)$ inequalities are because
$N<2^{M_0+1}$.
\eclaimprf

\Claim{\label{X.boolean.LCA.vs.partn.growth.thm.C6}
$\tlA(N) \ \leq \ 3\cdot N^\alp$.}
\bclaimprf
Since $M_0 > M_1 > \cdots  > M_J$, we know that $M_j \leq M_0 - j$ for
all $j\in\CC{1..J}$.  Thus,
\beq
\tlA(N)
 & \eeequals{(C\ref{X.boolean.LCA.vs.partn.growth.thm.eqn2})} &
  \frac{1}{N}\sum_{j=0}^J 2^j\cdot 3^{M_j}
\quad \leq \quad
  \frac{1}{N}\sum_{j=0}^J 2^j\cdot 3^{M_0-j}
\quad=\quad
\frac{3^{M_0}}{N}\sum_{j=0}^J \lb(\frac{2}{3}\rb)^j
\\&=&
\frac{3^{M_0}}{N}\lb(\frac{1-\lb(\frac{2}{3}\rb)^{J+1}}{1-\frac{2}{3}}\rb)
\quad=\quad
\frac{3^{M_0+1}}{N}\lb(1-\lb(\frac{2}{3}\rb)^{J+1}\rb)
\quad\leq\quad
\frac{3^{M_0+1}}{N}
\\&=&
 3\cdot \frac{3^{M_0}}{N}
\quad\leeeq{(\star)}\quad 3\cdot \frac{3^{M_0}}{2^{M_0}}
\quad = \quad 3\cdot 2^{\alp M_0}
\quad\leeeq{(\star)}\quad 3\cdot N^{\alp}
\eeq
{\bf(C\ref{X.boolean.LCA.vs.partn.growth.thm.eqn2})} is by 
Claim \ref{X.boolean.LCA.vs.partn.growth.thm.eqn2}, and the $(\star)$ inequalities are because
$2^{M_0} \ \leq \  N$.
\eclaimprf
$
\begin{array}{rcccl}
\mbox{Combining Claims \ref{X.boolean.LCA.vs.partn.growth.thm.C5} and
\ref{X.boolean.LCA.vs.partn.growth.thm.C6} yields}
\ \ \D \frac{1}{3} N^\alp & < & \tlA(N) & \leq & 3\cdot N^\alp; \\ \\
\mbox{Hence,} \qquad 
\alp\log(N)-\log(3) & < & \log\lb(\maketall \tlA(N)\rb) & \leq & 
\alp\log(N)+\log(3); \\ \\
\mbox{Hence,} \qquad
\D\alp-\frac{\log(3)}{\log(N)} & < & \D \frac{\log\lb(\maketall \tlA(N)\rb)}{\log(N)} & \leq & 
\D\alp+\frac{\log(3)}{\log(N)}.
\end{array}
$

Taking the limit as $N\goto\oo$, we conclude that $\D\lim_{N\goto\oo} 
\frac{\log\lb(\tlA(N)\rb)}{\log(N)}  =  \alp$, as desired.\nolinebreak
\ethmprf
{\em Remarks.} {\bf(i)}\quad Proposition \ref{X.boolean.LCA.vs.partn.growth.thm}{\bf(b)}
shows that, in general, we can only expect chopping to occur along
a subset of $\dJ\subset\Natur$ of density one.

 {\bf(ii)}\quad We proved Proposition \ref{X.boolean.LCA.vs.partn.growth.thm}
for the very simple  BLCA of Example \ref{X:ledrappier}.
Similar  results are probably 
true for arbitrary BLCA, but the appropriate version of Claim 
\ref{X.boolean.LCA.vs.partn.growth.thm.C2} will be much more complex in
general, leading to more complex formulae in parts {\bf(a)}, {\bf(b)} and
{\bf(c)}.

\paragraph*{\sc The Lipschitz pseudomeasure:}

Recall that $\D \setsize[L]{\bS} \ = \ \lim_{\eps\goto 0} \,
\frac{1}{2\eps}\, \lam\lb(\maketall \Ball(\bS,\eps)\rb)$. \
This is an increasing
limit, because for any $\eps>0$ and $n\in\Natur$,
\quad $\lam\lb(\maketall \Ball(\bS,\eps)\rb) \leq
n\cdot \lam\lb(\maketall \Ball(\bS,\frac{\eps}{n})\rb)$.  Thus:
\beqn
\label{set.size.defn.1}
\mbox{For all $\eps>0$,}\quad
\lam\lb(\maketall \Ball(\bS,\eps)\rb) \quad \leq \quad
 2\eps\cdot\setsize[L]{\bS}.
\eeqn
If $\sP\in\OP$, then $\setsize[L]{\partial\sP}$ affects the continuity properties of $\Splat$ as follows:

\Proposition{\label{splat.is.lipschitz}}
{
  Let $\sP\in\OP$, and endow $\sA^\Lat$ with the Besicovitch metric.  Then:
\bthmlist
  \item  $\Splat:\tlTor \into \sA^\Lat$ is Lipschitz, with
Lipschitz constant $\setsize[L]{\partial \sP}$.  That is, for any
$\tors,\tort\in\tlTor$, \quad $d_B\lb(\maketall \Splat(\tors),\Splat(\tort)\rb)
\ \leq \ \setsize[L]{\partial \sP} \cdot d(\tors,\tort)$.

  \item Let $\Tor = \Torus{1}$  {\rm(so $\setsize[L]{\bullet} = \#(\bullet)$)},
 and let $\Sep{\sP} = \min\set{\maketall d(\torb_1,\torb_2)}{\torb_1,\torb_2\in\partial\sP \ \mbox{distinct}}$. 
 For any
 $\tors,\tort\in\tlTor$, \quad
$\statement{ $d(\tors,\tort)\leq \Sep{\sP}$}\IMPLIES
\statement{ $d_B\lb(\maketall \Splat(\tors),\Splat(\tort)\rb) \ = \
 \#(\partial \sP) \cdot d(\tors,\tort)$}$.
\ethmlist
}
\bthmprf {\bf(a)}
  Identify $\Tor\cong \CO{0,1}^K$.  For simplicity assume $\tors = 0$,
and let $\tort=(t_1,\ldots,t_K)$, where
 $t_k\in\CO{0,1}$ for all $k\in\CC{1..K}$. 
 Assuming that $
\norm{t_1,\ldots,t_k}{} \ \leq \ \norm{(1-t_1),\ldots,(1-t_k)}{}$,
the shortest path from $0$ to $\tort$ is along the line segment
$  \bI \ = \ 
 \set{(rt_1,\ldots,rt_k)}{r\in\CC{0,1}}$.
  Thus, if $d(0,\tort)=\eps$, then $\length{\bI}=\eps$. 
Let  $\torm= \lb(\frac{1}{2}t_1,\ldots,\frac{1}{2}t_k\rb)$ be the
midpoint of $\bI$, and let $\Ball\lb(\torm,\frac{\eps}{2}\rb)$ be the ball of
radius $\eps/2$ about $\torm$. Then, for any $\ell\in\Lat$, 
\begin{eqnarray}
\nonumber
\hspace{-4em}\lefteqn{\statement{$\Splat(0)_\ell \ \neq \ \Splat(\tort)_\ell$}
 \iff  
\statement{$\sP\lb(\trans{\ell}(0)\rb) \ \neq \ 
\sP\lb(\trans{\ell}(\tort)\rb)$}
 \iiimplies{(\dagger)}
\statement{$\partial\sP\, \intsct\,\trans{\ell}(\bI) \ \neq \ \emptyset$}}
\qquad\quad
\\  \label{splat.is.lipschitz.eqn1}
 & \iiimplies{(\star)}&
\statement{$\partial\sP\, \intsct\, \Ball\lb(\trans{\ell}(\torm), \ \frac{\eps}{2}\rb) \ \neq \ \emptyset$}
\  \iff \ 
\statement{$\trans{\ell}(\torm) \ \in \ \Ball\lb(\partial\sP,\frac{\eps}{2}\rb)$}.
\end{eqnarray}
$(*)$ is because $\bI\subset \Ball\lb(\torm,\frac{\eps}{2}\rb)$.
  If $K=1$, then  $\bI=\Ball\lb(\torm,\frac{\eps}{2}\rb)$, and the ``$\iiimplies{(\star)}$'' is actually a ``$\iff$''; \  if also $\eps<\Sep{\sP}$, then 
the ``$\iiimplies{(\dagger)}$'' is also a ``$\iff$''.
Thus,
\beq
 d_B\lb(\maketall \Splat(0),\Splat(\tort)\rb)
&=&
\density{\ell\in\Lat \ ; \ 
  \Splat(0)_\ell \ \neq  \ \Splat(\tort)_\ell}
\\ &\leeeq{(\star)}&
\density{\ell\in\Lat \ ; \ 
 \trans{\ell}(\torm) \ \in \  \Ball\lb(\partial\sP,\mbox{$\frac{\eps}{2}$}\rb)}
\quad\eeequals{(\dagger)}\quad 
\lam\lb[\maketall\Ball\lb(\partial\sP,\mbox{$\frac{\eps}{2}$}\rb)\rb].
\eeq
$(\star)$ is by (\ref{splat.is.lipschitz.eqn1}), and is an equality
if $K=1$ and $\eps<\Sep{\sP}$.
  $(\dagger)$ is by Proposition \ref{tau.monothetic.uniqely.ergodic}(b).
  Thus, 
\beqn
\label{splat.is.lipschitz.e0}
 d_B\lb(\maketall \Splat(0),\Splat(\tort)\rb)
\quad\leeeq{(\star)}\quad
\lam\lb[\maketall\Ball\lb(\partial\sP,\mbox{$\frac{\eps}{2}$}\rb)\rb]
\quad\leeeq{(e\ref{set.size.defn.1})}\quad
\eps \cdot \setsize[L]{\partial\sP},
\eeqn
where {\bf(e\ref{set.size.defn.1})} is by eqn (\ref{set.size.defn.1}).
If $K=1$ and $\eps<\Sep{\sP}$, then $(\star)$ is again an equality.  

{\bf(b)}\quad 
  Suppose $K=1$ and $\partial\sP$ is finite.  If $\eps<\Sep{\sP}$, then
$\Ball\lb(\partial\sP,\frac{\eps}{2}\rb) \ = \ \D \Disj_{\torb\in\partial\sP}
\Ball\lb(\torb,\mbox{$\frac{\eps}{2}$}\rb)$.  Thus,
$\D
 d_B\lb(\maketall \Splat(0),\Splat(\tort)\rb)
\ \eeequals{by\,(\ref{splat.is.lipschitz.e0})}\ 
\lam\lb(\maketall\Ball\lb(\partial\sP,\mbox{$\frac{\eps}{2}$}\rb)\rb)
\ = \ 
\sum_{\torb\in\partial\sP} \lam\lb(\Ball\lb(\torb,\frac{\eps}{2}\rb)\rb)
\ = \ 
\sum_{\torb\in\partial\sP} \eps
\ = \
 \eps\cdot \#(\partial \sP)$. \ethmprf

{\sc Note:} \  The equality in Proposition  \ref{splat.is.lipschitz}(b) fails if
$d(\tors,\tort)>\Sep{\sP}$. For example, let $\sA=\{0,1\}$, and
let $\bP_0 = \OO{0,\frac{1}{4}}\disj
\OO{\frac{501}{1000},\frac{3}{4}}$, while $\bP_1 =
\OO{\frac{1}{4},\frac{501}{1000}}\disj\OO{\frac{3}{4},1}$.  Thus,
$\#(\partial\sP) \ = \ 4$. However, $\Sep{\sP} =
\frac{249}{1000}<\frac{1}{4}$, and if $d(\tors,\tort)=\frac{1}{2} > \Sep{\sP}$,
then $d_B(\Splat(\tors),\Splat(\tort)) \ = \  \frac{1}{1000}
\ \neq  \ 4\cdot d(\tors,\tort)$.

\paragraph*{\sc Sensitivity and $\setsize[L]{\bullet}$-Chopping}

Let $(\bX,d,\varphi)$ be a topological dynamical system, and
let $\xi>0$.  If $x,y\in\bX$, then $(x,y)$ is a {\dfn $\xi$-expansive pair} if $d\lb(\varphi^n(y),\varphi^n(y_1)\rb)>\xi$ for some
$n\in\Natur$.  We say $(\bX,d,\varphi)$ is {\dfn $\xi$-expansive}  (see \S\ref{S:expansive}) if
$(x,y)$ is $\xi$-expansive for every $x,y\in\bX$ with $x\neq y$. 
If $\bY\subset\bX$ is a subset (not necessarily
$\varphi$-invariant), 
then  $\varphi\restr{\bY}$ is {\dfn $\xi$-sensitive} if
for all $y\in\bY$, and all $\del>0$, there some $y_1\in\bY$ with
$d(y,y_1)<\del$ such that $(y,y_1)$ is a $\xi$-expansive pair. Thus,
if $(\bX,d,\varphi)$ is  $\xi$-expansive, then
$\phi\restr{\bY}$ is $\xi$-sensitive for any $\bY\subseteq\bX$.

\Proposition{\label{expansive.implies.chopping}}
{
\bthmlist
  \item 
If $\sP\in\OP$, and $\tlgP = \Splat(\tlTor)$, then
\[
\statement{$(\tlgP,d_B,\Phi)$ is sensitive}
\IMPLIES\statement{$\Phitor$ intermittently $\setsize[L]{\bullet}$-chops
 $\sP$}.
\]
  \item Thus,
$\statement{$(\QS,d_B,\Phi)$ is expansive}
\IMPLIES\statement{$\Phitor$ is intermittently $\setsize[L]{\bullet}$-chopping}$.
\ethmlist
}
\bthmprf {\bf(b)} follows from {\bf(a)}, so we'll prove {\bf(a)}.

{\bf Case 1:} ({\em  $\sP$ is simple})
\quad
  Let $\bp_0 \ = \ \Splat(0)$, and fix $\eps>0$.  
Proposition \ref{open.simple.partn.inj}(b) says
$\Splat$ is a homeomorphism from $(\tlTor,d)$ to $(\tlgP,d_B)$.
  Thus, there is some $\del>0$ such that, for any
$\tort\in\tlTor$,\quad $\statement{$d_B\lb(\Splat(0),\Splat(\tort)\rb)
\ < \ \del$} \IMPLIES \statement{$d(0,\tort)<\eps$}$.

Suppose $(\tlgP,d_B,\Phi)$ is $\xi$-sensitive for some
$\xi>0$.  Then there is some
$\bp_1\in\tlgP$ with $d_B(\bp_0,\bp_1) < \del$, but 
$d_B\lb( \maketall \Phi^n(\bp_0),\Phi^n(\bp_1)\rb) \ > \ \xi$ for some
$n\in\Natur$. 

  Suppose $\bp_1 = \Splat(\tort_1)$, where $\tort_1\in\tlTor$ and
$d(0,\tort_1)<\eps$.  Let $\sP^{(n)} = \Phitor^n(\sP)$.  Thus, we have
\[
 \frac{\xi}{\eps}
\quad \leq \quad
 \frac{\xi}{d(0,\tort_1)}
\quad \leq \quad
  \frac{d_B\lb(\Splat^{(n)}(0),\Splat^{(n)}(\tort_1)\rb)}
  {d(0,\tort_1)}
\quad \lt{(*)} \quad
 \setsize[L]{\partial \sP^{(n)}},
\]
where $(*)$ is by
Proposition \ref{splat.is.lipschitz}(a).
But $\eps$ can be made arbitrarily small. 
Hence, $\setsize[L]{\partial \sP^{(n)}}$  can become 
arbitrarily large as $n\goto\oo$.

{\bf Case 2:} ({\em $\sP$ is not simple})\quad Use
Lemma \ref{quotient.partition} to replace
$\sP$ with a `quotient' partition $\quosP$ on a quotient torus
$\quoTor$, with a quotient $\dL$-action $\quotrans{}$, such that
$\quosP$ {\em is} simple.

\Claim{\label{expansive.implies.chopping.C1}\hspace{-1em}
 Let $\bS$ be the symmetry group of $\sP$, and
let $q\colon\Tor\rightarrow\quoTor$ be the
quotient map, as in
 {\rm Lemma \ref{quotient.partition}}.
Let
 \[
\OP_\bS  :=  \set{\sQ\in\OP}{\mbox{all elements of $\bS$ are symmetries of
$\sQ$}}.  
\]
  {\bf(a)} \  For any $\sQ\in\OP_\bS$, there is a unique  
 partition $\quosQ\in\quoOP$ such that $\quosQ\circ q = \sQ$.  

  {\bf(b)} \   $\Phitor(\OP_\bS) \subset \OP_\bS$.  
 In particular, if $\sQ=\Phitor(\sP)$, then $\sQ\in\OP_\bS$. 

  {\bf(c)} \  Let $\quoPhitor:\quoOP\into\quoOP$ be
the induced map on partitions of $\quoTor$.  Then $\quoPhitor(\quosP)
= \overline{\Phitor(\sP)}$. 

{\bf(d)} \  There is a constant $C>0$ such that, for any $\sQ\in\OP_\bS$,
\quad $\setsize[L]{\partial \sQ} \ = \ C\cdot \setsize[L]{\partial\quosQ}$. 

{\bf(e)} \   Hence, 
$\statement{$\Phitor$ intermittently chops $\sP$}
\iff
\statement{$\quoPhitor$ intermittently chops $\quosP$}$.
}
\bclaimprf 
{\bf(a)} is like Lemma \ref{quotient.partition}{\bf(b)}.
\quad
{\bf(b)} follows by applying $\bS$-symmetries to eqn.(\ref{partition.image}),
and {\bf(c)} follows by combining  eqn. (\ref{partition.image}) with 
Lemma \ref{quotient.partition}{\bf(b,c)}.

{\bf(d)}\quad  $\bS$ is a closed subgroup of $\Tor$,
so $\bS$ is a smooth embedding of a Lie group $\Torus{J}\x\bS_0$, where
$J\leq K$ and $\bS_0$ is a finite abelian group
(if $\bS$ is finite, then $J=0$ and $\bS=\bS_0$).
Let $V_{K-J}(\eps)$ be the volume of a ball of radius $\eps$ in $\Real^{K-J}$, and let
$\D  C \  := \ \D\lim_{\eps\goto0}\
\frac{\lam\lb[\Ball(\bS,\eps)\rb]}{V_{K-J}(\eps)}$; \
then $C$ measures the `size' of $\bS$ (if $\bS$ is finite, then $J=0$, and
$C=\#(\bS)$.)  

  Now, $\Tor$ is a smooth locally
trivial fibre bundle over $\quoTor$, with fibre $\bS$
\cite[Thm 3.58, p.120]{Warner}.
Thus, $\partial\sQ$ is a 
 sub-bundle over $\partial\quosQ$, while
$\Ball(\partial\sQ,\eps)$ is a smooth sub-bundle over
$\Ball(\partial\quosQ,\eps)$ (both also with  fibre $\bS$).
Suppose $\quolam$ is the Haar measure on $\quoTor$ and $\lam_\bs$ is the
Haar measure on $\bS$, scaled so that $\lam_\bs(\bS) \ = \ C$.
Then locally,
\ $d\lam \ = \  d\lam_\bs\tensor d\quolam$.  Hence,
for any $\eps>0$, \quad
$\lam\lb(\maketall\Ball(\partial\sQ,\eps)\rb) \  = \ 
 \lam_\bs(\bS) \cdot 
\quolam\lb(\maketall\Ball(\partial\quosQ,\eps)\rb)
\  = \ 
 C\cdot \quolam\lb(\maketall\Ball(\partial\quosQ,\eps)\rb)$.

Part {\bf(d)} follows. \quad {\bf(e)} then follows from {\bf(c)} and {\bf(d)}.
\eclaimprf
Now, Lemma
\ref{quotient.partition}{\bf(e)} says that $\quosP_{\quotrans{}}(\tlTor)=
\tlgP$, and $(\tlgP,d_B,\Phi)$ is sensitive by hypothesis, so apply {\bf
Case 1} to conclude that $\quoPhitor$ chops $\quosP$. 
Then Claim \ref{expansive.implies.chopping.C1}{\bf(e)}
implies that $\Phitor$  chops $\sP$.
\ethmprf

\section{Injectivity of CA restricted to QS  shifts
\label{S:injective}}

  We now address the second question of Hof and Knill
\cite{HofKnill}: is a cellular automaton injective when it is
restricted to a quasisturmian shift?  We will prove:

\Theorem{\label{comeager.quasisturm.inject}}
{
Let $\Phi:\sA^\Lat\into\sA^\Lat$ be any cellular automaton.
There is a dense $G_\delta$ subset $\ISP\subset\MP$, with
$\ISP\subseteq \Phitor(\ISP)$, such that, for any $\sP\in\ISP$,
the following dichotomies hold:
\bthmlist
  \item If $\sP\in\MP$ and $\mu=\QM{\sP}$,
 then  {either} $\Phi$ is 
 constant \muae, {or} $\Phi$ is injective \muae.

  \item If $\sP\in\OP$  and $\gP=\Qshift{\sP}$,
then {either}  $\Phi\restr{\gP}$ is
constant, {or} $\Phi\restr{\gP}$ is injective.
\ethmlist
}
  To prove Theorem \ref{comeager.quasisturm.inject}, recall that
$\sP$ is {\dfn simple} if $\sP$ has no translational symmetries
(\S\ref{app:simple.inject}).   

\Proposition{\label{lem:CA.sturm.inject}}
{
  Suppose $\sP$ and $\sQ=\Phitor(\sP)$ are both simple.
\bthmlist
  \item  If $\sP\in\MP$ and $\mu=\QM{\sP}$,  then $\Phi$ is injective \muae.

  \item  If $\sP\in\OP$ and $\gP \ = \ \Qshift{\sP}$,
then $\Phi\restr{\gP}$ is injective.
\ethmlist
}
\bthmprf      
{\bf(a)}\quad Let $\nu = \QM{\sQ}$.  Then
 Lemma \ref{induced.map.on.partitions}{\bf(c)} says that $\nu = \Phi(\mu)$.  
 Let $\gM_\Phi \ := \ \set{\bq\in\sA^\Lat}{\mbox{$\bq$ has multiple
$\Phi$-preimages in $\gP$}}$.  We must show that $\nu[\gM_\Phi]=0$.

Let
$\gM_{\Sqlat}  :=  \set{\bq\in\sA^\Lat}{\mbox{ $\bq$ has
multiple $\Sqlat$-preimages in $\Tor$}}$.
Lemma \ref{simple.partn.inj} says
 $\nu[\gM_{\Sqlat}] =  0$.

We claim that $\gM_\Phi  \subset  \gM_{\Sqlat}$.
To see this, let $\bq\in\gM_\Phi$.  Thus, there are
$\bp_1,\bp_2\in\gP$ with $\Phi(\bp_1) = \bq = \Phi(\bp_2)$.  Let
$\tort_1 := \Splat^{-1}(\bp_1)$ and $\tort_2 := \Splat^{-1}(\bp_2)$.  Now,
$\bp_1\neq\bp_2$, so Lemma \ref{simple.partn.inj} says
$\tort_1\neq\tort_2$ ($\lam$-$\as$).  But 
Lemma \ref{induced.map.on.partitions}{\bf(a)}
says $\Sqlat(\tort_1) = \bq = \Sqlat(\tort_2)$.
Thus $\bq\in  \gM_{\Sqlat}$.   

Thus, $\gM_\Phi \ \subset \ \gM_{\Sqlat}$, so  $\nu[\gM_\Phi]=0$.

{\bf(b)} 
Let $\bp_1,\bp_2\in\gP$, and suppose $\Phi(\bp_1)=\Phi(\bp_2)$.
Let $\Pinv:\gP\into\Tor$ be as in Proposition \ref{open.simple.partn.inj}(d).
Let $\tort_j=\Pinv(\bp_j)$ for for $j=1,2$, and let
$\tors = \tort_2-\tort_1$.  If $\bU := \tlTor\intsct \rot{-\tors}(\tlTor)$,
then $\lam[\bU]=1$, because $\lam[\tlTor]=1= \lam\lb[\rot{-\tors}(\tlTor)\rb]$.
Thus, if $\gQ:=\Splat(\bU)$, then
$\mu[\gQ]=1$.  We'll show that $\Phi$ is many-to-one on $\gQ$,
thus contradicting {\bf(a)}.

  If $\bq_1\in\gQ$, then $\bq_1=\Splat(\toru_1)$  for some  $\toru_1\in\bU$.
Let $\toru_2 := \toru_1 + \tors$; then $\toru_2\in\tlTor$.
 Let $\bq_2 := \Splat(\toru_2)$.  Then $\bq_2\neq\bq_1$ 
(since $\Splat$ is injective by 
Proposition \ref{open.simple.partn.inj}), but we claim  $\Phi(\bq_1)=\Phi(\bq_2)$.  

To see this, use the minimality of  $(\gP,d_C,\shift{})$ to
get $\{\ell_n\}_{n=1}^\oo\subset\Lat$ with
$\D \dClim_{n\goto\oo}\shift{\ell_n}(\bp_1) \ = \ \bq_1$. 
Now, $(\gP,d_C)$ is compact, so 
drop to a subsequence such that $\{\shift{\ell_n}(\bp_2)\}_{n=1}^\oo$
converges in $\gP$.

\Claim{\label{lem:CA.sturm.inject.C2}
$\D \dClim_{n\goto\oo}\shift{\ell_n}(\bp_2) \ = \ \bq_2$.}
\bclaimprf 
Let $\bq'_2 \ := \ \D \dClim_{n\goto\oo}\shift{\ell_n}(\bp_2)$. To see that
$\bq'_2=\bq_2$,  first note that
\begin{eqnarray}
\nonumber
 \lim_{n\goto\oo}\trans{\ell_n}(\tort_1) &=&
 \lim_{n\goto\oo}\trans{\ell_n}(\Pinv(\bp_1)) \quad \eeequals{(*)} \quad
 \lim_{n\goto\oo}\Pinv(\shift{\ell_n}(\bp_1))
 \\&  \eeequals{(\dagger)} &
 \Pinv\lb(\lim_{n\goto\oo}\shift{\ell_n}(\bp_1)\rb)
\quad=\quad
\label{lem:CA.sturm.inject.e1}
 \Pinv(\bq_1)
\quad\eeequals{(\ddagger)}\quad \toru_1,
\end{eqnarray}
  where $(*)$ is  Proposition \ref{open.simple.partn.inj}(d)[i]; \ \ 
$(\dagger)$ is because $\Pinv$ is continuous; \ \  and $(\ddagger)$ is 
Proposition \ref{open.simple.partn.inj}(d)[ii].

  If $\toru'_2 := \Pinv(\bq'_2)$, then by similar reasoning,
$\D \lim_{n\goto\oo}\trans{\ell_n}(\tort_2) \ = \ \toru'_2$.
But $\tort_2=\tort_1+\tors$, so
$\D\lim_{n\goto\oo}\trans{\ell_n}(\tort_2)
\ = \ \lim_{n\goto\oo}\trans{\ell_n}(\tort_1) + \tors
\ \eeequals{(\ref{lem:CA.sturm.inject.e1})} \ \toru_1 + \tors \ = \ \toru_2$.
  Thus, $\toru'_2 = \toru_2$.
Thus, $\bq'_2 \ \eeequals{(*)} \ \Splat(\toru'_2) \ = \
 \Splat(\toru_2) \ = \ \bq_2$, where 
$(*)$ is  Proposition \ref{open.simple.partn.inj}(d)[ii].
\eclaimprf
 \beq
\mbox{ Thus,}\quad 
 \Phi(\bq_1) &\eeequals{(\ddagger)}&
\Phi(\dClim_{n\goto\oo}\shift{\ell_n}(\bp_1)) \quad\eeequals{(*)}\quad
\dClim_{n\goto\oo}\shift{\ell_n}\lb(\Phi(\bp_1)\rb)
 \\&\eeequals{(\dagger)}&
\dClim_{n\goto\oo}\shift{\ell_n}\lb(\Phi(\bp_2)\rb)
\quad\eeequals{(*)}\quad
\Phi(\dClim_{n\goto\oo}\shift{\ell_n}(\bp_2))
\quad\eeequals{(\diamond)}\quad
  \Phi(\bq_2),
\eeq
where $(\ddagger)$ is because 
$\D \lim_{n\goto\oo}\shift{\ell_n}(\bp_1) \ = \ \bq_1$;\quad 
$(*)$ is because $\Phi$ is continuous and $\shift{}$-commuting;
\quad $(\dagger)$ is because $\Phi(\bp_1)=\Phi(\bp_2)$ by hypothesis;
\quad and $(\diamond)$ is by Claim \ref{lem:CA.sturm.inject.C2}.
Thus $\Phi(\bq_1) =  \Phi(\bq_2)$.

This works for any $\bq_1\in\gQ$.   Thus, $\Phi\restr{\gQ}$ is noninjective,
but $\mu[\gQ]=1$, which contradicts {\bf(a)}.
\ethmprf
 When are the hypotheses of Proposition \ref{lem:CA.sturm.inject}
satisfied?  An element $\tort=(t_1,\ldots,t_K)\in\Tor\cong\CO{0,1}^K$ 
is called {\dfn totally rational} if $t_1,\ldots,t_K$ are all rational numbers.
It follows:
\beqn
\label{totally.rational.finite.order}
  \mbox{For any $\tort\in\Tor$,} \
\statement{$n\cdot \tort = 0$ for some $n\in\Natur$}
\iff
\statement{$\tort$ is totally rational}.
\eeqn
  Let
$\QTor \subset\Tor$ be the set of all {\em nonzero} totally rational
elements.

\Lemma{\label{C:symmetry.rational}}
{
If $\sP\in\MP$ is a nonsimple partition, 
then it has a nontrivial symmetry in $\QTor$.}
\bthmprf
  Suppose $\Tor=\Torus{K}$.  
  Let $\bS$ be the translational symmetry group of $\sP$; \ then  $\bS$ is a
closed subgroup of $\Tor$, so there is a topological group
isomorphism $\phi:\Torus{J}\x\bS_0\into\bS$ (where $0\leq J \leq K$ and $\bS_0$ is a finite group).  If $J>0$, let
$\bF\subset\Torus{J}$ be a nontrivial finite subgroup; 
if $J=0$, then $\bS_0$ is nontrivial, so let $\bF=\bS_0$.
Let $\bF' = \phi(\bF)\subset\bS$.  Then $\bF'$ is a finite subgroup of
$\Tor$, and therefore, by (\ref{totally.rational.finite.order}),
 can only contain  totally rational elements.
\ethmprf
Let $\SP$ be the set of simple partitions; let $\TP$ be the set of
trivial partitions (so that $|\TP|=|\sA|$ is finite), and let
$ \STP \ := \ \SP\union\TP \ = \ \set{\sP\in\MP}{\mbox{$\sP$ is either simple or trivial}}$.

\Corollary{\label{simple.dense.Gdelta}}
{
$\STP$ is a dense $G_\del$ subset of $\MP$.
}
\bthmprf
   If $\torq\in\QTor$, let $\qP\subset \MP$ be all $\torq$-symmetric
   partitions.  

 {\em $\qP$ is closed in $\MP$:} \ If $\{\sP_1,\sP_2,\ldots\}$ is a sequence of
$\torq$-symmetric partitions converging to some $d_\symdif$-limit $\sP$, then
$\sP$ is also $\torq$-symmetric.  

{\em  $\qP$ is nowhere dense in $\MP$:} \ 
Observe that the $\torq$-symmetry of any $\sP\in\qP$ can be disrupted
by a small `perturbation' ---ie. by removing a small piece from $\bP_a$
and adding it to $\bP_b$ for some $a,b\in\sA$.

  Thus, $(\MP\setminus\qP)$ is open and dense in $\MP$.   
Now, Lemma \ref{C:symmetry.rational} implies that
$ \SP  =  \D \Intsct_{\torq\in\QTor} (\MP\setminus\qP)$.
But $\QTor$ is countable, so $\SP$ is a countable
intersection of open dense sets, thus, dense $G_\del$.  

Since $\TP$
is finite, it is also $G_\del$; hence $\STP$ is dense and $G_\del$.
\ethmprf
  Corollary \ref{simple.dense.Gdelta} and Proposition
\ref{lem:CA.sturm.inject} do not suffice to prove Theorem
\ref{comeager.quasisturm.inject}:  \ even if $\sP$ is
simple, $\Phitor(\sP)$ may not be.
 We need conditions to ensure that both $\sP$ and $\Phitor(\sP)$ are simple.

  If $\sB$ is a Boolean algebra of subsets of $\Tor$, then $\sB$ is
{\dfn totally simple} if no nontrivial set in $\sB$ maps to itself
under any nontrivial totally rational translation.  That is: for any
$\bB\in
\sB$, with $0<\lam(\bB)<1$, and any $\torq\in\QTor$, \quad 
 $\rot{\torq}(\bB)\ \neq \ \bB$.

\example{\label{X.totally.simple}
Let $\tora\in\CO{0,1}\cong\Torus{1}$ be irrational, and let 
$\sB$ be the Boolean algebra consisting of all finite unions of
subintervals of $\Torus{1}$ of the form $\OO{n\tora, m\tora}$ for some
$n,m\in\Zahl$ (where $n\tora$ and $m\tora$ are interpreted mod 1).  Then
$\sB$ is totally simple.  To see this, suppose $\bB\in\sB$; then
$\bB = \D\Union_{j=1}^J \OO{n_j\tora, m_j\tora}$ (for some $n_1,\ldots,n_J,m_1,\ldots,m_J\in\Zahl$).  Assume WOLOG
that $n_1\tora$ is a boundary point of $\bB$.
If $\rot{\torq}(\bB)=\bB$ for
some $\torq\in\Rat$, then 
$\rot{\torq} (n_1\tora) \ = n_j\tora$ for some $j\in\CC{1..J}$, 
which means $\torq = n_j\tora - n_1\tora = (n_j-n_1)\tora$
is an integer multiple of $\tora$.  But $\torq$ is rational and $\tora$ is
irrational, so we must have $n_j-n_1=0$; hence $\torq=0$.}

  If $\sQ=\{\bQ_a\}_{a\in\sA}$ is a partition, we write ``$\sQ\prec \sB$''
if $\bQ_a\in\sB$ for all $a\in\sA$.

If $\sP=\{\bP_a\}_{a\in\sA}$ is a partition,
let $\Boole{\sP}$ be the Boolean algebra generated by the
set $\set{\trans{\ell}(\bP_a)}{\ell\in\Lat, \ a\in\sA}$ under all
finite unions and intersections.

\Lemma{\label{totally.simple.yields.simple}}
{
{\bf(a)} \  If $\sB$ is a totally simple boolean algebra, and $\sQ\prec \sB$ is a
nontrivial partition, then $\sQ$ is simple.

  {\bf(b)} \  If $\Phi$ is any cellular automaton, then 
$\Phitor(\sP) \ \prec \ \Boole{\sP}$. 

  {\bf(c)} \   Suppose $\Boole{\sP}$ is totally simple.  If  
$\Phitor(\sP)$ is nontrivial, then $\Phitor(\sP)$ is simple.  
}
\bthmprf
{\bf(a)} \quad
 $\sQ$ is nontrivial, and $\sQ\prec\sB$, so there is
some $a\in\sA$ such that  $0<\lam(\bQ_a)<1$, and $\bQ_a\in\sB$.
If $\sQ$ is nonsimple, then Lemma \ref{C:symmetry.rational} 
yields a nontrivial symmetry $\tort\in\QTor$ with
$\rot{\tort}(\sQ)=\sQ$.  Hence $\rot{\tort}(\bQ_a)= \bQ_a$,
which contradicts the total simplicity of $\sB$.

 {\bf(b)} follows from eqn.(\ref{partition.image})
defining $\Phitor(\sP)$.  \quad {\bf(c)} then follows from
{\bf(a)} and {\bf(b)}
\ethmprf

\example{\label{X.totally.simple.2}
 Let $\Tor=\Torus{1}\cong\CO{0,1}$ and $\dL=\Zahl$.  Let
 $\tora\in\CO{0,1}$ be irrational, and let $\Zahl$ act by
$\trans{z}(\tort) = \tort+ z\tora$ as in Example \ref{X:circle.rotation}.
Let $\sA=\{0,1\}$, and let $\sP \ = \ \lb\{ \OO{0,\tora}, \ \OO{\tora,1}\rb\}$ as in Example \ref{X:sturm1}.
Then $\Boole{\sP}$ is the totally simple
Boolean algebra from Example \ref{X.totally.simple}.  Thus, if $\Phi$ is any CA, and 
$\Phitor(\sP)$ is nontrivial [eg. Example \ref{X:ledrappier}], then $\Phitor(\sP)$ is simple. }

  When is $\Boole{\sP}$ totally simple?
  If $\sP=\{\bP_a\}_{a\in\sA}$ is an open partition of $\Tor$,
 and $\tors\in\Tor$, then we say $\sP$ has {\dfn
 $\tors$-local symmetry} if there are $a,b\in\sA$ and open sets
$\bO,\bO'\subset\Tor$ so
 that $(\partial\bP_a)\intsct\bO \neq \emptyset \neq
 (\partial\bP_b)\intsct\bO'$, and
 $\rot{\tors}\lb((\partial\bP_a)\intsct\bO\rb) \ = \
 (\partial\bP_b)\intsct\bO'$.
 We say $\sP$ is {\dfn primitive} if $\sP$ has no nontrivial $\tors$-local
symmetries for any $\tors\in\
\LQTor \ := \ \set{\trans{\ell}(\torq)}{\ell\in\Lat, \ \torq\in\QTor}$.

\example{Let $\trans{}$ and $\sP$ be 
as in Example \ref{X.totally.simple.2}.
Then $\sP$ has $\tora$-local symmetry,
but is still primitive,
because $\tora\not\in\Zahl\!\QTor
 =  \set{z\tora + \torq}{z\in\Zahl, \ \torq\in \QTor}$
(because $0\not\in\QTor$).
}

\Lemma{\label{primitive.yields.totally.simple}}
{
  If $\sP\in\OP$ is primitive, then  $\Boole{\sP}$ is totally simple.
}
\bthmprf
  For any $a\in\sA$ and $\ell\in\Lat$, let $\bP_a^\ell = \trans{\ell}(\bP_a)$.
Then any $\bU \in \Boole{\sP}$ has the form
\[
  \bU \quad = \quad \Union_{j=1}^J \, \Intsct_{k=1}^K  \bP_{a_{jk}}^{\ell_{jk}},
\]
 for some $J,K>0$ and some collections
 $\{a_{jk}\}_{j=1}^J\,_{k=1}^K \subset \sA$ and
 $\{\ell_{jk}\}_{j=1}^J\,_{k=1}^K \subset \Lat$.  Suppose $\rot{\torq}(\bU)=\bU$
for some $\torq\in\QTor$.  Then $\rot{\torq}(\partial\bU) = \partial\bU$.  
But 
\beq
\partial\bU \quad  \subset \quad 
\Union_{j=1}^J  \, \partial \lb(\Intsct_{k=1}^K \bP_{a_{jk}}^{\ell_{jk}}\rb),
&& \\ \mbox{and,  for all $ \ j\in\CC{1..J}$}, \quad
 \partial \lb(\Intsct_{k=1}^K \bP_{a_{jk}}^{\ell_{jk}}\rb) \
& = &
 \Union_{k_*=1}^K \, \lb(\lb(\partial  \bP_{a_{jk_*}}^{\ell_{jk_*}}\rb) \intsct
\Intsct_{k_*\neq k =1}^K \bP_{a_{jk}}^{\ell_{jk}}\rb).
\eeq
  Hence, we must have open sets $\bV,\bV'\subset\Tor$ such that 
\beqn
\label{primitive.yields.totally.simple.e1}
  \rot{\torq} \lb( \bV\intsct \lb(\partial  \bP_{a_{jk_*}}^{\ell_{jk_*}}\rb) \intsct
\Intsct_{k_*\neq k =1}^K \bP_{a_{jk}}^{\ell_{jk}}\rb)
\quad=\quad
\bV'\intsct  \lb(\partial  \bP_{a_{j'k_*'}}^{\ell_{j'k_*'}}\rb) \intsct
\Intsct_{k_*\neq k =1}^K \bP_{a_{j'k}}^{\ell_{j'k}},
\eeqn
 for some $j,j'\in\CC{1..J}$ and $k_*,k_*'\in\CC{1..K}$.
Let $a:=a_{jk_*}$ and $a' := a_{j'k_*}$;  let $\ell := \ell_{jk_*}$
and $\ell' := \ell_{j'k_*'}$.
If $\bU := \bV\intsct \D\Intsct_{k_*\neq k =1}^K \bP_{a_{jk}}^{\ell_{jk}}$
and $\bU' := \bV'\intsct  \D\Intsct_{k_*\neq k =1}^K \bP_{a_{j'k}}^{\ell_{j'k}}$,
then $\bU$ and $\bU'$
are open sets, and 
(\ref{primitive.yields.totally.simple.e1})
becomes: \ 
 $ \rot{\torq} \lb( (\partial  \bP_{a}^{\ell}) \intsct \bU\rb)
\ = \ 
 (\partial  \bP_{a'}^{\ell'}) \intsct \bU'$.
\ Equivalently,
 \  $ \rot{\torq} \lb(\trans{\ell}(\partial  \bP_{a}) \intsct \bU\rb)
 \ = \ 
\trans{\ell'} (\partial  \bP_{a'}) \intsct \bU'$.
\ Equivalently,
 \ $ \rot{\torq} \circ \trans{\ell-\ell'}
\lb((\partial  \bP_{a})  \intsct \bO \rb)
\ = \  (\partial  \bP_{a'}) \intsct \bO'$,
where $\bO := \trans{-\ell}(\bU)$ and
 $\bO' := \trans{-\ell'}(\bU')$.
\ Equivalently, \ $ 
\rot{\tors} \lb((\partial  \bP_{a}) \intsct \bO \rb)\ = \ 
 (\partial  \bP_{a'}) \intsct \bO'$, 
where $\tors := \ \trans{\ell-\ell'}(\torq) \ \in \ \LQTor$.
 Thus, $\sP$ has an $\tors$-local
symmetry, contradicting our hypothesis.
\ethmprf

\Corollary{\label{primitive.into.simple}}
{
 Let $\PP\subset\OP$ be the set of primitive partitions.  Then  $\PP \ \subset \ \Phitor^{-1}\lb(\STP\rb)$.
}
\bthmprf
If $\sP\in\PP$, then Lemma \ref{primitive.yields.totally.simple}
says $\Boole{\sP}$ is totally simple.  Hence,
Lemma \ref{totally.simple.yields.simple}{\bf(c)}
 says that either $\Phitor(\sP)$ is 
trivial or  $\Phitor(\sP)$ is
simple.  Hence, $\Phitor(\PP) \ \subset \ \STP$.
\ethmprf

\Lemma{\label{primitive.dense}}
{ 
$\PP$ is $d_\symdif$-dense in $\MP$.
}
\bthmprf 
  Let $\sM\in\MP$ be a measurable partition and let $\eps>0$.
 Lemma \ref{DP.dense.in.MP} 
says  $\DP$ is $d_\symdif$-dense in $\MP$,
so there is a dyadic partition $\sD\closeto{\eps/2}\sM$.
We claim there is a primitive partition $\sP \closeto{\eps/2}\sD$.
Thus,  $\sP \closeto{\eps/2}\sD \closeto{\eps/2} \sM$, so 
$\sP \closeto{\eps}\sM$, so we're done.

To construct the primitive  partition $\sP \closeto{\eps/2}\sD$,
suppose $\Tor = \Torus{K}$, and consider
two cases:

\begin{figure}
\psfrag{MM}[][]{$\sM$}
\psfrag{DD}[][]{$\sD$}
\psfrag{PP}[][]{$\sP$}
\psfrag{~~}[][]{$\closeto{\eps/2}$}
\centerline{\includegraphics[scale=0.8]{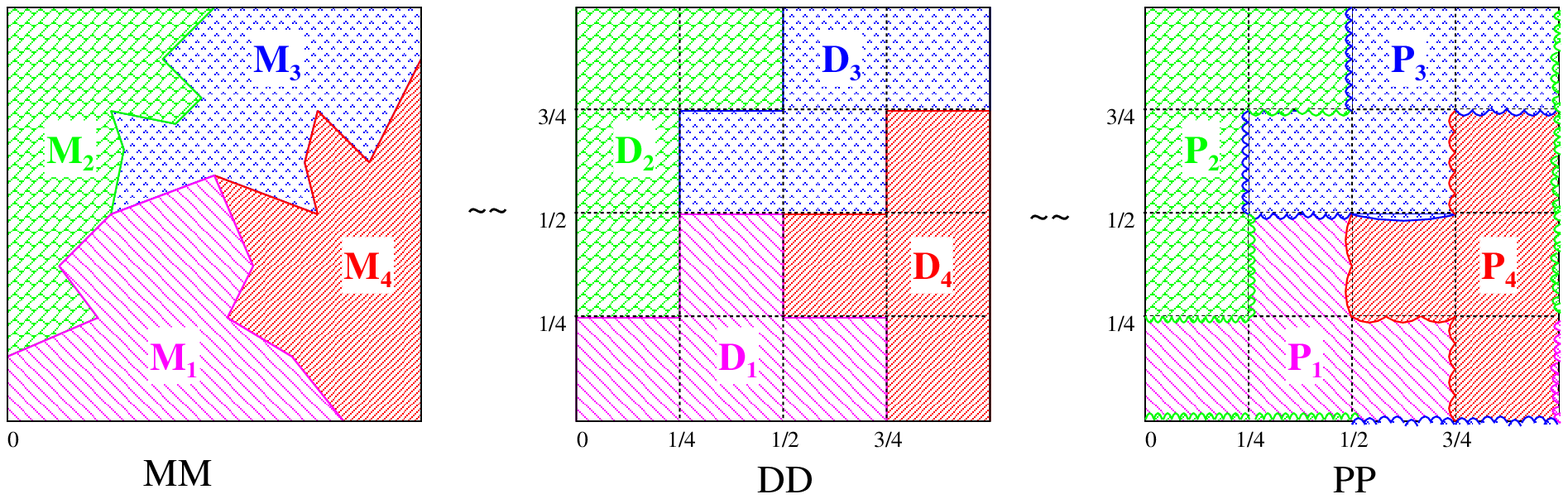}}
\caption{\footnotesize Approximating partition $\sM$ with a primitive
partition $\sP$.\label{fig:wavy}}
\end{figure}

{\em Case $K=1$:} \quad $\Tor=\Torus{1}\cong\CO{0,1}$, so
 $\partial\sD$ is a countable set of points.
By slightly perturbing these points, we can construct $\sP\in\OP$ with 
 $\sP\closeto{\eps/2}\sD$
such that, for any $\torb_1,\torb_2\in\partial\sP$, \ $\torb_1-\torb_2 \not\in
\LQTor$.  Hence, $\sP$ is primitive.

{\em Case $K>1$:} \quad 
  First, observe that, if $\tors\in\Tor$ and
 $\sD$ has a $\tors$-local symmetry,
then $\tors$ must be a dyadic element of $\Tor$.  Thus, it
suffices to `perturb' $\sD$ so as to disrupt all dyadic local symmetries,
without introducing any other local symmetries. 

  If $\bD_a\in\sD$, then $\partial \bD_a$ is a finite disjoint union of
$(K-1)$-dimensional dyadic cubic faces.  Let $\sK_a$ be the set
of these faces.  Let $\sK = \Union_{a\in\sA} \sK_a$; then $\sK$ has 
is a finite collection of dyadic cubic faces.
Now, let $\{\alp_\kap\}_{\kap\in\sK}\in\CC{0,1}$ be distinct frequencies.
If $\sD$ is $n$-dyadic, then let $\del < \frac{1}{2^{n+1}}$.  
For each $a\in\sA$, let $\sP_a$ be the set obtained as follows:
begin with $\sD_a$, and `corrugate' each face $\kap\in\sK_a$,
with a sine wave of frequency $\alp_\kap$ and amplitude $\del$
(Fig.\ref{fig:wavy})

Since all the frequencies are distinct, there can by no local symmetry
from any face to any other; hence $\sP$ is primitive.  If
$\del$ is small enough, then $\sP \closeto{\eps/2}\sD$.
\ethmprf

\Corollary{\label{preimage.of.SP.is.dense.GD}}
{
  For any $n\in\Natur$, \ $\Phitor^{-n}\lb(\STP\rb)$ is a dense $G_\del$ subset
of $(\MP,d_\symdif)$.
}
\bthmprf 
{\em G$\del$:}\quad 
 Corollary \ref{simple.dense.Gdelta} says $\STP$ is a $G_\del$ subset
of $\MP$.  Proposition \ref{Phi.star.continuous} says $\Phitor^n$ is
$d_\symdif$-continuous, so the $\Phitor^n$-preimage of any open set is open; 
hence,  $\Phitor^{-n}\lb(\STP\rb)$ is also $G_\del$.

{\em Dense:} \quad Corollary \ref{primitive.into.simple} says $\PP\subset
\Phitor^{-n}\lb(\STP\rb)$, while Lemma \ref{primitive.dense} says that
$\PP$ is dense in $\MP$; hence $\Phitor^{-n}\lb(\STP\rb)$ is also dense
in $\MP$.
\ethmprf

\bthmprf[Proof of Theorem \ref{comeager.quasisturm.inject}]
First, recall that $(\MP,d_\symdif)$ is a complete metric space 
(Proposition \ref{MP.complete}), and thus, a Baire space \cite[Corollary 25.4(b), p.186]{Willard}.

  Now, let $\ISP \ := \ \D \Intsct_{n=1}^\oo \Phitor^{-n} \lb(\STP\rb)$.
Then $\ISP$ is a countable intersection of dense $G_\del$ subsets
of $\MP$ (Corollary \ref{preimage.of.SP.is.dense.GD}), and thus, is
itself a dense $G_\del$ subset (because $\MP$ is Baire).
  By construction, $\Phitor(\ISP) \supseteq \ISP$.

If $\sP\in\ISP$ is nontrivial, then $\sP$ is simple, and 
either $\Phitor(\sP)$
is trivial (ie. $\Phi\restr{\gP}$ is constant) or $\Phitor(\sP)$
is also simple.  Now apply Proposition \ref{lem:CA.sturm.inject}{\bf(a)} and
{\bf(b)} respectively to get Theorem
\ref{comeager.quasisturm.inject}{\bf(a)} and {\bf(b)} respectively.
\ethmprf

\section{\label{app:custom.quasi}
Customized Quasisturmian Systems\protect\footnotemark}

\footnotetext{This section contains technical results which are used in 
\S\ref{S:surject} and \S\ref{S:fixedpoints}}

\begin{figure}[hbtf]
\centerline{\includegraphics[scale=0.65]{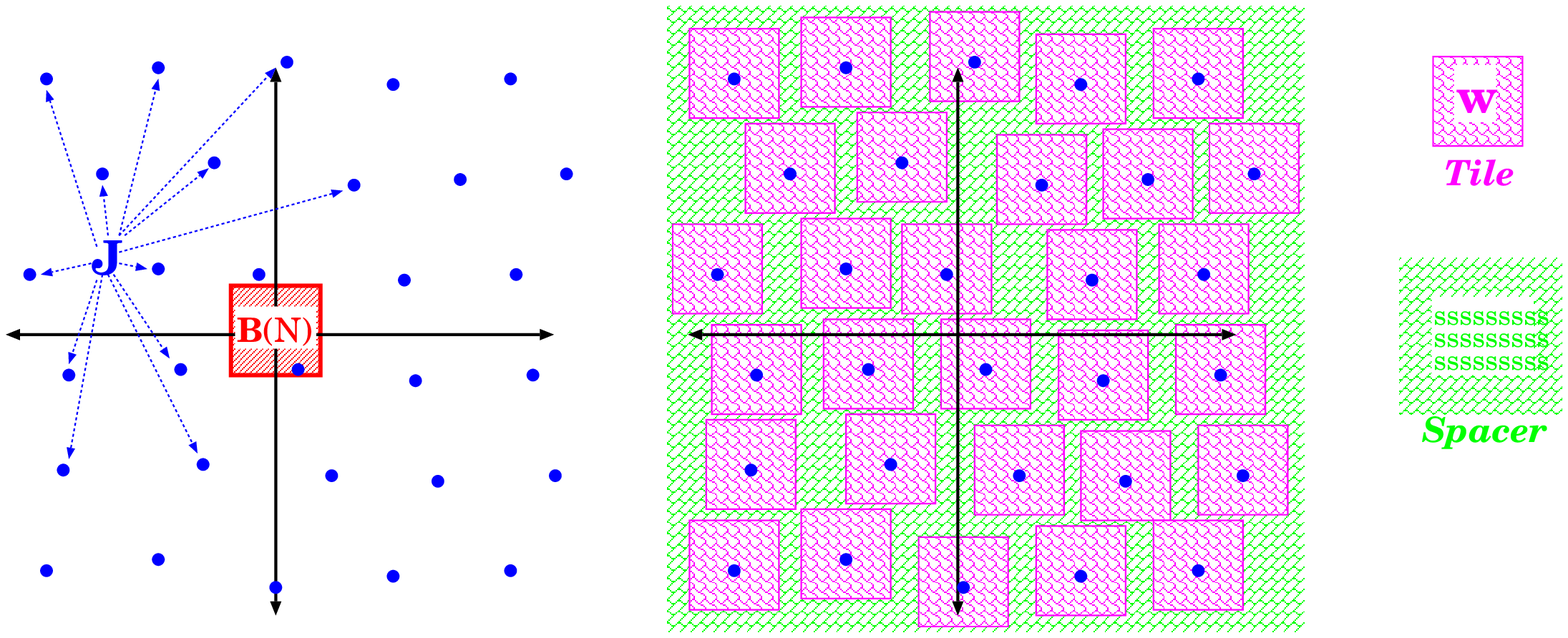}}
\caption{\footnotesize $\ba$ is $\eps$-tiled with tile $\bw$, with
skeleton $\dJ$ and spacer $s$.\label{fig:tiling}}
\end{figure}

 Let $\dB(N):=\CO{-N..N}^D\subset\Lat$.  Fix 
a word $\tile\in\sA^{\dB(N)}$ and a letter $s\in\sA$.
If $\bx\in\sA^\Lat$, and $0<\eps<1$, then we say $\bx$ is {\dfn $\eps$-tiled}
by $\tile$  with {\dfn spacer} $s$ if there is a subset $\dJ\subset\dL$ 
(the {\dfn skeleton} of the tiling) such that, as shown in 
Figure \ref{fig:tiling}:
\bdesc
  \item[(T1)] $\lb(\maketall\fj_1+\dB(N)\rb) \intsct 
\lb(\maketall\fj_2+\dB(N)\rb) \ = \ \emptyset$,
 for any distinct $\fj_1,\fj_2 \in \dJ$.

  \item[(T2)] $\density{\dJ}> \D \frac{1-\eps}{(2N)^D}$.  Thus,
 $\density{\dB(N)+\dJ \maketall } \ = \ (2N)^D \cdot \density{\dJ}
 \ > \ 1-\eps$.
  
  \item[(T3)] For every $\fj\in\dJ$, \quad $\bx\restr{\fj+\dB(N)} \ = \ \tile$.

  \item[(T4)] For any $\ell\not\in (\dB(N)+\dJ)$, \quad $x_\ell=s$.
\edesc

\Proposition{\label{painting.the.city}}
{
 Let $N\in\Natur$ and $\eps>0$. For any tile $\tile\in\sA^{\dB(N)}$
and  `spacer'  $s\in\sA$,  there is an open 
partition $\sP\in\OP$ such that every element
of $\Qshift{\sP}$ is $\eps$-tiled by  $\tile$ with spacer $s$.\qed
}

To prove Proposition  \ref{painting.the.city}, we use a $\Zahl^D$-action version of the Rokhlin-Kakutani-Halmos Lemma:

\Lemma{\label{Rokhlin.Lemma} \cite{OrnsteinWeissRokhlin}
}
{
 For any $N\in\Natur$ and $\eps>0$, there exists an open subset
$\bJ\subset\Tor$ such that:
\bthmlist
  \item The sets $\lb\{\trans{\fb}(\bJ)\rb\}_{\fb\in\dB(N)}$ are disjoint.
  \item $\D \lam(\bJ) \ > \  \frac{1-\eps}{(2N)^D}$, and thus,
$\D \lam\lb(\maketall\rb.\Disj_{\fb\in\dB(N)} \trans{\fb}(\bJ)
\lb.\maketall\rb) \ > \ 1-\eps$.
\qed
\ethmlist
}

  If $\Lat=\Zahl$, then $\dB(N) = \CO{-N..N}$, and
 $\lb\{\trans{\fb}(\bJ)\rb\}_{\fb\in\dB(N)}$ is called 
an {\dfn $(\eps,N)$-Rokhlin tower} of {\dfn height} $2N$, with {\dfn base} $\bJ$;
the sets  $\trans{\fb}(\bJ)$ are the {\dfn levels} of the tower.
 If $\Lat=\Zahl^D$ for $D\geq 2$, then we call the structure
$\lb\{\trans{\fb}(\bJ)\rb\}_{\fb\in\dB(N)}$ an {\dfn
$(\eps,N)$-Rokhlin city}, with {\dfn base} $\bJ$.  The sets
$\trans{\fb}(\bJ)$ are the {\dfn houses} of the city. 

 Think of the elements of $\sA$ as `colours'; then we can build an
$\sA$-labelled partition by `painting' the houses of the Rokhlin 
city in different $\sA$-colours, as follows.  Fix a word $\bw\in\sA^{\dB(N)}$,
and  a `spacer' letter $s\in\sA$.
Let $\bS := int(\D \Tor \setminus \Disj_{\fb\in\dB(N)} \trans{\fb}(\bJ))$.
If $\tile = [\til_\fb]_{\fb\in\dB(N)}$, then
for any $a\in\sA$, let  $\dP_a : =  \set{\fb\in\dB(N)}{\til_\fb \ = \  a}$.
Define partition $\sP:=\{\bP_a\}_{a\in\sA}$, where
\beqn
\label{painting.the.city.eqn}
\bP_s \quad := \quad \bS  \ \disj \ 
 \Disj_{\fp\in\dP_s} \trans{\fp}(\bJ),
\qquad
\mbox{and, for any $a\in\sA\setminus\{s\}$},
\qquad \bP_a \quad := \quad \Disj_{\fp\in\dP_a} \trans{\fp}(\bJ).
\eeqn
  We say  $\sP$ is obtained by {\dfn painting} 
Rokhlin city $\lb\{\trans{\fb}(\bJ)\rb\}_{\fb\in\dB(N)}$ with 
word $\bw$ and spacer $s$.

\Lemma{\label{painting.the.city.2}}
{
  Let $N\in\Natur$ and $\eps>0$, and
 let $\lb\{\trans{\fb}(\bJ)\rb\}_{\fb\in\dB(N)}$ be an 
$(\eps,N)$-Rokhlin city.  Suppose
$\sP$ is obtained by painting  $\lb\{\trans{\fb}(\bJ)\rb\}_{\fb\in\dB(N)}$
with $\bw$ and $s$, as in eqn.{\rm(\ref{painting.the.city.eqn})}.

If $\tort\in\Tor$ and
$\bp=\Splat(\tort)$, then $\bp$ is $\eps$-tiled by $\tile$
with skeleton $\dJ \ = \ \set{\ell\in\Lat}{\trans{\ell}(\tort) \in \bJ}$.
}
\bthmprf \!\!\!The Generalized Ergodic Theorem says
$\density{\dJ}  =  \lam(\bJ)   >  \frac{1-\eps}{(2D)^N}$.
Also, $\forall \ \fj\in\dJ$, \ 
$\bp\restr{\fj+\dB(N)}  =  \tile$.  Finally,
$\statement{$\ell\not\in\dJ+\dB(N)$}
\Rightarrow\statement{$\trans{\ell}(\tort)\in\bS$}
\Rightarrow\statement{$p_\ell=s$}$.
\ethmprf

\bthmprf[Proof of Proposition \ref{painting.the.city}:] 
 Let $\lb\{\trans{\fb}(\bJ)\rb\}_{\fb\in\dB(N)}$ be an 
 $(\eps,N)$-Rokhlin city, provided by Lemma \ref{Rokhlin.Lemma}.
Now paint the city with $\tile$ and $s$, and 
apply Lemma \ref{painting.the.city.2}.
\ethmprf

\Corollary{\label{QP.is.dense}}
{
 $\QS$ is dense in $\sA^\Lat$ in the Cantor topology.
}
\bthmprf
  Fix $\ba\in\sA^\Lat$ and $\eps>0$.  Let $\tort\in\Tor$;
 we'll build a partition $\sQ\in\MP$ so
that $\bq=\Sqlat(\tort)$ is $\eps$-close to $\ba$ in the Cantor metric. 
Let $N > -\log_2(\eps)$,
and let $\tile = \ba\restr{\dB(N)}$.  Use Proposition \ref{painting.the.city}
to find $\sP\in\OP$ such that $\bp=\Splat(\tort)$ is tiled by copies of
$\tile$.  Thus, there is some $\ell\in\dL$ such that
 $\bp\restr{\ell+\dB(N)} = \tile$.  Let $\sQ = \trans{\ell}(\sP)$; thus,
if $\bq=\Sqlat(\tort)$, then  $\bq = \shift{\ell}(\bp)$ by 
Proposition \ref{splat.homeo}(a),
so that $\bq\restr{\dB(N)} = \tile = \ba\restr{\dB(N)}$,
so that $d_C(\bq,\ba) <  2^{-N}  <  \eps$.
\ethmprf

Let $\ergM$ be the space of $\shift{}$-ergodic
probability measures on $\sA^\Lat$, with the weak* topology induced by
convergence along cylinder sets.  

\Corollary{\label{MQP.is.wkst.dense}}
{
  $\QSM$ is weak*-dense in $\ergM$.
}
\bthmprf
  Let $\nu\in\ergM$.  Fix $N>0$;  let $\bc_1,\ldots,\bc_J\in\sA^{\dB(N)}$,
and for each $j\in\CC{1..J}$, let $\gC_j:=\set{\ba\in\sA^\Lat}{\ba\restr{\dB(N)}=\bc_j}$ be the corresponding cylinder sets.  For any $\eps>0$, we'll
construct a quasisturmian measure $\mu\in\QSM$ such that
$\mu[\gC_j] \closeto{2\eps} \nu[\gC_j]$ for all $j\in\CC{1..J}$.

  The Generalized Ergodic Theorem yields some $\ba\in\sA^\Lat$ which
is $(\shift{},\nu)$-generic for $\gC_1,\ldots,\gC_J$.  Thus, if
$M$ is large enough, then 
\beqn
\label{MQP.is.wkst.dense.e1}
 \mbox{For all $j\in\CC{1..J}$,}\qquad
  \nu[\gC_j] \quad \closeto{\,\eps\,} \quad
\frac{1}{(2M)^D} \sum_{\fb\in\dB(M)} \ 
\chr{j}\lb(\maketall\shift{\fb}(\ba)\rb),
\eeqn
  where $\chr{j}$ is the characteristic function of $\gC_j$.
Now, $\chr{j}(\ba)$ is a function only of $\ba\restr{\dB(N)}$, so
the sum (\ref{MQP.is.wkst.dense.e1}) is a function only of 
$\bw:=\ba\restr{\dB(N+M)}$.
Use Proposition \ref{painting.the.city} to find some $\sP\in\OP$ so
that $\bp\in\Qshift{\sP}$ is $\eps$-tiled with $\bw$.  
If $\mu=\QM{\sP}$, then $\bp$ is $(\shift{},\mu)$-generic, so
for for all $ \ j\in\CC{1..J}$,

$\begin{array}[b]{rcl}
\mu[\gC_j]
&=&
\D \lim_{K\goto\oo} \ \frac{1}{(2K)^D} \sum_{\fb\in\dB(K)} \ 
\chr{j}\lb(\maketall\shift{\fb}(\bp)\rb)
\ \ \closeto{\,\eps\,}\ \
\frac{1}{(2M)^D} \sum_{\fb\in\dB(M)} \ 
\chr{j}\lb(\maketall\shift{\fb}(\bw)\rb)
\\ &=&
\D \frac{1}{(2M)^D} \sum_{\fb\in\dB(M)} \ 
\chr{j}\lb(\maketall\shift{\fb}(\ba)\rb)
\quad\closeto{\,\eps\,}\quad
\nu[\gC_j]. \ \hfill 
\end{array}$\nolinebreak\hspace{-0.7em}
 \ethmprf

\Lemma{\label{similar.tilings}}
{
Let $\bx,\bx'\in\sA^\Lat$ with $d_B(\bx,\bx')\leq \del$.
Suppose $N>0$, and that $\bx$ can be
$\eps$-tiled by some  $\bw\in\sA^{\dB(N)}$,
 with skeleton $\dJ\subset\Lat $.  Let $\del<(1-\eps)/(2N)^D$, so 
$\eps': = \eps+ (2N)^D \del<1$.  Then $\bx'$
can be $\eps'$-tiled by $\bw$, with skeleton $\dJ'\subset\Lat$,
such that $\density{\dJ\intsct\dJ'} \ \geq \  \density{\dJ}-\del$.
}
\bthmprf  Let $\dJ' = \set{\ell\in\Lat}{\bx'\restr{\ell+\dB(N)} \ = \ \bw}$,
and let $ \dJ_\Del \ = \ \dJ\setminus\dJ'$.  Thus, $\dJ\intsct\dJ'
\ = \ \dJ\setminus\dJ_\Del$.  

We'll show that $\density{\dJ_\Del} \ \leq \ \del$.  
Thus, $\density{\dJ\intsct\dJ'}
\ = \ \density{\dJ} - \density{\dJ_\Del}
\ \geq \ \density{\dJ} - \del$.

Let
 $\dK = \set{\ell\in\Lat}{x_\ell\neq x'_\ell}$.  Then
 $ \density{\dK}  =  d_B(\bx,\bx')  =  \del$, and
for any $\fj\in\dJ$,
\beq
\statement{$\fj\in\dJ_\Del$}
 &\iff&
\statement{$\bx'\restr{\fj+\dB(N)} \neq \bw$}
\quad\iff\quad
\statement{$\bx'\restr{\fj+\dB(N)} \neq \bx\restr{\fj+\dB(N)}$}
\\&\iff&
\statement{$x'_{\fj+\fb} \neq x_{\fj+\fb}$, for some $\fb \in \dB(N)$}
\\&\iff&
\statement{$\fj+\fb \in \dK$,  for some $\fb\in\dB(N)$}.
\eeq
  Thus, we can define a function $\bet:\dJ_\Del\into\dB(N)$ such that
$\fj+\bet(\fj)\in\dK$ for all $\fj\in\dJ_\Del$.  This defines
a function $\kap:\dJ_\Del\into\dK$ by $\kap(\fj) = \fj+\bet(\fj)$.

Note that $\kap$ is an injection, because
for any distinct $\fj_1,\fj_2\in\dJ_\Del\subset\dJ$,
tiling condition {\bf(T1)} says that
 $(\fj_1+\dB(N)) \intsct (\fj_2+\dB(N)) \ = \ \emptyset$;
 hence 
$\fj_1+\bet(\fj_1)\ \neq \ \fj_2+\bet(\fj_2)$.

\claim{Let $\dJ_\kappa \ = \ \kap(\dJ_\Del) \ \subseteq \ \dK$. Then
$\density{\dJ_\kappa}  \ = \  \density{\dJ_\Del}$.}
\bclaimprf
  For any $M>0$, every element of 
$\dJ_\Del\intsct \dB(M)$ must go to some element of $\dJ_\kappa\intsct \dB(M+N)$ under $\kappa$.  Likewise, 
every element of $\dJ_\kappa\intsct \dB(M-N)$ must `come from' 
$\dJ_\Del\intsct \dB(M)$.  Thus, \
$ \#\lb[\dJ_\kappa\intsct \dB(M-N)\rb] \ \leeeq{(a)} \
  \#\lb[\dJ_\Del\intsct \dB(M)\rb] \ \leeeq{(b)} \
  \#\lb[\dJ_\kappa\intsct \dB(M+N)\rb]$.
\beq
 \mbox{Thus,}\quad
 \density{\dJ_\kappa} &=&
\lim_{M\goto\oo}\,  \frac{\#\lb[\dJ_\kappa\intsct \dB(M-N)\rb]}{(2M-2N)^D}
\\ &=&
\lb(\lim_{M\goto\oo}\, \frac{(2M)^D}{(2M-2N)^D}\rb)\cdot\lb( \lim_{M\goto\oo}\,
  \frac{\#\lb[\dJ_\kappa\intsct \dB(M-N)\rb]}{(2M)^D}\rb)
\\&=&
 \lim_{M\goto\oo}\,
  \frac{\#\lb[\dJ_\kappa\intsct \dB(M-N)\rb]}{(2M)^D}
\ \ \leeeq{(A)} \ \ 
\lim_{M\goto\oo}\,  \frac{\#\lb[\dJ_\Del\intsct \dB(M)\rb]}{(2M)^D}
\\&=&
 \density{\dJ_\Del}
\quad\leeeq{(B)}\quad
 \lim_{M\goto\oo}\,  \frac{\#\lb[\dJ_\kap\intsct \dB(M+N)\rb]}{(2M)^D}
\\&=&
\lb(\lim_{M\goto\oo}\, \frac{(2M+2N)^D}{(2M)^D}\rb)\cdot \lb(\lim_{M\goto\oo}\,
  \frac{\#\lb[\dJ_\kappa\intsct \dB(M+N)\rb]}{(2M+2N)^D}\rb)
\\&=&
\lim_{M\goto\oo}\,  \frac{\#\lb[\dJ_\kappa\intsct \dB(M+N)\rb]}{(2M+2N)^D}
\quad=\quad
 \maketall \density{\dJ_\kappa},
\eeq
where {\bf(A)} is by {\bf(a)} and {\bf(B)} is by {\bf(b)}
\eclaimprf
Thus, Claim 1 implies $\density{\dJ_\Del}  =  \density{\dJ_\kappa}
  \leq  \density{\dK}   =  \del$, as desired.\nolinebreak
\ethmprf

\Corollary{\label{dB.dense.implies.dC.dense}}
{
  If $\gX\subset\QS$ is $\shift{}$-invariant and
$d_B$-dense in $\QS$, then 
$\gX$ is also $d_C$-dense in $\QS$.
}
\bthmprf
  Let $\bq\in\QS$ and $N>0$; we want $\bx\in\gX$ such that $\bx\restr{\dB(N)} =  \bq\restr{\dB(N)}$.  Let $\bw:=\bq\restr{\dB(N)}$. If $\eps>0$, then
Proposition \ref{painting.the.city} yields some $\bq'\in\QS$ which
is $\eps$-tiled by $\bw$.  Let $\del< (1-\eps)/(2N)^D$, and use
 the $d_B$-density
of $\QS$ to find  $\bx'\in\gX$ with $d_B(\bx',\bq')<\del$.
Then Lemma \ref{similar.tilings} says that $\bx'$ can be $\eps'$-tiled
by $\bw$.  Thus, there is some
$\ell\in\Lat$ such that $\bx'\restr{\dB(N)+\ell}  =  \bw$. 
If $\bx=\shift{-\ell}(\bx')$, then $\bx\restr{\dB(N)} = \bw$, and
$\bx\in\gX$ because $\gX$ is $\shift{}$-invariant.
\ethmprf

\section{Surjectivity and image density
\label{S:surject}}

A natural conjecture: {\em If $\Phi:\sA^\Lat\into\sA^\Lat$ is
surjective, then $\Phitor:\MP\into\MP$ is also surjective.}
Unfortunately, this is false.  For example, suppose
$\sA=\{0,1\}=\Zahlmod{2}$.   For any $\sP\in\MP$, if
$\sP=\{\bP_0,\bP_1\}$, then let $\barsP := \{\barbP_0,\barbP_1\}$, where
$\barbP_0:=\bP_1$ and $\barbP_1:=\bP_0$.

\Proposition{}
{
Let $\Phi$ be as in {\rm Example
\ref{X:ledrappier}}.  If $\sP\in\Phitor(\MP)$, then 
$\barsP\not\in\Phitor(\MP)$.
}
\bthmprf
 First let $\sU\in\OP$
be the `unity' partition $\sU\equiv 1$, ie: $\bU_0  =  \emptyset$, \ 
$\bU_1  =  \Tor$.  We claim $\sU\not\in\Phitor(\MP)$.  To
see this, suppose $\sQ\in\MP$ 
and $\sU = \Phitor(\sQ)$.
Treating $\sU$ and $\sQ$ as functions from $\Tor$ to $\sA$, we have:
$ \sU(\tort)  = \  \sQ(\tort) + \sQ\lb(\maketall\trans{}(\tort)\rb)$
(mod $2$) for all $\tort\in\Tor$ ---ie. $\sU \ = \ \sQ + \sQ\circ\trans{}$
(mod $2$).
 But  $\sU\equiv 1$, so $1   = \ \sQ + \sQ\circ\trans{}$,
which means
$\sQ\circ\trans{}  = \ 1 - \sQ$.  Thus,
$\sQ\circ \trans{2} \ = \ \sQ\circ \trans{} \circ\trans{} \ = \
 (1 - \sQ)\circ\trans{}
\ = \  1 - \sQ\circ\trans{}
 \ = \ 1 - (1-\sQ) = \sQ$.
Thus, $\sQ$ is $\trans{2}$-invariant.  But $\trans{}$ is totally ergodic,
so this means that $\sQ$ is a constant ---either $\sQ\equiv 1$ or
$\sQ\equiv 0$.  Neither of these partitions maps to $\sU$
under $\Phitor$, so $\sU\not\in\Phitor(\MP)$.
 
    Now, let $\sP\in\MP$ and suppose $\sP=\Phitor(\sQ)$ and
$\barsP=\Phitor(\sQ')$ for some $\sQ,\sQ'\in\MP$.
 Note that $\sP+\barsP\equiv\sU$.  Let $\sQ^\dagger:=\sQ+\sQ'$; then 
$\Phitor(\sQ^\dagger) \ \eeequals{(*)}  \ \Phitor(\sQ)+\Phitor(\sQ')
\ = \ \sP + \barsP = \sU$, where $(*)$ is because $\Phitor$ 
is $\Zahlmod{2}$-linear.
But we know that $\sU\not\in\Phitor(\MP)$. Contradiction.
\ethmprf
 Thus, $\Phitor$ is not surjective:  $\Phitor(\MP)$
fills at most `half' of $\MP$. Nevertheless, we will prove:

\Theorem{\label{CA.surjective.on.MP}}
{
  Let $\Phi:\sA^\Zahl\into\sA^\Zahl$ be any CA.  The
following are equivalent:

  {\bf(a)} \ $\Phi$ is surjective onto $\sA^\Zahl$. 

  {\bf(b)} \ $\Phitor(\MP)$ is  $d_\symdif$-dense in $\MP$;
and $\Phitor(\OP)$ is $d_\symdif$-dense in $\OP$.

  {\bf(c)} \  $\Phi(\QS)$ is $d_B$-dense in $\QS$.

  {\bf(d)} \  $\Phi(\QS)$ is $d_C$-dense in $\QS$.
}
If $\dL=\Zahl$ (as in Theorem \ref{CA.surjective.on.MP}), then there is
some irrational $\tora\in\Tor$ such that $\trans{z}(\tort) = \tort +
z\tora$ for any $z\in\Zahl$ and $\tort\in\Tor$.  The system
$(\Tor,\varsigma,\lam)$ is called an {\dfn irrational rotation}.  To
prove Theorem \ref{CA.surjective.on.MP}, we'll use the {\em rank
one} property of irrational rotations.

\Theorem{(del Junco) \cite{delJQuasiRank} \label{delJQuasiRankThm}}
{

  Any irrational rotation is {\dfn topologically rank one}. 
 That is, there is a sequence
$\{\bJ_i\}_{i=1}^\oo$ of open subsets of $\Tor$ such that:

 {\bf(a)} \  $\bJ_i$ is the base of an $(\eps_i,N_i)$-Rokhlin tower
(see \S\ref{app:custom.quasi}), where
$\eps_i \goto 0$ and $N_i \goto \oo$.

 {\bf(b)} \  Any measurable subset $\bW\subset\Tor$ can be approximated
arbitrarily well by a disjoint union of tower levels.  That is:
for any $\del>0$, there is some $i\in\Natur$ and some subset
 $\dM\subset\dB(N_i)$ such that, if $\tlbW 
\ =  \ \D\Disj_{\fm\in\dM} \varsigma^\fm(\bJ_i)$, then $\lam(\bW\symdif\tlbW)
\ < \ \del$.\qed
}

\bthmprf[Proof of Proposition \ref{CA.surjective.on.MP}:]
``{\bf(b)} $\IMPLIES$ {\bf(c)}'' 
First note that $\Phitor(\ZP)$ is dense in $\ZP$,
because Proposition \ref{Xi.iso}(a) says $\ZP$ is dense in $\OP$, so 
$\Phitor(\ZP)$ is dense in $\OP$ (and thus, in $\ZP$).
Now apply Proposition \ref{Xi.Phi.iso}.

``{\bf(c)} $\IMPLIES$ {\bf(d)}'' follows from Corollary \ref{dB.dense.implies.dC.dense}

``{\bf(d)} $\IMPLIES$ {\bf(a)}'' \ Corollary \ref{QP.is.dense}
implies  $\Phi(\QS)$ is $d_C$-dense in $\sA^\Lat$; 
thus,  $\Phi(\sA^\Lat)$ is $d_C$-dense in $\sA^\Lat$.  But
$\Phi$ is continuous and $\sA^\Lat$ is $d_C$-compact, so
$\Phi(\sA^\Lat)$ is also $d_C$-compact, thus, $d_C$-closed.  Hence,
$\Phi(\sA^\Lat) = \ \sA^\Lat$.

``{\bf(a)} $\IMPLIES$ {\bf(b)}''\quad
We'll show that $\Phitor(\OP)$ is dense in $\OP$.
It follows that $\Phitor(\MP)$ is dense in $\MP$,
because Proposition \ref{OP.dense.in.MP} 
says $\OP$ is dense in $\MP$.

  Let $\sP\in\OP$ and $\eps>0$. We'll construct 
 $\sQ^{\eps}\in\OP$ such that
$d(\Phitor(\sQ^{\eps}),\sP) \ < \ \eps$.

 Let $\del := \eps/5$.  Let $\{\bJ_i,\del_i,N_i\}_{i=1}^\oo$ be as in Theorem
\ref{delJQuasiRankThm}.  For any $i\in\Natur$ and any $\fb\in\dB(N_i)$, let
$\bJ_i^\fb \ := \ \trans{\fb}(\bJ_i)$. 
Then let $\bS_i := int(\Tor \setminus \D \Disj_{\fb\in\dB(N_i)} \bJ^\fb_i)$.
Fix $s\in\sA$.

\claim{\label{CA.surjective.on.MP.C1}
 There is $i\in\Natur$ and a word $\tile\in\sA^{\dB(N_i)}$ such that,
if $\tlsP$ is the partition obtained by painting 
city $\lb\{\bJ_i^{\fb}\rb\}_{\fb\in\dB(N_i)}$ with
$\bw$ and spacer $s$, then $d_\symdif(\sP, \ \tlsP) \ < \ \del$}.
\bclaimprf
Let $A:=\#(\sA)$.  For each $a\in\sA$, use Theorem
\ref{delJQuasiRankThm}{\bf(b)} to find some set $\tlbP_a$ (a union of
levels in  $\lb\{\bJ_i^{\fb}\rb\}_{\fb\in\dB(N_i)}$) such that
$\lam(\tlbP_a\symdif\bP_a) < \frac{\del}{2A}$.  Assume $\{\tlbP_a\}_{a\in\sA}$
are disjoint.
Enlarge $\tlbP_s$ by adjoining $\bS_i$ to it.  If $i$ is 
large enough, then $\lam(\bS_i)<\frac{\del}{2}$.  Thus, 
\[
d_\symdif(\sP,\tlsP) \ \  = \ \   \D \lam(\bP_s\symdif\tlbP_s) +
\sum_{s\neq a\in\sA} \lam(\bP_a\symdif\tlbP_a)
\ \  \leq \ \   \D \frac{\del}{2} +\frac{\del}{2A} +
\sum_{s\neq a\in\sA} \frac{\del}{2A} \ \  = \ \  \del.
\]
  We define $\tile$:  for any $a\in\sA$ and $\fb\in\dB(N_i)$,  
 $\statement{$\til_\fb  =  a$}
\iff \statement{$\bJ_i^\fb \subset \tlbP_a$}$.  Then $\tlsP$
results from painting $\lb\{\bJ_i^{\fb}\rb\}_{\fb\in\dB(N_i)}$ with
$\tile$ and spacer $s$, as in eqn.(\ref{painting.the.city.eqn}) of
\S\ref{app:custom.quasi}. 
\eclaimprf
  Fix $\tort\in\tlTor\intsct\bJ_i$ and let $\tlbp := \tl\Splat(\tort)$.
Thus, $\tlbp\restr{\dB(N_i)} \ = \ \tile$.
Now, $\Phi$ is surjective on $\sA^\Lat$, so find $\tlbq\in\sA^\Lat$ so
that $\Phi(\tlbq)=\tlbp$.  Suppose $\Phi$ has local map
$\phi:\sA^{\dB(n)}\into\sA$, and let $N_i' := N_i-n$.  If $\stile \ = \
\tlbq\restr{\dB(N_i)}$; then $\Phi(\stile) \ = \ \Phi(\tlbq)\restr{\dB(N'_i)}
\ = \ \tlbp\restr{\dB(N'_i)} \ = \ \tile\restr{\dB(N'_i)}$.
Build a partition $\sQ^{\eps} = \{\bQ^{\eps}_a\}_{a\in\sA}$ by
painting $\lb\{\bJ_i^{\fb}\rb\}_{\fb\in\dB(N_i)}$ with $\stile$,
and then let $\bq^{\eps} := \sQ^{\eps}_{\trans{}}(\tort)$.  Let $\bp^\eps  :=  \Phi(\bq^\eps)$.

\claim{\label{CA.surjective.on.MP.C2.2.e2}
If $i$ is made large enough, then $d_B(\bp^{\eps},\tlbp) < 2\del$.}
\bclaimprf
Let $\dJ_i
= \set{\ell\in\dL}{\trans{\ell}(0)\in\bJ_i}$.
For all $\fj\in\dJ_i$,
 Lemma \ref{painting.the.city.2} says
 $\tlbp\restr{\fj+\dB(N_i)} \ = \ \tile$ and
$\bq^\eps\restr{\fj+\dB(N_i)}  = \ \stile$.  Thus,
$\bp^{\eps}\restr{\fj+\dB(N'_i)}  = \ \phi(\stile) 
 = \   \tile\restr{\dB(N'_i)}
 = \ \tlbp\restr{\fj+\dB(N'_i)}$.  Thus,
\beqn
\label{CA.surjective.on.MP.C2.2.e1}
\bp^{\eps}\restr{\dJ_i+\dB(N'_i)} \quad = \quad \tlbp\restr{\dJ_i+\dB(N'_i)}.
\eeqn
\begin{eqnarray}
\nonumber
\mbox{Hence,} \ 
d_B(\bp^{\eps},\tlbp) &=&
\density{\ell\in\dL \ ; \  p^{\eps}_\ell \neq \tlp_\ell}
\ \ \leeeq{(e\ref{CA.surjective.on.MP.C2.2.e1})} \ \
\density{\maketall \Lat \setminus (\dJ_i+\dB(N'_i))}
\\ \nonumber &=&
 1 - \density{\maketall \dJ_i+\dB(N'_i)}
\quad=\quad
 1 - (2N'_i)^D \cdot \density{\dJ_i}
\\ \nonumber & \eeequals{(\dagger)} &
 1 -  \lb(\frac{N'_i}{N_i}\rb)^D (2N_i)^D \cdot \lam(\bJ_i)
\quad \lt{(*)}\quad 
 1 -  \lb(\frac{N'_i}{N_i}\rb)^D \cdot  (1-\del_i)
\end{eqnarray}
{\bf(e\ref{CA.surjective.on.MP.C2.2.e1})} is by
eqn.(\ref{CA.surjective.on.MP.C2.2.e1}); \ $(\dagger)$ is the Generalized 
Ergodic Theorem, and $(*)$ is because
$\bJ_i$ is the base of a $(\del_i,N_i)$-Rokhlin tower.

Make $i$ large enough that  $\del_i<\del$, and also 
$N_i \ > \D \frac{n}{1-\sqrt[D]{1-\del}}$, so that
$\lb(\frac{N'_i}{N_i}\rb)^D \ > \ (1-\del)$.   Hence
$\D
 1 -  \lb(\frac{N'_i}{N_i}\rb)^D \cdot  (1-\del_i)
\quad < \quad 1 - (1-\del)^2 \ = \ \del(2-\del) < 2\del$.
\eclaimprf
If $\sP^\eps := \Phitor(\sQ^{\eps})$,  then
\begin{eqnarray}
\label{CA.surjective.on.MP.C2.1}
\sP^\eps_{\trans{}}(\tort)
 & \eeequals{(\diamond)} & \Phi\lb(\sQ^\eps_{\trans{}}(\tort)\rb) \quad =\quad
 \Phi(\bq^\eps) \quad =\quad \bp^\eps. \\
\label{CA.surjective.on.MP.C2.2}
\mbox{Thus,}\quad 
d_\symdif(\sP^{\eps},\tlsP)
& \eeequals{(*)} &
 2 d_B\lb(\sP^{\eps}_{\trans{}}(\tort), \ \tl\Splat(\tort)\rb)
\quad \eeequals{(\dagger)} \quad
 2 d_B\lb(\bp^{\eps},\tlbp\rb)  \quad 
\lt{(C\ref{CA.surjective.on.MP.C2.2.e2})} \quad 4\del. \quad\qquad \\
\nonumber
\mbox{Thus,}\quad
d_\symdif(\sP^\eps,\sP)
& \leq & d_\symdif(\sP^\eps,\tlsP) + d_\symdif(\tlsP,\sP)
\quad \lt{(\ddagger)} \quad 4\del + \del \quad=\quad 5\del\quad <\quad \eps.
\end{eqnarray}
Here, $(\diamond)$ is 
Lemma \ref{induced.map.on.partitions}(a);\quad
 $(*)$ is by Proposition  \ref{Xi.iso}(b); \quad
$(\dagger)$ is by eqn.(\ref{CA.surjective.on.MP.C2.1});\quad
{\bf(C\ref{CA.surjective.on.MP.C2.2.e2})} is Claim
\ref{CA.surjective.on.MP.C2.2.e2};\quad
and $(\ddagger)$ is by
eqn.(\ref{CA.surjective.on.MP.C2.2}) and Claim \ref{CA.surjective.on.MP.C1}.

Thus, $d_\symdif(\Phitor(\sQ^{\eps}),\sP) < \eps$.
But $\eps$ is arbitrary.  So $\Phitor(\OP)$ is dense in $\OP$.
\ethmprf
{\em Remarks.}  \ (a) \ If $\Phi:\sA^\Zahl\into\sA^\Zahl$ is surjective,
and $\Phitor(\MP)$ is a {\em closed} subset of $\MP$,
then Theorem \ref{CA.surjective.on.MP}(b) implies that $\Phitor$ is
surjective.  When (if ever) is $\Phitor(\MP)$ closed?  

(b)  Extending Theorem \ref{delJQuasiRankThm} to $\Lat=\Zahl^D$ would
immediately extend Theorem \ref{CA.surjective.on.MP} to  $\Lat=\Zahl^D$.

\section{\label{S:fixedpoints}
Fixed points and Periodic Solutions}

  If $\Phi:\sA^\Lat\into\sA^\Lat$ is a cellular automaton, let
$\Fix{\Phi}:=\set{\ba\in\sA^\Lat}{\Phi(\ba)=\ba}$ be the set of all
{\dfn fixed points} of $\Phi$.  If $p\in\Natur$, then a {\dfn
$p$-periodic} point for $\Phi$ is an element of $\Fix{\Phi^p}$.   If
$\fv\in\Lat$, then a {\dfn $p$-periodic travelling wave} with {\dfn
velocity} $\fv$ is an element of $\Fix{\Phi^p\circ\shift{-\fv}}$.
When does $\Phi$ have quasisturmian fixed/periodic points?  
First, it is easy to verify:

\Proposition{}
{
  $\Fix{\Phi}$, \  $\Fix{\Phi^p}$
and $\Fix{\Phi^P\circ\shift{-\fv}}$ are subshifts of finite type {\rm(SFTs)}.\qed
}

  (See \cite{Kitchens,LindMarcus} for an introduction to subshifts of finite type).    Thus, we ask:  if $\gF\subset\sA^\Lat$
is an SFT, when is $\gF\intsct\QS$ nonempty?
  Let $\gF\subset\sA^\Lat$ be an SFT, and let $a\in\sA$.  
 We say $a$ is an {\dfn inert state} for $\gF$ if
$\bc\in\gF$, where $\bc$ is the constant sequence such that $c_\ell=a$ for all
$\ell\in\Lat$.

\example{ Let $\sA=\{0,1\}$.

$\inn{\rm a}$ \ Let $\Nh = \CC{-1..1}^D$, 
so $\#(\Nh)=3^D$.  The {\dfn  Voter Model} \cite{Vichniac} has local map:
\[
  \phi(\ba) \quad:=\quad
\choice{1 &\If& \gK(\ba) \ > \  {3^D}/{2}; \\
  0 &\If&  \gK(\ba) \ < \ {3^D}/{2}.}
\qquad
\mbox{where} \  \gK(\ba) \ := \  \sum_{\nh\in\Nh} a_b
\]
 Both $0$ and $1$ are inert states for $\Fix{\Phi}$.

$\inn{\rm b}$ \  Let  $\Nh \subseteq\Zahl^2$, and
let $0\leq s_0\leq b_0 \leq b_1 \leq s_1\leq \#(\Nh)$.
A {\dfn Larger than Life} (LtL) CA \cite{Evans1,Evans2,Evans3}
has local map
\[
  \phi(\ba) \quad:=\quad
\choice{1 &&\mbox{if} \  a_0 = 1 \AND  s_0\leq \gK(\ba) \leq s_1; \\
        1 &&\mbox{if} \  a_0 = 0 \AND  b_0\leq \gK(\ba) \leq b_1; \\
        0 &&\mbox{otherwise.}} \qquad
\qquad
\mbox{where} \  \gK(\ba) \ := \  \sum_{\nh\in\Nh} a_b.
\]
For example, J.H.Conway's {\em Game of Life} \cite{BerlekampConwayGuy}
 is an LtL CA, with
$\Nh=\CC{-1..1}^2$,  \ $s_0=b_0=b_1=3$, and $s_1=4$.
LtL CA have many fixed points, periodic points, and travelling waves
with various velocities (for example, the well-known {\em gliders} and
{\em fish} of {\em Life}).  The corresponding SFTs all have $0$ as
an inert state.
}

   Let $V>0$. An SFT $\gF$
with inert state $0$ is {\dfn Dirichlet} with {\dfn valence}
$V$ if, for any $r>V$, and any $\gF$-admissible configuration
$\ba\in\sA^{\dB(r)}$, there exists $\bb\in\gF$ such that
$\bb\restr{\dB(r-V)}  =  \ba\restr{\dB(r-V)}$
but $b_\ell  =  0$ for all $\ell\in\Lat\setminus \dB(r+V)$.
We call $\bb$ a {\dfn Dirichlet extension} of $\ba$.

\Proposition{}
{
 Suppose $\gF$ is a Dirichlet SFT.  If $\trans{}$ is any ergodic
$\Lat$-action on $\Tor$, then $\gF\intsct \QS \ \neq \ \emptyset$
and contains nonconstant elements.
}
\bthmprf
  Let $\tort\in\Tor$.
  We must construct an open partition $\sP\in\OP$ such that
$\Splat(\tort) \ \in \ \gF$.  We will do this by painting
a Rokhlin city (see \S\ref{app:custom.quasi}).
Let $r>0$ and let $V$ be the Dirichlet valence of $\gF$. 
 Let $N:=r+V$, and let 
$\lb\{\trans{\fb}(\bU)\rb\}_{\fb\in\dB(N)}$
be a Rokhlin city in $\Tor$.  Let $\ba\in\sA^{\dB(r)}$ be some 
$\gF$-admissible sequence, and let $\bb\in\sA^\Lat$ be a
Dirichlet extension of $\ba$.  Let $\tile :=  \bb\restr{\dB(N)}$,
let $\eps>0$;  then Corollary \ref{painting.the.city} yields
an open partition
$\sP\in\OP$ so that every element of $\Qshift{\sP}$ is $\eps$-tiled by 
$\tile$ with spacer $0$.  It follows that every 
element of $\Qshift{\sP}$ is $\gF$-admissible.
\ethmprf

 For example, the fixed point SFT of the Voter Model is Dirichlet, as are
the SFTs of fixed points, periodic points, and
travelling waves for any LtL CA.
Thus,  $\QS$ contains nonconstant fixed points for the Voter Model, and
nonconstant travelling waves for  LtL CA.

\section{\label{app:LCA} 
Background on Linear Cellular Automata\protect\footnotemark}

\footnotetext{This section contains technical results which are used
in \S\ref{S:expansive}, \S\ref{S:nilpotence.rigid}, and  \S\ref{S:recur}.}

 Suppose that $p$ is prime and $\sA=\Zahlmod{p}=\CO{0..p}$ is a cyclic group;
then $\sA^\Lat$ is also an abelian group under componentwise addition.
We say $\Phi$ is a {\dfn linear cellular automaton} (LCA) if $\Phi$ has
a local map $\phi:\sA^\Nh\into\sA$ of the form
\beqn
\label{LCA2}
  \phi(\ba) \quad=\quad \sum_{\nh\in\Nh} \varphi_\nh a_\nh \pmod{p},
\qquad\mbox{for any $\ba\in\sA^\Nh$},
\eeqn
where $\Nh\subset\Lat$ is finite, and $\varphi_\nh\in\CO{1..p}$
are constants for all $\nh\in\Nh$.  Equivalently, 
\beqn
\label{LCA}
  \Phi(\ba) \quad=\quad \sum_{\nh\in\Nh} \varphi_\nh\cdot\shift{\nh}(\ba),
\qquad\mbox{for any $\ba\in\sA^\Lat$}.
\eeqn
  If $p=2$, then $\sA=\Zahlmod{2}$,
and $\Phi$ is called a {\dfn boolean linear cellular automaton} (BLCA),
and eqns. (\ref{LCA2}) and (\ref{LCA}) become:
\beqn
\label{BLCA}
  \phi(\ba) \ = \ \sum_{\nh\in\Nh} a_\nh,
\ \mbox{for any $\ba\in\sA^\Nh$},
\qquad\AND\qquad
 \Phi(\ba) \ = \ \sum_{\nh\in\Nh} \shift{\nh}(\ba),
\ \mbox{for any $\ba\in\sA^\Lat$}.
\eeqn
  We call $\Nh$ the {\dfn neighbourhood} of $\Phi$. A BLCA is entirely
determined by its neighbourhood.

  The advantage of the `polynomial of shifts' notation in equations
(\ref{LCA}) and (\ref{BLCA}) is that iteration of $\Phi$ corresponds
to multiplying the polynomial by itself.  For example, if 
$\Phi \ = \ 1+\shift{}$, then the Binomial Theorem says
$\D \Phi^n \ = \ \sum_{n=0}^N \lb({N\atop n}\rb) \shift{n}$.

 For any $N\in\Natur$, let
$\mtrx{N^{(i)}}{i=0}{\oo}$ be the {\dfn $p$-ary expansion} of $N$,
so that $\D N \ =\ \sum_{i=0}^\oo N^{(i)} p^i$.  Let $\Lucas{N}  := 
\set{n\in\CC{0..N}}{n^{(i)} \leq N^{(i)}, \mbox{for} \ \forall \ i \in
\Natur}$.
To get binomial coefficients mod $p$, we use:

\paragraph{\sc Lucas' Theorem \cite{Lucas}:}
{\sl
 $\D \lb[N \atop n\rb]_p \ = \
\prod_{i=0}^\oo
  \lb[N^{(i)} \atop n^{(i)}\rb]_p$, \ 
where we define $\D \lb[N^{(i)} \atop n^{(i)}\rb]_p:=0$ if $n^{(i)}>N^{(i)}$,
and $\D\lb[0 \atop 0 \rb]_p:=1$.
\qquad Thus, $\D\lb[ N \atop n \rb]_p \not= \ 0$ iff $n \in \Lucas{N}$. \qed
}
\breath

  For example, if $p=2$ and $\Phi = 1+\shift{}$, then
Lucas' Theorem says that
$\D \Phi^N   =   \sum_{n\in\Lucas{N}} \shift{n}$.  More generally,
Lucas' Theorem and Fermat's Theorem \cite[\S6, Thm.1]{Dudley} together imply:

\Lemma{\label{p.power.of.LCA}}
{
 Let $p$ be prime, let $\sA=\Zahlmod{p}$, and let
$\Phi \ =\ \D \sum_{\nh\in\Nh} \varphi_\nh\cdot\shift{\nh}$, as
in eqn.{\rm(\ref{LCA})}. Then for any $m\in\Natur$, 
if $P:=p^m$, then
$\Phi^{P} \ =\ \D \sum_{\nh\in\Nh} \varphi_\nh\cdot\shift{P\cdot\nh}$.
\qed
}

\section{\label{app:torus}
Background on Torus Rotation Systems\protect\footnotemark}

\footnotetext{This section contains technical results which are used
in \S\ref{S:expansive} and \S\ref{S:nilpotence.rigid}}

Let $\tau:\Lat\into\Tor$ be a group monomorphism with dense image, and for any
$\ell\in\Lat$, let $\trans{\ell} = \rot{\tau(\ell)}$ denote the
corresponding rotation of $\Tor$.  This defines an ergodic torus rotation
system $(\Tor,d,\varsigma)$.   Torus
rotation systems are {\em minimal} in the sense that every
$\tort\in\Tor$ has dense $\trans{}$-orbit in $\Tor$.  For the
`generic' torus rotation system, an even stronger property holds.
 Suppose $\Lat=\Zahl$.  
 If $p\in\Natur$;
then the system $(\Tor,d,\varsigma)$ is
{\dfn minimal along powers of $p$} if, for any $\tort\in\Tor$,
and any $\eps>0$, there is some $m\in\Natur$ such that
\ $\trans{(p^m)}(0) \ \closeto{\eps} \  \tort$.

  More generally, suppose $\Lat=\Zahl^D$.  For each $d\in\CC{1..D}$, let 
$\fe_d \ := \ (\overbrace{0,\ldots,0}^{d-1},1,\ldots,0)\in\Lat$.
If $p\in\Natur$, then the system $(\Tor,d,\varsigma)$ is
{\dfn minimal along powers of $p$} if, for any $\tort_1,\ldots,\tort_D\in\Tor$,
and any $\eps>0$, there is some $m\in\Natur$ such that, for all $d\in\CC{1..D}$,
\quad $\trans{(p^m\cdot \fe_d)}(0) \ \closeto{\eps} \ \tort_d$.

  We'll show that the `generic' torus rotation is minimal along powers of $p$.
 First, note that for any monomorphism $\tau:\Lat\into\Tor$,
there are unique
$\tora_1,\tora_2,\ldots,\tora_D\in\Tor$ such that $\tau(\ell) \ = \
\ell_1\tora_1 + \ell_2\tora_2+\cdots+ \ell_D\tora_D$ for any
$\ell=(\ell_1,\ldots,\ell_D)\in\Lat$.  Treat $(\tora_1,\ldots,\tora_D)$
as an element of the Cartesian power $(\Tor)^D$,
 and let $\lam^D$ be the product Lebesgue measure on $(\Tor)^D$.

\newcommand{\digit}[4]{\presub{#2} {#1}_{#3}^{#4}}

\Proposition{\label{minimal.along.powers.of.p}}
{
Fix $p\in\Natur$.  For 
$\forall_{\lam^D} \ (\tora_1,\ldots,\tora_D) \ \in \ (\Tor)^D$ the system
 $(\Tor,d,\varsigma)$ is minimal along powers of $p$.
}
\bthmprf
Let $\Tor = \Torus{K} \cong \CO{0,1}^K$.
For each $d\in\CC{1..D}$, let $\tora_d \ =\ 
(\presub{1} a_d,\presub{2} a_d,\ldots,\, _K a_d)$, and
 let $\tort_d \ =\ 
(\presub{1} t_d,\presub{2} t_d,\ldots,\, _K t_d)$.
For each $k\in\CC{1..K}$, suppose
$\presub{k} a_d$ and $\presub{k} t_d$ have $p$-ary expansions:
\[
\presub{k}a_d \quad=\quad \sum_{i=1}^\oo \frac{\digit{a}{k}{d}{i}}{p^i},
\qquad
\AND
\qquad 
\presub{k}t_d \quad=\quad \sum_{i=1}^\oo \frac{\digit{t}{k}{d}{i}}{p^i},
\]
where $\digit{a}{k}{d}{i}\in\CO{0..p}$ and
 $\digit{t}{k}{d}{i}\in\CO{0..p}$ for all $i\in\{1,2,3,\ldots\}$.
 Let $L := \lb\lfloor\log_p(\eps/K)\rb\rfloor$. 

  Let $\sA=\CO{0..p}$.  Fix $d\in\CC{1..D}$ and $k\in\CC{1..K}$.
  For $\forall_\lam \ \tora_d\in\Tor$, the $p$-ary sequence
$\presub{k} \ba_d \ := \ \lb(\digit{a}{k}{d}{1},\digit{a}{k}{d}{2},\digit{a}{k}{d}{3},\ldots\rb)$ is
generic for the $(\frac{1}{p},\ldots,\frac{1}{p})$-Bernoulli measure
on $\sA^\Natur$.   Thus, the word 
$\presub{k}\bt_d  :=  \lb[\digit{t}{k}{d}{1},\digit{t}{k}{d}{2},\ldots,\digit{t}{k}{d}{L}\rb]$
occurs in $\presub{k} \ba_d$ with frequency $p^{-L}$.  Furthermore, for $\forall_{\lam^D} \ 
(\tora_1,\ldots,\tora_D)\in(\Tor)^D$ the collection of sequences
$\lb\{\presub{k} \ba_d \rb\}_{k=1}^{K}\, _{d=1}^D  \ \subset \ \sA^\Natur$
are independent.  Hence, the words $\presub{k}\bt_d$ occur {\em simultaneously}
in  $\presub{k} \ba_d$ for all $d\in\CC{1..D}$ and $k\in\CC{1..K}$
with frequency $p^{-L\cdot K\cdot D} \ > \ 0$.
Thus, there is some $m\in\Natur$ such that
\ $ \lb[\digit{a}{k}{d}{m+1},\digit{a}{k}{d}{m+2},\ldots,\digit{a}{k}{d}{m+L}\rb]
\ = \ \lb[\digit{t}{k}{d}{1},\digit{t}{k}{d}{2},\ldots,\digit{t}{k}{d}{L}\rb]$, \ 
for all $k\in\CC{1..K}$ and $d\in\CC{1..D}$.
Thus, $p^m\cdot  \presub{k}a_d \ \closeto{1/p^L} \  \presub{k}t_d$,
for all $k\in\CC{1..K}$ and $d\in\CC{1..D}$.  Thus, 
 for all  $d\in\CC{1..D}$, \
$\trans{(p^m \fe_d)}(0)  =  
p^m  \tora_d \ \closeto{ K/p^L} \  \tort_d$,
and  $\ K/p^L  \leq  K \frac{\eps}{K}  =  \eps$.
\nolinebreak\ethmprf

\section{Expansiveness
\label{S:expansive}}

  Let $\Phi:\sA^\Lat\into\sA^\Lat$ be a cellular automaton, and let
 $\xi>0$.  If $\bp\in\QS$, we say that the topological
 dynamical system $(\QS,d_B,\Phi)$ is (positively) {\dfn $\xi$-expansive}
 at $\bp$ if, for any $\bq\in\QS$ with $\bq\neq\bp$, there is some $n\in\Natur$ such that
 $d_B\lb(\Phi^n(\bp),\Phi^n(\bq)\rb)>\xi$.  We say $(\QS,d_B,\Phi)$ is
 (positively) {\dfn $\xi$-expansive}  if $(\QS,d_B,\Phi)$ is $\xi$-expansive at
 every $\bp\in\QS$.

 Proposition 13 of \cite{BlanchardFormentiKurka} states that a
cellular automaton $\Phi$  is never expansive in the Besicovitch
topology.  This is proved by constructing a configuration $\ba\in\sA^\Lat$
that is `nonexpansive' for $\Phi$.  However, $\ba\not\in\QS$,
so the proof  of Proposition  13 in \cite{BlanchardFormentiKurka} does {\em not}
apply to $(\QS,d_B,\Phi)$.

Expansiveness is the `opposite' of equicontinuity.  Torus rotation
systems are equicontinuous; \ hence, any shift map acts
equicontinuously on $\QS$.  It is natural to conjecture that {\em all} CA
act equicontinuously on $\QS$, especially in light of Proposition 13
in \cite{BlanchardFormentiKurka}.  We'll refute this by
showing that the boolean linear CA of Example \ref{X:ledrappier}
is expansive on $\QS$.

\Proposition{\label{ledrappier.expansive.on.quasisturm}}
{
   Fix $\tora\in\Tor$, and suppose $\Lat=\Zahl$ acts on $\Tor$ by 
$\trans{z}(\tort) \ = \ \tort+ z\cdot \tora$ for any $\tort\in\Tor$ and
$z\in\Zahl$.   Let $\sA=\Zahlmod{2}=\{0,1\}$, and let $\Phi = 1+\shift{}$
be the {\rm BLCA} of Example {\rm\ref{X:ledrappier}}.
  For $\forall_\lam \ \tora\in\Tor$, and for any $\xi<1$, the
system $(\QS,d_B,\Phi)$ is $\xi$-expansive.
}
\bthmprf  Let $\bo = (\ldots,0,0,0,\ldots)\in\sA^\Zahl$. 
Then $\Phi(\bo)=\bo$, because $\Phi$ is a linear CA.  
We will first show that $(\QS,d_B,\Phi)$ is $\xi$-expansive at $\bo$.

Let $\bp\in\QS$, so that $\bp = \Splat(\tort)$ for some
$\sP\in\OP$ and $\tort\in\tlTor$.  Let $\sO$ be the constant zero
partition (ie. $\sO=\{\bO_0,\bO_1\}$, where $\bO_0 := \Tor$ and $\bO_1
:= \emptyset$).  Thus, $\bo=\sO_{\trans{}}(\tort)$, and
Proposition \ref{Xi.iso}(a) and Lemma \ref{induced.map.on.partitions}(a)
imply that $d_B\lb(\Phi^n(\bp), \bo\rb) \ = \ 
\frac{1}{2} d_\symdif\lb(\Phitor^n(\sP), \ \sO\rb)$.
Thus, we seek $n\in\Natur$ such that
$d_\symdif\lb(\Phitor^n(\sP), \ \sO\rb)  \ > \ 2\xi$.

\Claim{\label{ledrappier.expansive.on.quasisturm.C2}
Let $\bP\subset\Tor$ be measurable, with $0<\lam[\bP]<1$.
There is some $\tort\in\Tor$ such that
 $\lam\lb[\bP \intsct \rot{\tort}(\bP)\maketall \rb] \ < \ \lam[\bP]^2$.}
\bclaimprf
First note that \quad
$\D \int_{\Tor}
 \lam\lb[\bP \intsct \rot{\tort}(\bP)\maketall \rb] \ d\lam[\tort]
\ \ = \ \
\int_{\Tor} \int_{\bP}
 \chr{\rot{\tort}(\bP)}(\tors) \ d\lam[\tors] \ d\lam[\tort]$
\begin{eqnarray}\nonumber
&\eeequals{(*)}& \nonumber
 \int_{\bP} \int_{\Tor}  
\chr{\rot{-\tors}(\bP)}(-\tort) \ d\lam[\tort] \ d\lam[\tors]  
\quad\eeequals{(H)}\quad
 \int_{\bP}  \lam\lb[\rot{-\tors}(\bP)\rb] \ d\lam[\tors]  
\\&\eeequals{(H)}&
 \int_{\bP}  \lam\lb[\bP\rb] \ d\lam[\tors]  
\quad=\quad
 \lam\lb[\bP\rb] \cdot  \lam\lb[\bP\rb]
\quad = \quad \lam[\bP]^2.
\label{ledrappier.expansive.on.quasisturm.e1}
\end{eqnarray}
$(*)$  is because
$\chr{\rot{\tort}(\bP)}(\tors)
\ = \  \chr{\bP}(\tors-\tort)
\  = \  \chr{\rot{-\tors}(\bP)}(-\tort)$ for any
$\tort,\tors\in\Tor$. \
{\bf(H)} is because $\lam$ is the Haar measure on $\Tor$.

  But  $\lam\lb[\bP \intsct \rot{0}(\bP)\maketall \rb] 
 \ = \ \lam[\bP] > \lam[\bP]^2$.  Thus,
eqn.(\ref{ledrappier.expansive.on.quasisturm.e1}) implies
there must be some $\tort\neq 0$ such that 
 $\lam\lb[\bP \intsct \rot{\tort}(\bP)\maketall \rb] \ < \ \lam[\bP]^2$.
\eclaimprf

\Claim{\label{ledrappier.expansive.on.quasisturm.C3}
Let $\bP\subset\Tor$ be measurable, with $0<\lam[\bP]<1$.
For $\forall_\lam\, \tora\in \Tor$, there is some $m\in\Natur$ such
that $\lam\lb[\maketall \bP \symdif \trans{2^m}(\bP)\rb] \ > \ 
\lb(\maketall 2-\lam[\bP]\rb)\cdot \lam[\bP]$.}
\bclaimprf
  This follows from Claim \ref{ledrappier.expansive.on.quasisturm.C2} and
Proposition \ref{minimal.along.powers.of.p}.
\eclaimprf

\Claim{\label{ledrappier.expansive.on.quasisturm.C4}
There is an increasing sequence  $\{m_k\}_{k=1}^\oo \subset \Natur$,
such that, if $\sP^0  :=  \sP$, and
for all $k\in\Natur$, \ $\sP^k  :=  \Phitor^{(2^{m_k})}(\sP^{k-1})$,
then for all $k\in\Natur$, \ $\lam[\bP^{k+1}_1] \ > \ 
\lb(\maketall 2-\lam[\bP^k_1]\rb)\cdot \lam[\bP^k_1]$. }
\bclaimprf
Suppose $\sP=\{\bP_0,\bP_1\}$.
Let $m_1$ be the
result of setting $\bP = \bP_1$ in
  Claim \ref{ledrappier.expansive.on.quasisturm.C3}.
Lemma \ref{p.power.of.LCA} says that
 $\Phi^{(2^m)}  =  1 + \shift{(2^m)}$, so
$\Phitor^{(2^m)}  =  1 + \trans{(2^m)}$.  Thus, if $\sP^1
 \ := \ \Phitor^{(2^m)}(\sP)$, then
 $\bP^1_1 \ = \ \bP_1 \symdif \trans{2^m}(\bP_1)$; hence
Claim \ref{ledrappier.expansive.on.quasisturm.C3} says
$\lam[\bP_1^1] \ > \ 
\lb(\maketall 2-\lam[\bP_1]\rb)\cdot \lam[\bP_1]$.  

  Inductively, suppose we have  $m_k\in\Natur$
and partition $\sP^k  :=  \Phitor^{(2^{m_k})}(\sP^{k-1})$.  
Let $m_{k+1}$ be the result of
setting $\bP = \bP^{k}_1$ in Claim 
\ref{ledrappier.expansive.on.quasisturm.C3}.
Then $\Phi^{(2^{m_{k+1}})}  =  1 + \shift{(2^{m_{k+1}})}$, so
$\Phitor^{(2^{m_{k+1}})}  =  1 + \trans{(2^{m_{k+1}})}$.  Thus, if $\sP^{k+1}
 \ := \ \Phitor^{(2^{m_{k+1}})}(\sP)$, then
 $\bP^{k+1}_1 \ = \ \bP^k_1 \symdif \trans{2^{m_{k+1}}}(\bP^k_1)$; hence
Claim \ref{ledrappier.expansive.on.quasisturm.C3} says
$\lam[\bP^{k+1}_1] \ > \ 
\lb(\maketall 2-\lam[\bP^{k}_1]\rb)\cdot \lam[\bP^{k}_1]$.  
\eclaimprf

\Claim{\label{ledrappier.expansive.on.quasisturm.C5}
If $\{r_k\}_{k=1}^\oo  \subset  \OO{0,1}$, and $r_{k+1} > (2 - r_k)\cdot r_k$
for all $k\in\Natur$, then $\D \lim_{k\goto\oo}\,  r_k  =  1$.}
\bclaimprf
The sequence is increasing
(because  $r_k<1$, so $(2-r_k)  >   (2-1)  =  1$, so
 $r_{k+1}  >  (2-r_k)\cdot r_k  >   r_k$, \ for all $ \ k\in\Natur$).
  We claim that  $\D\sup_{k\in\Natur}\,  r_k  =  1$.
Suppose not; then  $\exists \ y<1$ such that
$r_k<y$, \ for all $ \ k\in\Natur$.  But then  $r_{k+1}  =  (2-r_{k})\cdot r_k
  >  (2-y) \cdot r_k$, \ 
for all $\, k$; hence, $r_k \ > \ (2-y)^k \cdot r_1$.  Thus,
$\D \lim_{k\goto\oo}\, r_k \ \geq \  r_1\cdot \lim_{k\goto\oo}\, (2-y)^k 
 \ = \ \oo$. Contradiction.
\eclaimprf

Now, for any $k\in\Natur$, \quad
$d_\symdif(\sO,\sP^k)
\ = \  2\cdot \lam[\bO_1\symdif \bP^k_1]
\ = \    2\cdot\lam[\bP^k_1]$.  Thus,
combining Claims \ref{ledrappier.expansive.on.quasisturm.C4} and
\ref{ledrappier.expansive.on.quasisturm.C5}, we conclude that
$\D  \lim_{k\goto\oo}\, d_\symdif(\sO,\sP^k)
\ = \ 
 2\cdot \lim_{k\goto\oo} \, \lam[\bP^k_1] \quad = \quad 2$.
  But observe that 
\beq
\sP^k &=&
\Phitor^{(2^{m_k})}(\sP_{k-1})
\quad = \quad \cdots \quad = \quad
  \Phitor^{(2^{m_k})} \circ  \Phitor^{(2^{m_{k-1}})} \circ
\cdots \circ  \Phitor^{(2^{m_1})} (\sP)
\\&=&
\Phitor^{n_k}(\sP),
\qquad \mbox{where $n_k \ = \ 2^{m_k}+ 2^{m_{k-1}} + \cdots + 2^{m_1}$.}
\eeq
Hence, $\D \lim_{k\goto\oo} \, d_B\lb(\bo, \Phi^{n_k}(\bp)\rb)
\ \eeequals{(*)} \ \frac{1}{2} \lim_{k\goto\oo} \, d_\symdif(\sO,\sP^k)
\ = \ 1 \ > \ \xi$, where $(*)$ is by
Proposition  \ref{Xi.iso}(b) and Lemma \ref{induced.map.on.partitions}(a).

  It follows that $(\Phi,\QS,d_B)$ is $\xi$-expansive at $\bo$.
To see that  $(\Phi,\QS,d_B)$  is
$\xi$-expansive everywhere, let $\bp,\bq\in\QS$;
  we must find some $n\in\Natur$ such that
 $d_B\lb(\Phi^n(\bp),\Phi^n(\bq)\rb)>\xi$.  Let $\br=\bp-\bq$;
find $n\in\Natur$ such that
 $d_B\lb(\Phi^n(\br),\ \bo\rb) \ = \
 d_B\lb(\Phi^n(\br),\Phi^n(\bo)\rb) \ > \ \xi$.
Thus, $d_B\lb(\Phi^n(\bp),\Phi^n(\bq)\rb)
\ = \ d_B\lb(\Phi^n(\bp)-\Phi^n(\bq),\bo\rb)
\ \eeequals{(L)} \ d_B\lb(\Phi^n(\bp-\bq), \ \bo\rb)
\ = \ d_B\lb(\Phi^n(\br), \ \bo\rb) \ > \ \xi$, where {\bf(L)} is because
$\Phi$ is linear.
\ethmprf

\section{\label{S:nilpotence.rigid}
Niltropism and Rigidity}

  Let $\{n_k\}_{k=1}^\oo\subseteq\Natur$.  The system $(\QS,d_B,\Phi)$
  is {\dfn rigid} along $\{n_k\}_{k=1}^\oo$ if, for all
  $\bq\in\QS$,\quad $\D \dBlim_{k\goto\oo}\, \Phi^{n_k}(\bq) = \bq$.
  Suppose $\sA$ is an abelian group with identity element $0$.  Let
  $\bo\in\sA^\Lat$ be the constant zero configuration (ie. $o_\ell=0$,
  for all $ \ell\in\Lat$).  Then $(\QS,d_B,\Phi)$ is {\dfn niltropic}
  along $\{n_k\}_{k=1}^\oo$ if, for all $\bq\in\QS$, \quad $\D
  \dBlim_{k\goto\oo}\, \Phi^{n_k}(\bq) = \bo$.

  Likewise, $(\MP,d_\symdif,\Phitor)$ is {\dfn rigid} along 
 $\{n_k\}_{k=1}^\oo$ if, for all $\sQ\in\MP$,\quad $\D
\symlim_{k\goto\oo}\, \Phitor^{n_k}(\sQ)  =  \sQ$, and
$(\MP,d_\symdif,\Phitor)$ is {\dfn niltropic} along  $\{n_k\}_{k=1}^\oo$ if for 
all $\sQ\in\MP$,
\quad $\D \symlim_{k\goto\oo}\, \Phitor^{n_k}(\sQ)  =  \sO$,
where $\sO\in\MP$ is the trivial
partition (ie. $\bO_0:=\Tor$ 
and $\bO_a:=\emptyset$ if $0\neq a\in\sA$).

\example{\label{X:shifts.are.rigid}
Let $\Lat=\Zahl$, and let $\tora\in\Tor$ be such that
$\trans{z}(\tort) = \tort+z\cdot \tora$ for all $z\in\Zahl$ and $\tort\in\Tor$.  Then $(\QS,d_B,\shift{})$ is rigid.  To see this, find a sequence
 $\{n_k\}_{k=1}^\oo\in\Natur$ such that $\D\lim_{k\goto\oo}\, n_k\cdot\tora
\ = \ 0$ in $\Tor$.  Thus, if $\tort\in\Tor$, then
\beqn
\label{X:shifts.are.rigid.eqn1}
\lim_{k\goto\oo}\, \trans{n_k}(\tort)\quad=\quad
\lim_{k\goto\oo}\, n_k\tora +\tort
\quad=\quad
\tort.
\eeqn
Now, for any $\bp\in\QS$, there is some $\sP\in\OP$
and some $\tort\in\tlTor$ such that $\bp=\Splat(\tort)$;
 hence,
$\D\lim_{k\goto\oo}\, \shift{n_k}\lb(\bp\rb) \ = \ 
\lim_{k\goto\oo}\, \shift{n_k}\lb(\Splat(\tort)\rb) 
\  \eeequals{(*)} \
\lim_{k\goto\oo}\, \Splat\lb(\trans{n_k}(\tort)\rb) \ \eeequals{(\dagger)} \ 
\Splat(\tort) \ = \ \bp$.
Here, $(*)$ is by Proposition \ref{splat.homeo}(a), and $(\dagger)$ is by
Proposition \ref{splat.homeo}(b) and
eqn.(\ref{X:shifts.are.rigid.eqn1}).
}

 If $\sA=\Zahlmod{p}$ and $\Phi$ is the LCA
of eqn.(\ref{LCA}) in \S\ref{app:LCA}, we
define $\trace{\Phi}  :=  \sum_{\nh\in\Nh} \varphi_\nh$ (mod $p$).
Suppose $\Lat=\Zahl^D$.  Let 
$\tau:\Lat\into\Tor$
be a monomorphism such that $\trans{\ell}(\tort) \ = \ \rot{\tau(\ell)}(\tort)$
for all $\ell\in\Lat$ and $\tort\in\Tor$.   Then there are unique
$\tora_1,\tora_2,\ldots,\tora_D\in\Tor$ such that $\tau(\ell) \ = \
\ell_1\tora_1 + \ell_2\tora_2+\cdots+ \ell_D\tora_D$ for any
$\ell=(\ell_1,\ldots,\ell_D)\in\Lat$.

\Theorem{\label{LCA.nilpotence.rigid}}
{   
   Let  $\sA=\Zahlmod{p}$ ($p$ prime).
  For $\forall_\lam \ \tora_1,\ldots,\tora_D\in\Tor$, 
there is a sequence $\{m_j\}_{j=1}^\oo\in\Natur$ such that,
 if $n_j:=p^{m_j}$ for all $j\in\Natur$, then  for any linear cellular
automaton $\Phi$,

  {\bf(a)} \ If $\trace{\Phi} \equiv 0$, then $(\MP,d_\symdif,\Phitor)$ and
$(\QS,d_B,\Phi)$ are niltropic along
$\{n_k\}_{k=1}^\oo$.

  {\bf(b)} \  If $\trace{\Phi}\not\equiv 0$, then $(\MP,d_\symdif,\Phitor)$ and
$(\QS,d_B,\Phi)$ are rigid along 
$\{(p-1)\cdot n_k\}_{k=1}^\oo$.
 }
\bthmprf We'll show rigidity/niltropism for $(\MP,d_\symdif,\Phitor)$;
 rigidity/niltropism for $(\QS,d_B,\Phi)$ follows from
 Proposition \ref{Xi.Phi.iso}.
 Let $\{\eps_j\}_{j=1}^\oo$ be a sequence decreasing to zero.
For  all $j>0$,
Proposition \ref{minimal.along.powers.of.p} yields some  $m_j\in\Natur$
such that $d(p^{m_j}\cdot\tora_d, \  0) \ < \ \eps_j$ for 
 all $d\in\CC{1..D}$.
Let $n_j:=p^{m_j}$.
 Let $\Phi$ be the LCA of eqn.(\ref{LCA}) in \S\ref{app:LCA}, and
 let $M \ := \ \D \max_{\nh\in\Nh} |\nh|$.
Let $\sP\in\MP$ and fix $\del>0$.

\Claim{\label{LCA.nilpotence.rigid.C2}
There exists $J\in\Natur$ such that, if $j>J$, then,
 for all $ \ \nh\in\Nh$,
\ \ $d_\symdif\lb(\maketall \shift{(n_j\cdot\nh)}(\sP), \ \sP \rb)  < 
 \del$.
}
\bclaimprf  
  The function $\Tor\ni\tort\mapsto\rot{\tort}(\sP)\in\MP$ is 
$d_\symdif$-continuous \cite[Prop 8.5, p.229]{Folland}.
  Thus, find some $\eps_J>0$ such that, 
\beqn
\label{LCA.nilpotence.rigid.C2.e1}
\mbox{For any $\tort\in\Tor$,} \quad
\statement{$d(\tort,0) \ < \ M\cdot \eps_J$}\IMPLIES
\statement{$d_\symdif\lb(\rot{\tort}(\sP), \  \sP\rb) \ < \ \del$}.
\eeqn
If $\fb=(b_1,\ldots,b_D)$, then 
$\trans{(n_j\cdot\nh)} \ = \ \rot{(b_1 n_j\tora_1+\cdots+ b_D n_j \tora_D)}$.
If $j>J$, then 
$d\lb[(n_j\cdot \tora_d), \  0\rb] \ < \ \eps_j$, for  all $d\in\CC{1..D}$. 
Hence, $d\lb[ (b_1 n_j\tora_1+\cdots+ b_D n_j \tora_D), \ 0  \rb ] \ < \ 
(|b_1|+\cdots+|b_D|)\cdot \eps_j \ = \ |\nh|\cdot\eps_j \ < \ M\eps_J$,
so $d_\symdif\lb(\maketall \trans{(n_j\cdot\nh)}(\sP), \ \sP\rb)
\ = \ d_\symdif\lb(\maketall \rot{(b_1 n_j\tora_1+\cdots+ b_D n_j \tora_D)}(\sP), \ \sP\rb) \ \lt{(\ref{LCA.nilpotence.rigid.C2.e1})} \ \del$.
\eclaimprf
Let $\Phi$ be as in eqn.(\ref{LCA}), and let $R:=\trace{\Phi}$.
Let $B:=\#(\Nh)$.  If $j>J$, then
\[
\Phitor^{n_j}(\sP)  \quad\eeequals{(\dagger)}\quad
 \sum_{\nh\in\Nh} \varphi_\nh\cdot \trans{(n_j\nh)}(\sP)
\quad\stackrel{(*)}{\closeto{B\del}} \quad  \sum_{\nh\in\Nh} 
\varphi_\nh\cdot \sP
\quad=\quad
\choice{\sO &\If & R=0; \\
    R\cdot \sP &\If & R\neq 0.}
\]
 $(\dagger)$ is by Lemma \ref{p.power.of.LCA}, and $(*)$ is
Claim \ref{LCA.nilpotence.rigid.C2}.

  This works for any $\del>0$.
  Thus, if $R=0$, then 
$\D  \lim_{j\goto\oo} \ d_\symdif\lb(\maketall\Phitor^{n_j}(\sP), \ \sO\rb)
 \ = \ 0$.

If $R\neq 0$, then 
$\D \Phitor^{(p-1)\cdot n_j}(\sP) 
\ \closeto{\del_1} \   
R^{p-1}\cdot \sP
\ \eeequals{(*)} \  \sP$, where $\del_1 \ := \ (p-1) B\cdot \del$, and
 $(*)$ is by Fermat's Theorem \cite[\S6, Thm.1]{Dudley}.
 Thus,
$\D  \lim_{j\goto\oo} \ d_\symdif\lb(\maketall\Phitor^{(p-1)\cdot n_j}(\sP), \ \sP\rb)
  =  0$.\nolinebreak
\ethmprf

\Examples{\label{X:rigid.nilpotent} {\rm Let $\Lat=\Zahl$.}}
{\item \label{X:rigid.nilpotent.Z2}
If $\sA=\Zahlmod{2}$,  then $1+\shift{}$
is niltropic, but $1+\shift{} + \shift{2}$ is rigid.

\item \label{X:rigid.nilpotent.Z3}
If $\sA=\Zahlmod{3}$, then $1+\shift{}$
is rigid, while $1+\shift{} + \shift{2}$ and  $1+2\shift{}$ are niltropic.
}

\section{\label{S:inv.quasi.measure} CA-Invariant Quasisturmian Measures}

  If $\Phi\colon\sA^\Lat\into\sA^\Lat$ is a CA, are there any
$\Phi$-invariant quasisturmian measures on $\sA^\Lat$?  The answer is
reminiscent of J. King's Weak Closure Theorem \cite{KingWkClos}.

\Theorem{\label{meas.phi.inv.iff.orbit.phitor.inv}}
{
  Let $\Phi:\sA^\Lat\into\sA^\Lat$ be a CA.
If $\sP\in\MP$ and $\mu = \QM{\sP}$, then

$\statement{$\mu$ is $\Phi$-invariant}
\iff
\statement{$\Phitor(\sP) \ = \ \rot{\tort}(\sP)$ for some $\tort\in\Tor$}$.
}

  To prove this, we must first characterize 
the set $\QSM$ of $\trans{}$-quasisturmian measures.

\Proposition{\label{QS.measure.iff.QP}}
{
  Let $\mu\in\sM(\sA^\Lat)$.  Then
$ \statement{$\mu$ is quasisturmian}$

$\ \iff \ 
\statement{The MP$\Lat$S $(\sA^\Lat,\shift{},\mu)$ is isomorphic to
a torus rotation}$.
}

  To prove Proposition \ref{QS.measure.iff.QP}, in turn, we use the
  following lemma:

\Lemma{\label{splat.quasiperiodic.factor.map}}
{
 Let $\sP\in\MP$ and  let $\trans{}$ be an $\Lat$-action on $\Tor$.  Let 
$\mu = \QM{\sP}$.  

  {\bf(a)} \ $\Splat:\Tor\into\sA^\Lat$ is an 
{\rm(MP$\Lat$S)} epimorphism from $(\Tor,\trans{},\lam)$ to $(\sA^\Lat,\shift{},\mu)$.

  {\bf(b)} \  If $\sP$ is simple {\rm(see \S\ref{app:simple.inject})}, then
$\Splat$ is an isomorphism from $(\Tor,\trans{},\lam)$ to $(\sA^\Lat,\shift{},\mu)$.
}
\bthmprf
{\bf(a)}\quad
$\Splat:\Tor\into\sA^\Lat$ is a measurable function, and
Proposition \ref{splat.homeo}(a) says
$\shift{\ell}\circ\Splat
\ = \ \Splat\circ\trans{\ell}$ for any $\ell\in\Lat$. Hence $\Splat$ is
a epimorphism from $(\Tor,\trans{},\lam)$ to $(\sA^\Lat,\shift{},\mu)$.

{\bf(b)}\quad  If $\sP$ is {\em simple} then
Lemma \ref{simple.partn.inj}
 says that $\Splat$ is injective \lamae.
Thus, $\Splat$ is an isomorphism  from $(\Tor,\trans{},\lam)$ to $(\sA^\Lat,\shift{},\mu)$.\ethmprf

\bthmprf[Proof of Proposition {\rm\ref{QS.measure.iff.QP}:}]
 `$\IMPLIES$' \  Let $\sP\in\MP$, let $\trans{}$ be an $\Lat$-action on $\Tor$, and let $\mu=\QM{\sP}$.  If
 $\sP$ is {simple} then Lemma
 \ref{splat.quasiperiodic.factor.map}(b) says $\Splat$ is an MP$\Lat$S
 isomorphism from  $(\Tor,\trans{},\lam)$ to
 $(\sA^\Lat,\shift{},\mu)$.  If $\sP$ is {not} simple, then 
Lemma \ref{quotient.partition}(e)
yields a simple 
 partition $\quosP\in\quoMP$ with $\mu = \QM[\quotrans{}]{\quosP}$; then
 Lemma \ref{splat.quasiperiodic.factor.map}(b) says  $\quoSplat$
 is an isomorphism from $(\quoTor,\quotrans{},\quolam)$ to
 $(\sA^\Lat,\shift{},\mu)$.

`$\seilpmi$'\quad  Suppose $\Psi:\Tor\into\sA^\Lat$ is an isomorphism
from $(\Tor,\trans{},\lam)$ to $(\sA^\Lat,\shift{},\mu)$.  We claim
that $\Psi = \ \Splat$ for some measurable partition $\sP\in\MP$.
For any $\ell\in\Lat$, let $\pr{\ell}:\sA^\Lat\into\sA$ be projection
 onto the $\ell$th
coordinate ---ie. $\pr{\ell}(\ba) \ := \ a_\ell$, for any $\ba\in\sA^\Lat$.
Define $\sP \ := \ \pr{0}\circ\Psi:\Tor\into\sA$.  Then $\sP$ is a measurable
$\sA$-valued function ---ie. an $\sA$-labelled partition ---on $\Tor$.
Observe that, for $\forall_\lam \ \tort\in\Tor$ and all $\ell\in\dL$,
\[
  \Splat(\tort)_\ell \ \ = \ \ 
\sP\circ\trans{\ell}(\tort)
\ \ = \ \
\pr{0}\circ\Psi\circ\trans{\ell}(\tort)
\ \ = \ \
\pr{0}\circ\shift{\ell}\circ\Psi(\tort)
\ \ = \ \
\pr{\ell}\circ\Psi(\tort)
\ \ = \ \
\Psi(\tort)_\ell.
\]
  Hence $\Splat(\tort) \ = \ \Psi(\tort)$.  This holds for 
 $\forall_\lam \ \tort\in\Tor$, so $\Splat = \Psi$ \lamae.
Hence, $\mu \ = \ \Psi(\lam) \ = \ \Splat(\lam) \ = \ \QM{\sP}$,
so $\mu$ is quasisturmian.
\ethmprf

\bthmprf[Proof of Theorem \ref{meas.phi.inv.iff.orbit.phitor.inv}:]
``$\seilpmi$'' is obvious. We must prove ``$\IMPLIES$''.

{\bf Case 1:} ({\em $\sP$ is simple})  
If $\sQ=\Phitor(\sP)$, then $\QM{\sQ}  =  
\QM{\Phitor(\sP)} \ \eeequals{(*)} \ \Phi(\QM{\sP})  =  \Phi(\mu) 
\ \eeequals{(\dagger)} \  \mu$,
where $(*)$ is Lemma \ref{induced.map.on.partitions}(c),
and $(\dagger)$ is because $\mu$ is $\Phi$-invariant.
Thus it  suffices to show:

\Claim{If $\sQ\in\MP$ and $\QM{\sQ} \ = \ \mu$, then
 $\sQ = \rot{\tort}(\sP)$ for some $\tort\in\Tor$.} 
\bclaimprf
   $\sP$ is simple, so Lemma \ref{splat.quasiperiodic.factor.map}(b)
says $\Splat$ is an MP$\Lat$S isomorphism from $(\Tor,\trans{},\lam)$
to $(\sA^\Lat,\shift{},\mu)$.
 $\sQ$ may not be simple, but
  Lemma \ref{splat.quasiperiodic.factor.map}(a)
says that $\sQ_{\trans{}}:\Tor\into\sA^\Lat$ is an MP$\Lat$S  epimorphism.
Thus, $\psi \ := \ \sP_{\trans{}}^{-1}\circ\sQ_{\trans{}}:\Tor\into\Tor$
is a endomorphism from  $(\Tor,\trans{},\lam)$ to itself,
(ie. $\psi:\Tor\into\Tor$ and $\psi\circ\trans{\ell}
\ = \ \trans{\ell}\circ\psi$ for all $\ell\in\Lat$).

\subclaim{ There is some $\tora\in\Tor$ such that $\psi=\rot{\tora}$.}
\bsubclaimprf
  Define $\alp:\Tor\into\Tor$ by $\alp(\tort) := \psi(\tort)-\tort$.
We claim that $\alp$ is constant \lamae.

Why? $\alp$ is measurable because $\psi$ is measurable.
 $\alp$ is $\trans{}$-invariant, because for
any $\ell\in\Lat$ and $\tort\in\Tor$,\quad
$\alp\lb(\trans{\ell}(\tort)\rb)
\ = \ \psi\lb(\trans{\ell}(\tort)\rb) - \trans{\ell}(\tort)
\ = \ \trans{\ell}\lb(\psi(\tort)\rb) - \trans{\ell}(\tort)
\ = \ \lb(\psi(\tort)+\tau(\ell)\rb) - \lb(\tort+\tau(\ell)\rb)
\ = \  \psi(\tort)-\tort \ = \ \alp(\tort)$.  But $\trans{}$ is ergodic,
so $\alp$ must be constant \lamae.

Hence, $\psi=\rot{\tora}$,
where $\tora$ is the constant value of $\alp$.  
\esubclaimprf
But then \ 
$
\statement{$\sP_{\trans{}}^{-1}\circ\sQ_{\trans{}} \ = \ \rot{\tora}$}
\iff
\statement{$\sQ_{\trans{}} \ = \ \sP_{\trans{}}\circ\rot{\tora}$}
\iff
\statement{$\sQ \ = \ \sP\circ\rot{\tora}$, \ \lamae}$. \
In other words, $\sQ = \rot{\tora}(\sP)$.\eclaimprf

\paragraph{Case 2:} ({\em $\sP$ is not simple})
 Lemma \ref{quotient.partition}(e) yields
a simple partition $\quosP$ such that $\mu=\quoQM{\quosP}$.
 When applied to $\quosP$,  {\bf Case 1} implies that
 $\quoPhitor(\quosP) \ = \ \rot{\quotort}(\quosP)$
for some $\quotort\in\quoTor$.  Let $q:\Tor\into\quoTor$ be
the quotient map, and $\tort\in q^{-1}(\quotort)$; it follows that
$\Phitor(\sP) \ = \ \rot{\tort}(\sP)$.
\ethmprf

\example{Suppose $\bq\in\QS$ is a quasisturmian travelling wave for $\Phi$,
with period 1 (see \S\ref{S:fixedpoints}). 
 Thus, $\bq=\Sqlat(\tort)$, for some $\sQ\in\OP$ and $\tort\in\tlTor$.
 Hence, $\mu:=\QM{\sQ}$ is a $\Phi$-invariant
QS measure.  If $\bq$ has velocity  $\fv\in\Lat$, then
 $\Phitor(\sQ) \ = \ \trans{\fv}(\sQ)$.}

Let $\Lat=\Zahl^D$, and let $\tora_1,\ldots,\tora_D\in\Tor$
be as defined prior to Theorem \ref{LCA.nilpotence.rigid}
 in \S\ref{S:nilpotence.rigid}.

\Proposition{\label{no.LCA.inv.quasi.measures}}
{ 
  Let $\sA=\Zahlmod{p}$ ($p$ prime); \ let $\Phi$ be a
linear CA with $\trace{\Phi}= 0$ {\rm (see \S\ref{S:nilpotence.rigid})}.
Then for $\forall_\lam \ 
\tora_1,\ldots,\tora_D\in\Tor$, there are no nontrivial $\Phi$-invariant
measures in $\QSM$.
}
\bthmprf
  Suppose $\mu=\QM{\sP}$ for some $\sP\in\MP$. By Lemma
 \ref{splat.quasiperiodic.factor.map}(c), we can assume $\sP$ is simple.
 If $\mu$ is $\Phi$-invariant,
then Corollary \ref{meas.phi.inv.iff.orbit.phitor.inv} says that
$\Phitor(\sP) \ = \ \rot{\tort}(\sP)$ for some $\tort\in\Tor$.
But $\trace{\Phi}\equiv 0$ (mod $p$), so 
Theorem \ref{LCA.nilpotence.rigid}(a) yields
a sequence $\{n_j\}_{j=1}^\oo\in\Natur$ such that
\beqn
\label{no.LCA.inv.quasi.measures.e1}
  \lim_{j\goto\oo} \rot{n_j \tort}(\sP)
\quad=\quad
  \lim_{j\goto\oo} \Phitor^{n_j}(\sP)
\quad=\quad
\sO,
\eeqn
  where $\sO$ is the trivial partition.
  However $\rot{\tort}$ acts isometrically
on $\MP$, so for any $j\in\Natur$
\beqn
\label{no.LCA.inv.quasi.measures.e2}
  d\lb(\rot{n_j \tort}(\sP), \ \sO\rb)
\quad=\quad
  d\lb(\sP, \ \rot{-n_j \tort}(\sO)\rb)
\quad=\quad
  d(\sP, \ \sO).
\eeqn
Combining equations (\ref{no.LCA.inv.quasi.measures.e1}) 
and (\ref{no.LCA.inv.quasi.measures.e2}) we conclude that
$d(\sP, \ \sO)  =  0$.  Thus $\sP=\sO$.\nolinebreak
\ethmprf
\example{\label{X:no.LCA.inv.quasi.measures}
 Let $\Lat:=\Zahl$, and recall Example \ref{X:rigid.nilpotent}.

 If $\sA=\Zahlmod{2}$, then 
Example \ref{X:rigid.nilpotent.Z2} implies that
$1+\shift{}$ has no QS invariant measures.

 If  $\sA=\Zahlmod{3}$, then Example \ref{X:rigid.nilpotent.Z3} 
implies that $1+2\shift{}$ has no QS invariant measures.\nolinebreak
}

\section{Asymptotic Nonrandomization
\label{S:recur}}

 Let $\sA=\Zahlmod{2}$, and let $\eta$ be the
$(\frac{1}{2},\frac{1}{2})$ Bernoulli measure on $\sA^\Zahl$.  If
$\mu\in\sM(\sA^\Zahl)$ and $\Phi$ is a CA, then $\Phi$ {\dfn
asymptotically randomizes} $\mu$ if there is a set $\dJ\subset\Natur$
of density 1 such that $\D\wkstlim_{\dJ\ni j\goto\oo} \Phi^j(\mu) \ = \
\eta$.  Linear cellular automata
asymptotically randomize a wide range of measures, including most
Bernoulli measures \cite{MaassMartinez,Lind} and Markov chains
\cite{PivatoYassawi1,FerMaassMartNey}, and also Markov random fields
 supported on the full shift, subshifts of finite type, or sofic shifts
\cite{PivatoYassawi2,PivatoYassawi3}. 
 Indeed, the current literature has no examples of nonperiodic measures
on $\Zahlmod{2}^\Zahl$ 
that are {\em not} asymptotically randomized by LCA.  We will now show that
a broad class of quasisturmian measures are not asymptotically randomized by
the linear CA $\Phi=1+\shift{}$.

A measure $\mu\in\sM(\sA^\Zahl)$ is {\dfn dyadically recurrent} if there
is a sequence $n_k\goesto{k\goto\oo}{} \oo$ and constants
$\eps,\del>0$ such that, if $k\in\Natur$, then   $\mu[\gR_k] > \eps$, 
where we define
\beqn
\label{dyadic.recurrence.set}
  \gR_k \quad:=\quad \set{\ba\in\sA^\Zahl}{a_\ell = a_{\ell+2^{n_k}} \ \mbox{for all} \ \ell \in \CO{0...\del 2^{n_k}}},
\eeqn
 Thus, if
$N_k := 2^{n_k}$, then
$\ba_{\CO{0...\del N_k}}  =  \ba_{\CO{N_k...N_k+\del N_k}}$ for 
all $\ba\in\sA^\Zahl$ in a subset of measure $\eps$.

\Proposition{\label{DR.no.AR}}
{
 If $\mu$ is dyadically recurrent, and $\Phi:=1+\shift{}$, then $\Phi$
cannot asymptotically randomize $\mu$.
}
\bthmprf
Let $E:=\lb\lceil-\log_2(\eps)\rb\rceil +1$, such that 
$2^{-E}\leq\frac{\eps}{2}$.
  Find $K\in\Natur$ such that,
if $k\geq K$, then $\del 2^{n_k-1}>E$.
Let $
\dJ \ := \ \set{2^{n_k}+j}{k\geq K \AND j\in \CO{0...\del 2^{n_k-1}}}$.
Then 
\[
 \density{\dJ} \ \  \geq \ \  
\lim_{k\goto\oo} \frac{\#\CO{2^{n_k}...2^{n_k}\!+\!\del 2^{n_k-1}}}
{\#\CC{0...2^{n_k+1}}} \ \  = \ \  
\lim_{k\goto\oo} \frac{\del 2^{n_k-1}}{2^{n_k+1}} \ \  = \ \  \frac{\del}{4} \ \  > \ \   0.
\]
If $J=(2^{n_k}+j)\in\dJ$, then Lucas' Theorem (\S\ref{app:LCA}) implies that
$\D\Phi^J 
\ = \ \sum_{\ell\in\Lucas{j}} \shift{\ell} \ + \ \sum_{\ell\in\Lucas{j}} \shift{\ell+ 2^{n_k}}$. \   Thus, if $\ba\in\gR_k$, then, for all $e\in\CO{0..E}$,
\beq
\Phi^J(\ba)_e & = &
 \sum_{\ell\in\Lucas{j}} a_{e+\ell} \ + \ \sum_{\ell\in\Lucas{j}} a_{e+\ell+ 2^{n_k}}
\quad = \quad 
\sum_{\ell\in\Lucas{j}} \lb(\maketall a_{e+\ell} \ + \ a_{e+\ell+2^{n_k}}\rb)
\\ &\eeequals{(*)} &  2 \cdot \sum_{\ell\in\Lucas{j}} a_{e+\ell}
\quad = \quad 0.
\eeq
$(*)$: if $\ell\in\CC{0..j}$, then
$e+\ell  <  E + j  <  2\cdot \del 2^{n_k-1}
 =  \del 2^{n_k}$, so $a_{e+\ell} = a_{e+\ell+2^{n_k}}$, by
eqn.(\ref{dyadic.recurrence.set}).

Thus, if $\gO:=\set{\bb\in\sA^\Zahl}{b_{\CO{0..E}}=[0..0]}$, then
$\Phi^J(\ba)\in\gO$ for all $\ba\in\gR_k$.  Hence
$\Phi^J(\mu)[\gO] \ \geq \ \mu[\gR_k] > \eps > \eps/2 \geq 2^{-E} = \eta[\gO]$.
Thus, $\Phi^j(\mu)[\gO]$ can't converge to $\eta[\gO]$ along elements in
$\dJ$.  But $\density{\dJ}  >  0$;  thus,
$\Phi^j(\mu)$ can't weak* converge to $\eta$ along a set of density 1.
\ethmprf
 Now we'll construct a dyadically recurrent quasisturmian measure.
  If $\Lat=\Zahl$, then a $\Zahl$-action $\trans{}$ is {\dfn dyadically recurrent} if there is a sequence $n_k\goesto{k\goto\oo}{}\oo$ such that,
for all $\tort\in\Tor$, \ 
$d\lb(\trans{2^{n_k}}(\tort), \tort\rb) \ < \ 2^{-n_k}$.
Recall the definition of the Lipschitz pseudomeasure $\setsize[L]{\bullet}$
from \S\ref{S:boundary}.

\Proposition{\label{DRact.2.DRmeas}}
{
 Let $\sP\in\OP$, with $\setsize[L]{\partial\sP} \ < \ \oo$.
 If $\trans{}$ is dyadically recurrent, then
$\mu=\QM{\sP}$ is dyadically recurrent.
}
\bthmprf Fix $\del \ < \  \D \frac{1}{2\setsize[L]{\sP}}$;  then
$\eps' := 2\del\cdot\setsize[L]{\sP} <1$, so that
$\eps:= 1-\eps' > 0$.

If $\bC = \Ball(\partial\sP,2^{-n_k})$, then
$\lam[\bC] \ \leq \ 2^{1-n_k}\cdot\setsize[L]{\sP}$ by eqn.(\ref{set.size.defn.1}) in
\S\ref{S:boundary}.
Observe that, for any $\tort\in\Tor$,
$\statement{$\sP(\tort) \neq \sP\lb(\trans{2^{n_k}}(\tort)\rb)$}
\IMPLIES \statement{$\tort\in\bC$}$.  
Now, for any  $\ell\in\CC{0...2^{n_k-m}}$, define
\beq
\bB_\ell & := & \set{\tort\in\dT}{\sP(\tort)_\ell \neq \sP(\tort)_{2^{n_k}+\ell}}
\quad = \quad 
\set{\tort\in\dT}{\sP\lb(\trans{\ell}(\tort)\rb)
 \neq \sP\lb(\trans{\ell+2^{n_k}}(\tort)\rb)}
\\
&\subset& \set{\tort\in\dT}{\trans{\ell}(\tort) \in \bC} 
\quad = \quad \trans{-\ell}(\bC).
\eeq
Thus, $\lam[\bB_\ell] < \lam[\bC]   \leq 2^{1-n_k}\setsize[L]{\sP}$.
So if 
$\D \bB  :=  \Union_{\ell=0}^{\lfloor \del 2^{n_k}\rfloor } \bB_\ell$,  
 then $\lam[\bB] \ \leq \  \del 2^{n_k} \cdot 2^{1-n_k}\setsize[L]{\sP}
= 2\del \setsize[L]{\sP}\  = \  \eps' < 1$. 
  So if $\bR_k := \dT\setminus\bB$ and
$\gR_k := \Splat(\bR_k)$, then $\mu[\gR_k] = \lam[\bR_k] 
 \ \geq \  1 - \eps' \ = \  \eps$.

 If $\tort\in\bR_k$, then $\sP(\tort)_\ell =
\sP(\tort)_{\ell+2^{n_k}}$, for  all $\ell\in
\CO{0..\del 2^{n_k}}$.  So if $\ba\in\gR_k$, then $a_\ell =
a_{\ell+2^{n_k}}$, for
 all $\ell\in \CO{0..\del 2^{n_k}}$.
\ethmprf
  When is a $\Zahl$-action $\trans{}$ dyadically recurrent?  
 Let $\Tor:=\Torus{1}$;
 then there is some $\torr\in\Tor$ such that
$\trans{z}(\tort) = \tort+z\torr$ for all $\tort\in\Tor$ and $z\in\Zahl$.
Identify $\Torus{1}\cong\CO{0,1}$, and
suppose that $\torr$ has binary expansion $\torr =
0.r_1r_2r_3\ldots$.  Say that $\torr$ is {\dfn dyadically
recurrent} if there is a sequence $n_k\goesto{k\goto\oo}{}\oo$ 
such that $r_j = 0$ for all
$j\in\OC{n_k...2n_k}$.   For example, let $n_1=1$, and
inductively define $n_{k+1} = 2n_k + 1$.  This gives the sequence
$\{1,3,7,15,31,\ldots\}$.  Then the number
\[
 \torr \ = \ \sum_{k=1}^\oo 2^{n_k}
\ = \ 
0.1\underbrace{0}_{1} \ 1 \ 
  \underbrace{000}_{3} \ 1 \ 
 \underbrace{0000000}_{7} \ 1 \ 
  \underbrace{00\ldots0}_{15} \ 1  \ldots
\]
is dyadically recurrent.

\Proposition{\label{DRnum.2.DRact}}
{
  Let $\bR \subset\Torus{1}$ be the set of dyadically recurrent elements.
\bthmlist
  \item  If $\torr\in\bR$, then $\trans{z}(\tort) := \tort+z\torr$
is a dyadically recurrent $\Zahl$-action.

  \item  $\bR$ is a dense $G_\delta$ subset of $\Torus{1}$, but $\lam[\bR]=0$.
\ethmlist
}
\bthmprf
{\bf(a)}\quad For any $k\in\Natur$,  observe that
\[
2^{n_k}\cdot\torr  \ \  = \ \   0.r_{(n_k+1)} \ r_{(n_k+2)}\ldots 
\ \  = \ \  0.\underbrace{000\ldots00}_{n_k} \ r_{(2n_k+1)} \ r_{(2n_k+2)}\ldots
\ \  <  \ \   2^{-n_k}.
\]
  Thus, for any $\tort\in\Tor$,
\quad
$d\lb(\trans{2^{n_k}}(\tort), \tort\rb) \ = \
d\lb(2^{n_k}\torr + \tort, \  \tort\rb) \ = \  |2^{n_k}\cdot\torr| \ < \  2^{-n_k}$.

{\bf(b)} \quad For all
$N\in\Natur$, let $\bI_N := \D \OO{0,\frac{1}{2^{2N}}}$, and let $\D
\bU_N := \Union_{n=0}^{2^N} (\frac{n}{2^N} + \bI_N)$.  Thus, if
$\toru\in\bU_N$ and $\toru:=0.u_1u_2u_3\ldots$, then
$u_n=0$ for $n\in\CO{N..2N}$.
Now $\bU_N$ is open (it's a union of open intervals), and $\bU_N$
is $(\frac{1}{2^N})$-dense in $\Torus{1}$.  Thus, for any
$N\in\Natur$, \quad
$\bB_N \ := \ \D
\Union_{n=N}^\oo \bU_n$ is open and dense in $\Torus{1}$.
But $\bR \ = \ \D \Intsct_{N=1}^\oo \ \Union_{n=N}^\oo \bU_n
\ = \ \Intsct_{N=1}^\oo \bB_N$, so $\bR$ is dense and $G_\del$.
Finally, $\lam[\bR]=0$ by  the Borel-Cantelli Lemma, because
 $\D \sum_{n=1}^\oo \lam[\bU_n] \ = \ \sum_{n=1}^\oo \frac{1}{2^n}  \ = \  1
 \   < \  \oo$.
\ethmprf
{\sc Questions:} (a) \ Proposition \ref{DR.no.AR} says there must be
other weak* cluster points of the set
 $\{\frac{1}{N} \sum_{n=1}^N \Phi^n(\mu)\}_{N=1}^\oo$
besides $\eta$. What are they?  Any cluster point $\mu_\oo$ will be a
$\Phi$-invariant measure, so if $\mu_\oo$ has nonzero entropy then
$\mu_\oo=\eta$ by \cite[Theorem 4.1]{MaassHostMartinez}.  Thus, we know that
$\mu_\oo$ has zero entropy. However, Example \ref{X:no.LCA.inv.quasi.measures} says that $\mu_\oo$ cannot be
quasisturmian.

(b) \ Propositions \ref{DR.no.AR}, \ref{DRact.2.DRmeas} and
\ref{DRnum.2.DRact} together imply asymptotic nonrandomization 
for a class of irrational rotations which is comeager but of measure zero.
Is there {\em any} quasisturmian measure which is asymptotically randomized 
by {\em any} cellular automaton?

\breath
{\em Acknowledgments.} \ 
This research was partially supported by NSERC Canada.  The author would
like to thank the referee for many perceptive and valuable suggestions.

\newpage

{
\bibliographystyle{plain}
\bibliography{bibliography}
}

\end{document}